\documentclass[a4paper,11pt]{amsart}

\usepackage{amscd}
\usepackage{amsmath}
\usepackage{graphicx}
\usepackage{txfonts}
\usepackage{dsfont}
\usepackage{tikz-cd}
\usepackage[titletoc]{appendix}

\usepackage{geometry}
\usepackage{hyperref}
\usepackage{tikz}

\theoremstyle{definition}
\newtheorem{Defn}{Definition}[section]

\theoremstyle{plain}
\newtheorem{Lemma}[Defn]{Lemma}
\newtheorem{prop}[Defn]{Proposition}
\newtheorem{theorem}[Defn]{Theorem}
\newtheorem{Corollary}[Defn]{Corollary}
\newtheorem{main}{Theorem}

\theoremstyle{remark}
\newtheorem{remark}[Defn]{Remark}
\newtheorem{Example}[Defn]{Example}
\newtheorem{notation}[Defn]{Notation}

\DeclareMathOperator{\Rac}{\G_{\mathcal{R}}}
\DeclareMathOperator{\G}{\mathbf{F}}
\DeclareMathOperator{\Aug}{\mathfrak{I}}
\DeclareMathOperator{\Dim}{\mathbf{N}}

\definecolor{dPurple}{rgb}{0.7,0,0.1}

\title{Polynomiality ideals, with applications to polynomial functors over free nilpotent groups}

\author{Minkyu Kim}
\address{KIAS, Seoul, South Korea}
\email{kimminq@kias.re.kr}
\urladdr{https://minq92.github.io/Gaeul-Autumn/}

\date{}

\begin{document}

\maketitle

\begin{abstract}
    Let $\mathds{k}$ be a unital commutative ring.
    For a monoidal category $\mathcal{C}$ whose unit is a zero object, we introduce a descending filtration of the $\mathds{k}$-linearization of $\mathcal{C}$ by two-sided ideals, called the \emph{polynomiality ideals}: the $d$-th ideal encodes polynomiality of degree at most $d$ for functors from $\mathcal{C}$ to $\mathds{k}$-modules.
    These ideals yield general criteria for the existence of polynomial functors of a given degree and for equivalences between the categories of polynomial functors of degree at most $d$ over different base categories.
    For Lawvere theories satisfying a mild condition, we identify the polynomiality ideals with powers of augmentation ideals of monoids; this identification makes the criteria computable by means of classical group ring theory.
    As the main application, we study polynomial functors over finitely generated free nilpotent groups, comparing them across nilpotency classes and polynomial degrees.
    We obtain refinements of results of Baues and Pirashvili in the case of class at most $2$, recover several folklore results for free groups and free abelian groups, and establish a modular analogue formulated with dimension subgroups over a field of positive characteristic in place of the lower central series.
    Finally, we show that, in contrast to polynomial functors, analytic functors admit no analogous ideal in most cases.
\end{abstract}

\section{Introduction}

\subsection{Motivating questions}
\label{202510211607}
For a category $\mathcal{C}$, a $\mathcal{C}$-module is a functor from $\mathcal{C}$ to the category of $\mathds{k}$-modules.
In this paper, $\mathds{k}$ denotes a unital commutative ring unless otherwise specified.
When $\mathcal{C}$ is a monoidal category with zero unit, i.e. the unit for the monoidal structure is a zero object, one has the notion of the polynomial degree, or the degree, of a $\mathcal{C}$-module \cite{hartl2015polynomial}; it serves as an invariant that allows for a systematic study of $\mathcal{C}$-modules.
A polynomial $\mathcal{C}$-module is a $\mathcal{C}$-module with a finite polynomial degree.
The origin of these notions goes back to Eilenberg and Mac Lane \cite{EML} where $\mathcal{C}$ is taken to be an additive category.

Some fundamental questions arise for a given non-negative integer $d \in \mathds{N}$ and $\mathcal{C},\mathcal{D}$ monoidal categories with zero unit:
\begin{enumerate}
    \item[(Q1)] does there exist a $\mathcal{C}$-module of degree equal to $d$?
    \item[(Q2)] for a strong monoidal functor $\xi : \mathcal{C} \to \mathcal{D}$, does the functor $\xi$ induce an equivalence between the category of polynomial $\mathcal{D}$-modules of degree at most $d$ and that of polynomial $\mathcal{C}$-modules of degree at most $d$?
\end{enumerate}
We address these questions by developing a general framework.

\vspace{1mm}
Regarding question (Q1), let $\mathrm{D}(\mathcal{C})$ denote the set of $d \in \mathds{N}$ for which there exists a $\mathcal{C}$-module of degree equal to $d$.
If the ground ring $\mathds{k}$ is the integer ring, then, as a consequence of \cite{PolyPira}, we have $\mathrm{D}(\mathbf{P}_{R}) = \mathds{N}$ where $R$ is a unital ring and $\mathbf{P}_{R}$ denotes the category of finitely generated projective left $R$-modules.
In contrast, there exist additive categories $\mathcal{C}$ for which no non-constant polynomial $\mathcal{C}$-modules exist; equivalently $\mathrm{D}(\mathcal{C}) = \{ 0 \}$.
It was proved in \cite[Proposition 2.13]{DTV2023} that $\mathrm{D}(\mathcal{C}) = \{ 0\}$ if and only if, for any objects $X,Y$ in $\mathcal{C}$, we have $\mathds{k} \otimes_{\mathds{Z}} \mathcal{C}(X,Y) \cong 0$.

Let $\G$ denote the category of finitely generated free groups.
For a positive integer $c \in \mathds{N}^\ast$, let $\mathbf{N}_{c}$ be the category of finitely generated free nilpotent groups of class at most $c$.
In particular, $\mathbf{N}_{1}$ is the category of finitely generated free abelian groups.
For $\mathcal{C} \in \{ \G , \mathbf{N}_1 \}$ or the opposite, it is folklore that, if $\mathds{k} \neq 0$, then $\mathrm{D}(\mathcal{C}) = \mathds{N}$.

\vspace{1mm}
With respect to question (Q2), let $\Gamma(\xi)$ be the set of $d \in \mathds{N}$ for which $\xi$ induces an equivalence of the categories of $\mathcal{C}$-modules and of $\mathcal{D}$-modules, both of polynomial degree at most $d$.

For suitable $\mathcal{C}$, polynomial $\mathcal{C}$-modules of degree at most 2 admit a concise description in terms of square groups \cite{BJ1994,BP1999,HV2011}.
Let $[n]$ denote the set of non-negative integers at most $n$.
By the results of Baues and Pirashvili \cite{BP1999}, one obtains $[2] \subset \Gamma(\G \to \mathbf{N}_2)$ together with $[1] \subset \Gamma(\G \to \mathbf{N}_1)$ and $2 \not\in \Gamma(\G \to \mathbf{N}_1)$ when $\mathds{k} = \mathds{Z}$.
Here, $\G \to \mathbf{N}_c$ denotes the nilpotentization functor.

Let $\mathbf{W}$ be the category of finitely generated free monoids.
Hartl, Pirashvili and Vespa proved in \cite{hartl2015polynomial} that, for any $d \in \mathds{N}$, the category of polynomial $\mathbf{W}$-modules is equivalent to that of polynomial $\G$-modules, both of $\deg \leq d$.
In other words, $\Gamma(\mathbf{W} \to \G) = \mathds{N}$.

\subsection{Main results}
\label{202510011345}

In this section, we present some general criteria, addressing the previous questions (Q1) and (Q2), that lead to the main results.

The first main theorem concerns question (Q1) and provides a dichotomy for the set $\mathrm{D}(\mathcal{C})$.
Let $\mathcal{C}^{\mathsf{o}}$ denote the opposite category of $\mathcal{C}$.
\begin{main}[see Theorem \ref{202609142120}, Proposition \ref{202512231721}] \label{2026009152033}
    Let $\mathcal{C}$ be a monoidal category with zero unit.
    Then $\mathrm{D}(\mathcal{C})$ is either $\{ 0 \}$ or $\mathds{N}$, and $\mathrm{D}(\mathcal{C}^{\mathsf{o}}) = \mathrm{D}(\mathcal{C})$.
\end{main}

To determine this dichotomy together with (Q2), we introduce a canonical filtration of two-sided ideals that is intrinsic to $\mathcal{C}$, described below.

The notion of a polynomial map from a monoid to an abelian group---one that vanishes on a power of the associated augmentation ideal---can be traced back at least to Passi \cite{MR537126}.
We treat polynomial $\mathcal{C}$-modules of degree at most $d$ by means of analogues of the powers of augmentation ideals.
These analogues are two-sided ideals $\mathtt{I}^{(d)}_{\mathcal{C}}$, $d \in \mathds{N}$, of the monad $\mathtt{L}_{\mathcal{C}}$ naturally associated with $\mathcal{C}$, regarded as a monad in the bicategory of matrices in the sense of \cite{betti1983variation}; we refer to $\mathtt{I}^{(d)}_{\mathcal{C}}$ as the {\it $d$-th polynomiality ideal}, and precise definitions will appear in Section \ref{202408041517}.
The polynomiality ideals yield a canonical filtration:
\begin{align} \label{202609152026}
    \mathtt{I}^{(\bullet)}_{\mathcal{C}} : \quad \cdots \subset \mathtt{I}^{(d+1)}_{\mathcal{C}} \subset \mathtt{I}^{(d)}_{\mathcal{C}} \subset \cdots \subset \mathtt{I}^{(0)}_{\mathcal{C}} \subset  \mathtt{L}_{\mathcal{C}} .
\end{align}
The following theorems show that questions (Q1) and (Q2) reduce to an analysis of the ideals $\mathtt{I}^{(d)}_{\mathcal{C}}$, which is a problem intrinsic to the category $\mathcal{C}$.

\begin{main}[see Section \ref{202404111426}] \label{202510011337}
    Let $\mathcal{C}$ be a small monoidal category with zero unit.
    The following conditions are equivalent:
    \begin{enumerate}
        \item $\mathrm{D}(\mathcal{C}) = \mathds{N}$.
        \item The filtration (\ref{202609152026}) is strictly decreasing, i.e. $\mathtt{I}^{(d+1)}_{\mathcal{C}} \subsetneq \mathtt{I}^{(d)}_{\mathcal{C}}$ for every $d \in \mathds{N}$.
        \item 
        $\mathtt{I}^{(1)}_{\mathcal{C}} \subsetneq \mathtt{I}^{(0)}_{\mathcal{C}}$.
    \end{enumerate}
    Likewise, the following conditions are equivalent:
    \begin{enumerate}
        \item[(1')] $\mathrm{D}(\mathcal{C}) = \{0\}$; in other words, every polynomial $\mathcal{C}$-module is constant.
        \item[(2')] The filtration (\ref{202609152026}) collapses, i.e. $\mathtt{I}^{(d)}_{\mathcal{C}} = \mathtt{I}^{(0)}_{\mathcal{C}}$ for every $d$.
        \item[(3')] 
        $\mathtt{I}^{(1)}_{\mathcal{C}} = \mathtt{I}^{(0)}_{\mathcal{C}}$.
    \end{enumerate}
    By Theorem \ref{2026009152033}, exactly one of these two cases occurs.
\end{main}

\vspace{2mm}

We now turn from a single category to a comparison between categories.
The filtration by polynomiality ideals also controls how polynomial functors behave under base change along a functor.
Let $\mathcal{C},\mathcal{D}$ be small monoidal categories with zero unit.
Let $\xi : \mathcal{C} \to \mathcal{D}$ be a full strong monoidal functor that is the identity on objects.
Let 
$$\mathrm{B}(\xi) {:=} \{  d \in \mathds{N} \mid \xi:  \mathtt{I}^{(d-1)}_{\mathcal{C}}/\mathtt{I}^{(d)}_{\mathcal{C}} \stackrel{\cong}{\longrightarrow}  \mathtt{I}^{(d-1)}_{\mathcal{D}}/\mathtt{I}^{(d)}_{\mathcal{D}} \}.$$
The set $\Gamma(\xi)$ introduced in Section \ref{202510211607} can be investigated by using $\mathrm{B}(\xi)$:
\begin{main}[see Theorem \ref{202509291857}] \label{202510011336}
    The set $\Gamma(\xi)$ coincides with the maximal subset of $\mathrm{B}(\xi)$ consisting of consecutive integers starting at $0$.
\end{main}

In the rest of this paper, we mainly consider categories with finite products and a zero object, which are monoidal categories with zero unit via the cartesian structure.
What finite products provide is a useful description of the polynomiality ideals by generators (Theorem \ref{202609281707}): concretely, the generators are constructed from the diagonal map $X \to X^{\times(d+1)}$.
This is what underlies all the applications of this paper.

\subsection{Polynomial functors on Lawvere theories}

We now give several applications illustrating the main results of this paper.
An algebraic theory in the sense of Lawvere \cite{Lawvere1963}, a {\it Lawvere theory} for short, is a small category with finite products, whose objects are $\mathds{N}$, with the products on objects given by the addition of natural numbers.
In this paper, applying the preceding main theorems, we study Lawvere theories satisfying further conditions, one of which requires that each morphism set $\mathcal{C}_n {:=} \mathcal{C}(n,1)$ carries a natural monoid structure.
Details will be given in Section \ref{202509031755}.
All the base categories appearing here are Lawvere theories satisfying the conditions.
For such a Lawvere theory $\mathcal{C}$, the following structural property holds:
\begin{main}[see Theorem \ref{202408011134} for details]
\label{202509241101}
    For $d \in \mathds{N}$, we have a natural isomorphism
    $$\mathtt{I}^{(d)}_{\mathcal{C}} (m,n) \cong \Aug ( \mathcal{C}_n^{\times m})^{d+1},~ n,m \in \mathds{N},$$
    where $\Aug(M)$ denotes the augmentation ideal of a monoid $M$ over $\mathds{k}$.
\end{main}

This result is well-known for additive categories $\mathcal{C}$.
For $\mathcal{C} = \G^{\mathsf{o}}$, this is equivalent to \cite[Proposition 2.7]{djament2016cohomologie}.
Our result can be applied to numerous examples presented in Sections \ref{202509031755} and \ref{202510211610}.

The following theorems are the first outcome of applying the general theory of Section \ref{202510011345} to Lawvere theories, and the polynomiality ideals play the central role throughout.
The following theorem is concerned with (Q1) for the category of finitely generated free nilpotent groups:
\begin{main} \label{202510171045}
    Suppose that $\mathds{k} \neq 0$.
    For $c \in \mathds{N}^\ast$ and $\mathcal{C} \in \{ \mathbf{N}_{c} , \mathbf{N}_{c}^{\mathsf{o}}\}$, we have $$\mathrm{D}(\mathcal{C}) = \mathds{N}.$$
\end{main}
The proof can be found in Section \ref{202510271047} where we give more general results.

The following theorem addresses question (Q2) for the categories $\mathbf{N}_{c}$.
To state the theorem, we recall that the nilpotentization yields the following sequence of full functors:
\begin{align} \label{202510191757}
    \G \longrightarrow \cdots \longrightarrow \mathbf{N}_{c}  \longrightarrow \cdots \longrightarrow \mathbf{N}_2 \longrightarrow \mathbf{N}_1 .
\end{align}

\begin{main}
\label{202510101544}
    Suppose that $\mathds{k} \neq 0$.
    For $c_0, c_1 \in \mathds{N}^\ast \cup \{\infty\}$ such that $c_0 < c_1$, we have
    $$
    \Gamma( \mathbf{N}_{c_1} \to \mathbf{N}_{c_0}) = \Gamma( \mathbf{N}_{c_1}^{\mathsf{o}} \to \mathbf{N}_{c_0}^{\mathsf{o}}) = [c_0] ,   
    $$
    where we set $\mathbf{N}_{\infty} {:=} \G$.
\end{main}
The proof follows from Theorem \ref{202509021231} and Corollary \ref{202510141158}.
The application to $c_0\in \{1,2\}$ and $c_1=\infty$ recovers the results of Baues and Pirashvili.

\begin{remark}
    Polynomial functor theory has played an important role in the work of Aurélien Djament and Christine Vespa on stable homology of automorphism groups of free groups with coefficients in functors.
    While the present work does not contribute to this line of research, it suggests that the case of nilpotent groups involves features of a different nature, depending on the relationship between the nilpotency level and the polynomial degree.
\end{remark}

\vspace{2mm}
Let $R$ be a unital ring and $\mathbf{M}_{R}$ denote the category of free left $R$-modules of finite rank.
The study of polynomial $\mathbf{M}_{R}$-modules has found important applications in algebraic topology and representation theory, as explained in the lecture series \cite{franjou2003rational}, which focuses on the case where $R$ is a finite field.

In analogy with the sequence (\ref{202510191757}), we consider the following sequence of full functors connecting $\G$ and $\mathbf{M}_{\mathds{F}_{p}}$ where $p$ is a prime and $\mathds{F}_{p}$ is the field of $p$ elements:
$$
\G \longrightarrow \cdots \longrightarrow \Dim_{c,p} \longrightarrow \cdots \longrightarrow \Dim_{2,p} \longrightarrow \Dim_{1,p} = \mathbf{M}_{\mathds{F}_{p}} .
$$
For $c \in \mathds{N}^\ast$ and a group $G$, let $D_{c,p} (G)$ be the $c$-th dimension subgroup over $\mathds{F}_{p}$.
Then $\Dim_{c,p}$ denotes the full subcategory of groups whose objects are the quotients $\mathsf{F}_n / D_{c+1,p}(\mathsf{F}_n)$, $n \in \mathds{N}$, where $\mathsf{F}_n$ is the free group of rank $n$.
In this paper, using the categories $\Dim_{c,p}$, we investigate a modular analogue of the above results.

The relation to the category $\mathbf{N}_c$ is given by the full functor $\mathbf{N}_c \to \Dim_{c,p}$, induced by the fact that the dimension subgroups contain the lower central series.
The following theorem highlights a distinction between polynomial functors over $\mathbf{N}_c$ and $\Dim_{c,p}$:
\begin{main} \label{202510191826}
    Suppose that $\mathds{k}$ is a field of positive characteristic $p$.
    Let $c_0 \in \mathds{N}^\ast$ and $r_0 = \lfloor \log_p (c_0) \rfloor +1$.
    We then have
    \begin{align*}
        [c_0] \subset \Gamma( \mathbf{N}_{c_0} \to \Dim_{c_0,p} ) = \Gamma( \mathbf{N}_{c_0}^{\mathsf{o}} \to \Dim_{c_0,p}^{\mathsf{o}} ) \subset [p^{r_0}-1] .
    \end{align*}
\end{main}
The proof can be found in Section \ref{202512301553}.
As a consequence, for instance, if $d \geq p^{r_0}$, then the functor $\mathbf{N}_{c_0} \to \Dim_{c_0,p}$ does not induce an equivalence of the category of polynomial $\mathbf{N}_{c_0}$-modules of degree at most $d$ and that of polynomial $\Dim_{c_0,p}$-modules of degree at most $d$.

The following theorem gives a modular analogue of Theorem \ref{202510171045}:

\begin{main} \label{202510171046}
    Suppose that $\mathds{k}$ is a field of characteristic $p \geq 0$.
    For a prime $q$ and $c \in \mathds{N}^\ast$, consider $\mathcal{C} \in \{  \Dim_{c,q}, \Dim_{c,q}^{\mathsf{o}} \}$.
    If $p = q$, then $\mathrm{D}(\mathcal{C}) = \mathds{N}$; otherwise, $\mathrm{D}(\mathcal{C}) =\{0\}$.
\end{main}
For the proof, the reader is referred to Section \ref{202510271047}.
We also present an analog of Theorem \ref{202510101544}, which is proved in Section \ref{202512301553}:
\begin{main} \label{202510161105}
    Suppose that $\mathds{k}$ is a field of positive characteristic $p$.
    For $c_0,c_1 \in \mathds{N}^\ast \cup \{ \infty\}$ such that $c_0 < c_1$, we have
    $$
    \Gamma (\Dim_{c_1,p} \to \Dim_{c_0,p}) = \Gamma (\Dim_{c_1,p}^{\mathsf{o}} \to \Dim_{c_0,p}^{\mathsf{o}}) = [c_0],
    $$
    where we regard $\Dim_{\infty,p}$ as $\G$.
\end{main}

\vspace{2mm}
We now explain how the results above fit together.
By Theorem \ref{202510011336}, the (non-)equivalences in Theorems \ref{202510101544}, \ref{202510191826} and \ref{202510161105} reduce to a comparison of the associated polynomiality ideals; likewise, by Theorem \ref{202510011337}, Theorems \ref{202510171045} and \ref{202510171046} reduce to deciding whether $\mathtt{I}^{(1)}_{\mathcal{C}} = \mathtt{I}^{(0)}_{\mathcal{C}}$.
Theorem \ref{202509241101} then identifies these ideals with powers of augmentation ideals of monoids, so that the problem becomes one of classical group ring theory.
At this final step, we rely on classical theorems on augmentation ideals of groups, in particular Quillen's theorem \cite{Quillen1968} and the Sandling--Tahara theorem \cite{sandling1979augmentation}.

This method applies beyond the categories considered above.
In the body of the paper, we refine it to study polynomial functors on the Lawvere theory associated with a radical functor for groups.
In this generality, Theorem \ref{202512231657} refines Theorems \ref{202510171045} and \ref{202510171046}, Theorem \ref{202510281425} generalizes Theorem \ref{202510101544}, and Theorem \ref{202510021437} generalizes Theorems \ref{202510191826} and \ref{202510161105}.

\subsection{Analyticity and further studies}
\label{202410041335}

Let $\mathcal{C}$ be a monoidal category with zero unit.
The results in Section \ref{202510011345} show that polynomial $\mathcal{C}$-modules of degree at most $d$ are governed by the $d$-th polynomiality ideal.
A $\mathcal{C}$-module is {\it analytic} if its filtration by the polynomial approximations converges to itself.
A natural question is whether there exists a two-sided ideal that governs analytic $\mathcal{C}$-modules.
To make this question precise, in Section \ref{202510131439}, we establish a general relationship among two-sided ideals, adjunctions and properties of modules.
Within this framework, we show that, for most Lawvere theories $\mathcal{C}$, no such two-sided ideal can exist, in contrast to the case of polynomial functors.
Further details will be provided in Section \ref{202509301437}.

\vspace{2mm}
This is the first in a series of papers devoted to a systematic study of isomorphism-invariant properties of modules and associated adjunctions.
In various contexts, when one studies modules with certain properties, it is useful to consider adjunctions between module categories that approximate arbitrary modules by modules with that property.
We develop a method that unifies such adjunctions by means of ideals, and illustrate the approach with new examples.
This paper treats two-sided ideals with a focus on polynomial functors.
In the sequels, all the general constructions given here are extended to left and right ideals; the first of these has appeared as \cite{kimprop}.

\vspace{3mm}
This paper is outlined as follows.
In Section \ref{202609211329}, we review polynomial functors and study the set $\mathrm{D}(\mathcal{C})$.
In Section \ref{202410161743}, we introduce the bicategory of matrices $\mathsf{Mat}_{\mathds{k}}$, and give an overview of its monad theory.
In Section \ref{202510131439}, we give a general relationship among two-sided ideals, adjunctions of module categories and properties of modules.
Within this framework, in Section \ref{202408041517}, we introduce the polynomiality ideal which encodes the property of being of finite degree.
Using this, we prove Theorem \ref{202510011336}.
In Section \ref{202510211820}, we study the polynomiality ideal for a Lawvere theory.
In Section \ref{202510071642}, we introduce the notion of a radical functor for groups, and study the associated Lawvere theory.
In Section \ref{202510231558}, we give an overview of the Sandling--Tahara theorem and Quillen's theorem.
In Section \ref{202510211823}, based on our framework, we study polynomial $\mathbf{N}_{c}^{\mathsf{o}}$-modules and prove Theorem \ref{202510101544}.
In Section \ref{202510211825}, we study polynomial $\Dim_{c,p}^{\mathsf{o}}$-modules, and prove Theorems \ref{202510191826} and \ref{202510161105}.
In Section \ref{202509301437}, we discuss the existence of a core $\mathtt{L}_{\mathcal{C}}$-internalizer for analytic $\mathtt{L}_{\mathcal{C}}$-modules.
In Appendix \ref{202510142110}, we give general arguments underlying the constructions of Section \ref{202510131439}.

\section*{Notation}

\begin{itemize}
    \item $\mathds{N}$ and $\mathds{N}^\ast$ respectively denote the set of nonnegative integers and that of positive integers.
    \item $\mathds{Z}$ denotes the ring of integers.
    \item $\mathds{Q}$ denotes the rational field.
    \item $[n]$ denotes the set of non-negative integers at most $n$.
\end{itemize}

\section{Polynomial functors on monoidal categories with zero unit} \label{202609211329}

In this section, we briefly review polynomial functors and fix notation.
We also study the set $\mathrm{D}(\mathcal{C})$, which is related to question (Q1) from the introduction.

\subsection{Recollections on polynomial functor theory}
\label{202509301224}

We give a recollection of polynomial functor theory, as in \cite{hartl2015polynomial}.
A {\it monoidal category with zero unit} is a monoidal category $(\mathcal{C}, \vee, \ast)$, with monoidal product $\vee$ and unit object $\ast$, such that $\ast$ is a zero object.
For an object $X$, denote by $\varepsilon_X : X \to \ast$ and $\eta_X : \ast \to X$ the unique morphisms.
For an $r$-tuple $X_1,\cdots,X_r$ of objects and $1 \leq k \leq r$, we set
$$\epsilon_{X_1,\cdots,X_r;k} {:=} \bigvee^{k-1}_{i=1} \mathrm{id}_{X_i} \vee \varepsilon_{X_k} \vee \bigvee_{i=k+1}^{r} \mathrm{id}_{X_i} , ~\mathrm{and}\quad \eta_{X_1,\cdots,X_r;k} {:=} \bigvee^{k-1}_{i=1} \mathrm{id}_{X_i} \vee \eta_{X_k} \vee \bigvee_{i=k+1}^{r} \mathrm{id}_{X_i}.$$
The composition $\epsilon_{X_1,\cdots,X_r;k} \circ \eta_{X_1,\cdots,X_r;k}$ is the identity, so that
$$
e_{X_1,\cdots,X_r;k} {:=} \eta_{X_1,\cdots,X_r;k} \circ \epsilon_{X_1,\cdots,X_r;k}
$$
is an idempotent on $\bigvee_{i=1}^r X_i$.
We also write $e_{X} {:=} e_{X;1}$.

A functor $F : \mathcal{C} \to \mathds{k}\mbox{-} \mathsf{Mod}$ is {\it polynomial of degree at most $d$}, denoted by $\deg F \leq d$, if its $(d+1)$-th cross effect vanishes.
We recall that, for $d \in \mathds{N}$, the $(d+1)$-th cross effect is a functor $$\mathrm{cr}_{d+1} (F) : \mathcal{C}^{\times d+1} \to \mathds{k}\mbox{-}\mathsf{Mod},$$ whose value on objects $X_1, \cdots, X_{d+1}$ is defined to be the joint kernel of the maps $F(\epsilon_{X_1,\cdots,X_{d+1};k} )$ for $1 \leq k \leq d+1$.

A $\mathcal{C}$-module is {\it polynomial of degree $d$}, denoted by $\deg F = d$, if $\deg F \leq d$ but $\deg F \not\leq d-1$.

Cross effects admit the following alternative description.
Note that the idempotents $e_{X_1,\cdots,X_{d+1};k}$ on $\bigvee_{i=1}^{d+1} X_i$ for $1 \leq k \leq d+1$ commute pairwise, and so do the idempotents $F (e_{X_1,\cdots,X_{d+1};k})$ on $F(\bigvee_{i=1}^{d+1} X_i)$.
It follows that the composite
\begin{align} \label{202609141928}
    P_{X_1,\cdots,X_{d+1};F} {:=} (\mathrm{id}-F (e_{X_1,\cdots,X_{d+1};1})) \circ \cdots \circ (\mathrm{id}-F (e_{X_1,\cdots,X_{d+1};d+1})) 
\end{align}
is again an idempotent and its image equals $\mathrm{cr}_{d+1} (F) (X_1, \cdots, X_{d+1})$.

A functor $F : \mathcal{C} \to \mathds{k}\mbox{-}\mathsf{Mod}$ is {\it linear} if $F$ is polynomial of degree at most $1$, i.e. $\mathrm{cr}_2(F) \cong0$, and $F(\ast) \cong 0$.
For a linear functor $F$, we have $F(X \vee Y) \cong F(X) \oplus F(Y)$ (see e.g.  \cite[Proposition 2.4]{hartl2015polynomial}).
For $r \geq 1$, a functor $H : \mathcal{C}^{\times r} \to \mathds{k}\mbox{-}\mathsf{Mod}$ is {\it $r$-multilinear} if it is linear in each variable.

\begin{remark} \label{202609282124}
    For an endomorphism $f: X \to X$ and $\epsilon \in \{0,1\}$, we set $f^\epsilon {:=} f$ if $\epsilon = 1$; and $f^\epsilon {:=} \mathrm{id}_{X}$ if $\epsilon = 0$.
    Expanding \eqref{202609141928}, one may verify that $\deg F \leq d$ if and only if, for any objects $X_1, \cdots, X_{d+1}$ of $\mathcal{C}$,
    \begin{align*}
        P_{X_1,\cdots,X_{d+1};F} = \sum_{\epsilon \in \{0,1\}^{d+1}} (-1)^{|\epsilon|} \cdot F( e_{X_1}^{\epsilon_1} \vee e_{X_2}^{\epsilon_2} \vee \cdots \vee e_{X_{d+1}}^{\epsilon_{d+1}} ) = 0 .
    \end{align*}
    Here, $|\epsilon| = \sum^{d+1}_{j=1} \epsilon_j$.
    The description above is analogous to the definition of polynomial maps in the sense of Passi \cite[Definition 1.1, Chapter V]{MR537126}.
    We will show that, for a monoidal category $\mathcal{C}$ satisfying some conditions, we can describe a polynomial functor by using Passi's definition (see Corollary \ref{202509251010}).
\end{remark}

\begin{Lemma} \label{202609150934}
    Let $X_{1},\cdots,X_{n+m} \in \mathrm{Obj}(\mathcal{C})$.
    For a $\mathcal{C}$-module $F$, define a $\mathcal{C}$-module $G$ by
    $$G(X) = \mathrm{cr}_{m+1}  (F) (X,X_{n+1}, \cdots,X_{n+m}).$$
    For $Y = \bigvee^n_{i=1} X_i$, the following diagram commutes
    $$
    \begin{tikzcd}
        F(\bigvee^{n+m}_{i=1} X_i) \ar[r, "P_{X_1\cdots,X_{n+m};F}"] \ar[d, "P_{Y,X_{n+1},\cdots,X_{n+m};F}"'] & F(\bigvee^{n+m}_{i=1} X_i) \\
        G(\bigvee^{n}_{i=1} X_i) \ar[r,"P_{X_1\cdots,X_{n};G}"'] & G(\bigvee^{n}_{i=1} X_i) \ar[u,hookrightarrow]
    \end{tikzcd} .
    $$
\end{Lemma}
\begin{proof}
    By definition, $$G(f) = F(f \vee \mathrm{id}_{X_{n+1} \vee \cdots \vee X_{n+m}}) \circ P_{Z,X_{n+1},\cdots,X_{n+m};F} = P_{W,X_{n+1},\cdots,X_{n+m};F} \circ F(f \vee \mathrm{id}_{X_{n+1} \vee \cdots \vee X_{n+m}})$$ for any morphism $f:Z \to W$. 
    Hence, 
    \begin{align*}
        &P_{X_1\cdots,X_{n};G} \circ P_{Y,X_{n+1},\cdots,X_{n+m};F} = (\mathrm{id}-G (e_{X_1, \cdots,X_n;1})) \circ \cdots \circ (\mathrm{id}-G (e_{X_1,\cdots,X_n;n})) \circ P_{Y,X_{n+1},\cdots,X_{n+m};F} , \\
        =& (\mathrm{id}-F (e_{X_1, \cdots,X_n;1}  \vee \mathrm{id}_{X_{n+1} \vee \cdots \vee X_{n+m}})) \circ \cdots \circ (\mathrm{id}-F (e_{X_1,\cdots,X_n;n} \vee \mathrm{id}_{X_{n+1} \vee \cdots \vee X_{n+m}})) \circ P_{Y,X_{n+1},\cdots,X_{n+m};F} , \\
        =& (\mathrm{id}-F (e_{X_1, \cdots,X_{n+m};1})) \circ \cdots \circ (\mathrm{id}-F (e_{X_1,\cdots,X_{n+m};n})) \circ P_{Y,X_{n+1},\cdots,X_{n+m};F} .
    \end{align*}
    For $2 \leq j \leq m+1$, $e_{Y, X_{n+1},\cdots,X_{n+m};j} = e_{X_1,\cdots,X_{n+m};n+j-1}$.
    On the other hand, the obvious relation $e_{X_1,\cdots,X_{n+m};n} \circ e_{Y, X_{n+1},\cdots,X_{n+m};1} = e_{Y, X_{n+1},\cdots,X_{n+m};1}$ implies
    \begin{align*}
    ( \mathrm{id} - F (e_{X_1,\cdots,X_{n+m};n})) \circ  (\mathrm{id}-F (e_{Y, X_{n+1},\cdots,X_{n+m};1})) 
    = \mathrm{id} - F (e_{X_1,\cdots,X_{n+m};n}) .
    \end{align*}
    Therefore, by the above computation of $P_{X_1\cdots,X_{n};G} \circ P_{Y,X_{n+1},\cdots,X_{n+m};F}$, we obtain the commutative diagram in the statement.
\end{proof}

\begin{Lemma} \label{202609142110}
    Under the notation and assumptions of the preceding lemma, we have
    $$\mathrm{cr}_n ( G)(X_1, \cdots, X_n) \cong  \mathrm{cr}_{n+m}(F) (X_1,\cdots, X_{n+m}).$$
\end{Lemma}
\begin{proof}
    It follows from Lemma \ref{202609150934} and the definition of cross effects.
\end{proof}

\begin{Lemma} \label{202609281530}
    The full subcategory of $\mathcal{C}$-modules of degree $\leq d$ is closed under subobjects, quotient objects, direct sums and products.
\end{Lemma}
\begin{proof}
    The idempotent \eqref{202609141928} is natural in $F$, and it is compatible with direct sums and products, which are computed objectwise.
    Since $\deg F \leq d$ if and only if $P_{X_1,\cdots,X_{d+1};F}=0$ for all $X_1,\cdots,X_{d+1}$, the assertion follows.
\end{proof}

For $d \in \mathds{N}$, {\it the $d$-th polynomial approximation of $F$} is defined as the largest
subfunctor $P_{d}(F) \subset F$ of degree at most $d$; it exists since the class of functors of
degree at most $d$ is closed under quotients and direct sums (see Lemma \ref{202609281530}).
This makes $P_d$ a right adjoint to the inclusion of the category of $\mathcal{C}$-modules of degree
at most $d$ into that of all $\mathcal{C}$-modules.
\begin{prop} \label{202609281508}
    A $\mathcal{C}$-module of degree $\leq d$ is of degree $\leq d+1$.
\end{prop}
\begin{proof}
    Let $F$ be a $\mathcal{C}$-module of degree $\leq d$.
    Lemma \ref{202609142110} applied to this $F$ with $m= d$ and $n=2$ proves $\mathrm{cr}_{d+2}(F) = 0$, since $G (X) = \mathrm{cr}_{d+1}  (F) (X,X_{3}, \cdots,X_{d+2}) =0$.
\end{proof}

Hence by Proposition \ref{202609281508}, there is an increasing filtration 
\begin{align} \label{202609281811}
    P_0(F) \subseteq P_1(F) \subseteq \cdots \subseteq F .
\end{align}
A $\mathcal{C}$-module $F$ is {\it analytic} if $F = \bigcup_{d\geq 0} P_d (F)$.

\subsection{A dichotomy for $\mathrm{D}(\mathcal{C})$}

Throughout this section, let $\mathcal{C}$ be a monoidal category with zero unit.
In this section, we address question (Q1) from the introduction. 
To this end, we first give the precise definition of the set $\mathrm{D}(\mathcal{C})$.
\begin{Defn} \label{202511251137}
    We define $\mathrm{D}(\mathcal{C})$ to be the set of $d \in \mathds{N}$ for which there exists a polynomial $\mathcal{C}$-module of degree equal to $d$.
\end{Defn}

\begin{Example} \label{202609211106}
    Constant functors have degree zero, so $0 \in \mathrm{D}(\mathcal{C})$.
    In particular, $\mathrm{D}(\mathcal{C}) \neq \emptyset$.
\end{Example}

\begin{Lemma}
    If $d \in \mathrm{D}(\mathcal{C})$, then $d=0$ or $(d-1) \in \mathrm{D}(\mathcal{C})$.
\end{Lemma}
\begin{proof}
    It suffices to consider $d>0$.
    By Example \ref{202609211106}, we can assume that $d \geq 2$.
    Choose a functor $F$ of $\deg = d$.
    Then $\mathrm{cr}_{d}(F) \neq 0$, so that there exist $Y_1,\cdots,Y_d$ such that $\mathrm{cr}_{d}(F)(Y_1,\cdots,Y_d) \not\cong 0$.
    Let $G(X) = \mathrm{cr}_2 (F)(X,Y_d)$.
    By Lemma \ref{202609142110}, we obtain $\mathrm{cr}_d (G) (X_1,\cdots,X_d) \cong \mathrm{cr}_{d+1}(F) (X_1,\cdots,X_d,Y_d) \cong 0$, while $\mathrm{cr}_{d-1}(G) \neq 0$, since
    $$\mathrm{cr}_{d-1}(G) (Y_1,\cdots,Y_{d-1}) \cong \mathrm{cr}_{d}(F) (Y_1,\cdots,Y_{d-1},Y_d) \neq 0.$$
    Hence, $\deg G = d-1$.
    This proves $d-1 \in \mathrm{D}(\mathcal{C})$.
\end{proof}

By the preceding lemma, $\mathrm{D}(\mathcal{C})$ consists of consecutive nonnegative integers starting at zero. 
It is then natural to ask what its supremum is; the following proposition answers this question and proves part of Theorem \ref{2026009152033}.

\begin{theorem} \label{202609142120}
    The set $\mathrm{D}(\mathcal{C})$ is either $\mathds{N}$ or $\{0\}$.
\end{theorem}

Both cases of this dichotomy occur: see Theorem \ref{202510171045} for $\mathrm{D}(\mathcal{C})=\mathds{N}$, and Examples \ref{202509241443} and \ref{202509301321} for $\mathrm{D}(\mathcal{C})=\{0\}$.
See also Theorem \ref{202510171046}.
The proof of the theorem relies on the following lemma:
\begin{Lemma} \label{202609142131}
    For $d \in \mathds{N}^\ast$, $d \in \mathrm{D}(\mathcal{C})$ if and only if there exists a nonzero $d$-multilinear functor $\mathcal{C}^{\times d} \to \mathds{k}\mbox{-}\mathsf{Mod}$.
\end{Lemma}
\begin{proof}
    Let $d \in \mathds{N}^\ast$.
    For a polynomial $\mathcal{C}$-module $F$ of degree at most $d$, the $d$-th cross effect $\mathrm{cr}_d(F)$ is $d$-multilinear.
    This is well known when the monoidal structure on $\mathcal{C}$ is given by (co)products; the general case can be proved in the same way, using Lemma~\ref{202609142110}.
    Therefore, if $\deg F = d$, then $\mathrm{cr}_d (F)$ is a nontrivial $d$-multilinear functor.
    
    We now prove the converse.
    Let $H$ be a nonzero $d$-multilinear functor, and choose $Y_1,\cdots,Y_d$ with $H(Y_1,\cdots,Y_d) \neq 0$.
    Define a $\mathcal{C}$-module $F$ by $F(X) := H(X,\cdots,X)$; we show $\deg F = d$.
    Fix $X_1,\cdots,X_n \in \mathrm{Obj}(\mathcal{C})$ and let $\iota_k : X_k \to \bigvee_{i=1}^n X_i$ be the insertion obtained from $\eta_{X_i}$ for $i \neq k$ and $\mathrm{id}_{X_k}$.
    By the linearity of $H$ in each variable, we obtain an isomorphism
    $$
    \bigoplus_{\mathbf{i} \in [n]^d} H(X_{i_1},\cdots,X_{i_d}) \xrightarrow{\ \cong\ } H\Big(\bigvee_{j=1}^n X_j,\cdots,\bigvee_{j=1}^n X_j\Big) = F\Big(\bigvee_{j=1}^n X_j\Big), \qquad \mathbf{i}=(i_1,\cdots,i_d),
    $$
    whose restriction to the $\mathbf{i}$-th summand is $H(\iota_{i_1},\cdots,\iota_{i_d})$.
    For $1 \leq k \leq n$, we have $e_{X_1,\cdots,X_n;k} \circ \iota_i = \iota_i$ if $i \neq k$, whereas $e_{X_1,\cdots,X_n;k} \circ \iota_k$ factors through $\ast$.
    Since $F(e_{X_1,\cdots,X_n;k}) = H(e_{X_1,\cdots,X_n;k},\cdots,e_{X_1,\cdots,X_n;k})$ and $H$ vanishes whenever one of its arguments is $\ast$, the idempotent $F(e_{X_1,\cdots,X_n;k})$ acts on the $\mathbf{i}$-th summand as the identity if $k \notin \{i_1,\cdots,i_d\}$, and as zero otherwise.
    Hence, $P_{X_1,\cdots,X_n;F}$ is the projection onto the summands indexed by $\mathbf{i}$ with $\{i_1,\cdots,i_d\} = \{1,\cdots,n\} {=:} \underline{n}$, that is,
    $$
    \mathrm{cr}_n(F)(X_1,\cdots,X_n) \cong \bigoplus_{\mathbf{i} : \underline{d} \twoheadrightarrow \underline{n}} H(X_{i_1},\cdots,X_{i_d}).
    $$
    For $n=d+1$, there is no surjection $\underline{d} \to \underline{d+1}$, so $\mathrm{cr}_{d+1}(F)=0$ and $\deg F \leq d$.
    For $n=d$, the surjections $\underline{d} \to \underline{d}$ are the bijections, and the summand indexed by the identity gives $\mathrm{cr}_d(F)(Y_1,\cdots,Y_d) \supseteq H(Y_1,\cdots,Y_d) \neq 0$.
    Therefore $\deg F = d$.
\end{proof}

\begin{proof}[Proof of Theorem \ref{202609142120}]
    Suppose that $\mathrm{D}(\mathcal{C}) \neq \{0\}$.
    We show that $\mathrm{D}(\mathcal{C}) = \mathds{N}$.
    Let $0 \neq d \in \mathrm{D}(\mathcal{C})$.
    By Lemma \ref{202609142131}, there exists a nonzero $d$-multilinear functor $H$.
    Choose objects $Y_1,\cdots,Y_d$ in $\mathcal{C}$ such that 
    $$M {:=} H(Y_1,\cdots,Y_d) \neq 0.$$
    We define a $d$-multilinear functor $H'$ by $H' (X_1,\cdots,X_d) = \hom_{\mathds{k}} (M, H(X_1,\cdots,X_d))$.
    Here, $\hom_{\mathds{k}}$ denotes the $\mathds{k}$-module consisting of $\mathds{k}$-linear maps.
    Then $H'$ is nonzero since $H'(Y_1,\cdots,Y_d) = \hom_{\mathds{k}} (M,M)$ is a nonzero $\mathds{k}$-algebra by the condition $M \neq 0$.
    Let 
    $$G(X_1,\cdots,X_{d+1}) {:=} H'(X_1,\cdots,X_d) \otimes H'(Y_1,\cdots,Y_{d-1},X_{d+1}).$$
    Then $G$ is $(d+1)$-multilinear, since $H'$ is linear in each variable.
    Moreover, $G(Y_1,\cdots,Y_d,Y_d) = \hom_{\mathds{k}} (M,M)^{\otimes 2}$ is not trivial, since the composition $\hom_{\mathds{k}} (M,M)^{\otimes 2} \to \hom_{\mathds{k}} (M,M) \neq 0$ is surjective.
    Hence, $G$ is a nonzero $(d+1)$-multilinear functor.
    By Lemma \ref{202609142131}, $d+1 \in \mathrm{D}(\mathcal{C})$.
    This result recursively proves that $\mathrm{D}(\mathcal{C}) = \mathds{N}$.
\end{proof}

We conclude this section by giving another application of Lemma \ref{202609142131}.
Let $\mathcal{C}^{\mathsf{o}}$ denote the opposite category of $\mathcal{C}$, which is again a monoidal category with zero unit.

\begin{prop} \label{202512231721}
    $\mathrm{D}(\mathcal{C}) = \mathrm{D}(\mathcal{C}^{\mathsf{o}})$.
\end{prop}
\begin{proof}
    It suffices to prove that $\mathrm{D}(\mathcal{C}) \subset \mathrm{D}(\mathcal{C}^{\mathsf{o}})$, since the reverse inclusion follows by $(\mathcal{C}^{\mathsf{o}})^{\mathsf{o}} = \mathcal{C}$.
    By Theorem \ref{202609142120}, $\mathrm{D}(\mathcal{C})$ is either $\{0\}$ or $\mathds{N}$.
    In the former case, the proof is complete, since $\mathrm{D}(\mathcal{C})$ and $\mathrm{D}(\mathcal{C}^{\mathsf{o}})$ contain $0$ (see Example \ref{202609211106}).
    Suppose that $\mathrm{D}(\mathcal{C})=\mathds{N}$.
    Let $0 < d \in \mathrm{D}(\mathcal{C})$.
    By Lemma \ref{202609142131}, there exists a nonzero $d$-multilinear functor $H$.
    For some $Y_1,\cdots,Y_d$, we have $M {:=} H(Y_1,\cdots,Y_d) \neq 0$.
    We define a functor $H' : (\mathcal{C}^{\mathsf{o}})^{\times d} \to \mathds{k}\mbox{-}\mathsf{Mod}$ by $$H'(X_1,\cdots,X_d) {:=} \hom_{\mathds{k}} (H(X_1,\cdots,X_d), M).$$
    Since $H'(Y_1,\cdots,Y_d) = \hom_{\mathds{k}} (M, M)$ is the endomorphism algebra of $M\neq 0$, it is nonzero.
    So, the functor $H'$ is nonzero.
    As $H$ is linear in each variable, so is $H'$.
    Indeed, $\hom_{\mathds{k}}(-, M)$ is additive, so it sends the decomposition $H(X \vee X', \dots) \cong H(X, \dots) \oplus H(X', \dots)$ induced by the canonical maps in $\mathcal{C}$ to the corresponding decomposition of $H'$ induced by the canonical maps in $\mathcal{C}^{\mathsf{o}}$.
    By Lemma \ref{202609142131}, $d \in \mathrm{D}(\mathcal{C}^{\mathsf{o}})$, so $\mathrm{D}(\mathcal{C}) \subset \mathrm{D}(\mathcal{C}^{\mathsf{o}})$.
\end{proof}

\section{Bi-indexed modules}
\label{202410161743}

The objective of this section is to present a toolkit for the description of our results.
We introduce a bicategory $\mathsf{Mat}_{\mathds{k}}$ formed by bi-indexed modules, and investigate monad theory within the bicategory $\mathsf{Mat}_{\mathds{k}}$.
Monads in that bicategory are equivalent to $\mathds{k}$-linear categories; and modules over such a monad are equivalent to $\mathds{k}$-linear functors from the corresponding $\mathds{k}$-linear category to $\mathds{k}\mbox{-}\mathsf{Mod}$.
One advantage of working within the bicategory is that it enables us to treat $\mathds{k}$-linear categories and modules over them in a way that parallels classical ring–module theory.

The reader is referred to \cite{benabou2006introduction} for notions of bicategory theory.
The notion of a monad in a bicategory is a generalization of that of an algebra.

\subsection{The bicategory of matrices $\mathsf{Mat}_{\mathds{k}}$}

For a set $\mathcal{Z}$, a {\it $\mathcal{Z}$-indexed ($\mathds{k}$-)module} is a family of $\mathds{k}$-modules $\mathtt{F} = \{ \mathtt{F} (Z) \}_{Z\in \mathcal{Z}}$ indexed by the set $\mathcal{Z}$.
For $\mathcal{Z}$-indexed modules $\mathtt{F}$ and $\mathtt{F}^\prime$, a {\it $\mathcal{Z}$-homomorphism}, or merely a {\it map}, $\xi : \mathtt{F}^\prime \to \mathtt{F}$ is a family of $\mathds{k}$-linear maps $\xi = \xi_{Z} : \mathtt{F}^\prime (Z) \to \mathtt{F}(Z)$ parametrized by $Z \in \mathcal{Z}$.
It is obvious that the category of $\mathcal{Z}$-indexed modules and their homomorphisms constitutes an abelian category, with the (co)kernel objects calculated indexwise.
A homomorphism $\xi : \mathtt{F}^\prime  \to \mathtt{F}$ is an {\it epimorphism} (a {\it monomorphism}, resp.) if each map $\xi_{Z} : \mathtt{F}^\prime (Z) \to \mathtt{F} (Z)$ is surjective (injective, resp.).

Most concepts associated with $\mathds{k}$-modules are readily extended to indexed modules:
\begin{itemize}
    \item A {\it submodule} of $\mathtt{F}$ is a $\mathcal{Z}$-indexed $\mathds{k}$-module $\mathtt{G}$ such that $\mathtt{G}(Z)$ is a submodule of $\mathtt{F} (Z)$ for $Z \in \mathcal{Z}$.
    We write $\mathtt{G} \subset \mathtt{F}$ in this case.
    \item A {\it quotient module} $\mathtt{F}/\mathtt{G}$ is defined to be a $\mathcal{Z}$-indexed module such that $\left( \mathtt{F}/\mathtt{G} \right) (Z) {:=} \mathtt{F} (Z)/\mathtt{G} (Z)$.
    For $Z \in \mathcal{Z}$ and $f \in \mathtt{F}(Z)$, we denote the induced element of $(\mathtt{F}/\mathtt{G}) (Z)$ by $f \mod{\mathtt{G}}$.
\end{itemize}
In this paper, we freely make use of these concepts and symbols.

We now focus on bi-indexed modules.
The reader may find it helpful to think of a $(\mathcal{Y} \times \mathcal{X})$-indexed module $\mathtt{F}$ as a matrix whose entries are $\mathds{k}$-modules, with rows indexed by $\mathcal{Y}$ and columns indexed by $\mathcal{X}$; thus $\mathtt{F}(Y,X)$ is the $(Y,X)$-entry.
Let $\mathtt{F}$ be a $(\mathcal{Y}\times\mathcal{X})$-indexed $\mathds{k}$-module and $\mathtt{G}$ a $(\mathcal{Z}\times\mathcal{Y})$-indexed $\mathds{k}$-module.
The \emph{matrix product} of $\mathtt{G}$ and $\mathtt{F}$, denoted by $\mathtt{G} \otimes \mathtt{F}$, is the $(\mathcal{Z}\times\mathcal{X})$-indexed $\mathds{k}$-module defined by
$$
(\mathtt{G} \otimes \mathtt{F}) (Z,X) := \bigoplus_{Y \in \mathcal{Y}} \mathtt{G}(Z,Y) \otimes \mathtt{F}(Y,X), \qquad (Z,X) \in \mathcal{Z} \times \mathcal{X},
$$
in analogy with the usual matrix product.
The matrix product is associative up to coherent natural isomorphism.
Using it, we introduce the following bicategory.
\begin{Defn} \label{202410021236}
    We define the bicategory $\mathsf{Mat}_{\mathds{k}}$ as follows:
    \begin{itemize}
        \item Its {\it objects} are sets.
        We use symbols $\mathcal{X},\mathcal{Y},\mathcal{Z}$ to denote its objects.
        \item For objects $\mathcal{X},\mathcal{Y}$, a {\it 1-morphism} from $\mathcal{X}$ to $\mathcal{Y}$ is given by a $(\mathcal{Y} \times \mathcal{X})$-indexed $\mathds{k}$-module.
        The notation $\mathtt{F} : \mathcal{X} \rightsquigarrow \mathcal{Y}$ denotes such a 1-morphism $\mathtt{F} = \{ \mathtt{F} (Y,X) \}_{(Y,X)\in\mathcal{Y}\times\mathcal{X}}$.
        \item
        For 1-morphisms $\mathtt{F} : \mathcal{X} \rightsquigarrow \mathcal{Y}$ and $\mathtt{G} : \mathcal{Y} \rightsquigarrow \mathcal{Z}$, the {\it composition} is given by the matrix product $\mathtt{G} \otimes \mathtt{F} : \mathcal{X} \rightsquigarrow \mathcal{Z}$.
        \item
        For each object $\mathcal{X}$, the {\it identity morphism} $\mathds{I}_{\mathcal{X}} : \mathcal{X} \rightsquigarrow \mathcal{X}$ on $\mathcal{X}$ is defined to be
        \begin{align*}
            \mathds{I}_{\mathcal{X}} (Y,X) {:=} 
            \begin{cases}
                \mathds{k} , ~ (Y=X) , \\
                0 , ~ (Y \neq X) .
            \end{cases}
        \end{align*}
        \item For 1-morphisms $\mathtt{F}$ and $\mathtt{G}$ from $\mathcal{X}$ to $\mathcal{Y}$, a {\it 2-morphism} $\Phi$ from $\mathtt{F}$ to $\mathtt{G}$, denoted as $\Phi : \mathtt{F} \to \mathtt{G}$, is defined to be a $(\mathcal{Y}\times \mathcal{X})$-homomorphism.
        \item 
        The {\it identity 2-morphism} $\mathrm{id}_{\mathtt{F}}$ on $\mathtt{F}$ is given by the identity on the $(\mathcal{Y} \times \mathcal{X})$-indexed module $\mathtt{F}$.
    \end{itemize}
\end{Defn}

In this paper, we mainly use the terminology {\it homomorphism} or {\it map} ({\it isomorphism}, resp.) of indexed modules instead of 2-morphisms (2-isomorphism, resp.) in $\mathsf{Mat}_{\mathds{k}}$.

\begin{remark}    
    The bicategory $\mathsf{Mat}_{\mathds{k}}$ is a special case of the \emph{bicategory of $\mathbb{W}$-matrices} introduced in \cite[Section 1]{betti1983variation}, where $\mathbb{W}$ is the bicategory with one object, $\mathds{k}$-modules as 1-morphisms and $\mathds{k}$-linear maps as 2-morphisms.
    Our convention of writing the target first, $\mathtt{F}(Y,X)$ for $\mathtt{F} : \mathcal{X} \rightsquigarrow \mathcal{Y}$, is that of \cite{betti1983variation} and of linear algebra: it is the only convention under which the indices in $\mathtt{G} \otimes \mathtt{F}$ are adjacent, as in the matrix product, while $\mathtt{G} \otimes \mathtt{F}$ still means ``$\mathtt{F}$ followed by $\mathtt{G}$''.
    Note that this is opposite to the usual notation $\mathcal{C}(X,Y)$ for morphism sets, where the source is written first.
    We keep the latter for ordinary categories.
\end{remark}
As for ordinary matrices, we have a notion of transpose:
\begin{Defn} \label{202410131009}
    For a 1-morphism $\mathtt{F} : \mathcal{X} \rightsquigarrow \mathcal{Y}$, the \emph{transpose} of $\mathtt{F}$ is the 1-morphism $\mathtt{F}^{\mathrm{t}} : \mathcal{Y} \rightsquigarrow \mathcal{X}$ defined by $\mathtt{F}^{\mathrm{t}}(X,Y) := \mathtt{F}(Y,X)$ for $X \in \mathcal{X}$, $Y \in \mathcal{Y}$.
\end{Defn}

The following are clear:
\begin{align*}
    \mathds{I}_{\mathcal{X}}^{\mathrm{t}} = \mathds{I}_{\mathcal{X}}, \quad (\mathtt{G} \otimes \mathtt{F})^{\mathrm{t}} \cong \mathtt{F}^{\mathrm{t}} \otimes \mathtt{G}^{\mathrm{t}} , \quad (\mathtt{F}^{\mathrm{t}})^{\mathrm{t}} = \mathtt{F} .
\end{align*}
In particular, the transposition gives an anti-automorphism on the bicategory $\mathsf{Mat}_{\mathds{k}}$ of order $2$.

\subsection{Monads in $\mathsf{Mat}_{\mathds{k}}$}
\label{202408121736}

A monad in $\mathsf{Mat}_{\mathds{k}}$ on an object $\mathcal{X}$ consists of a 1-morphism $\mathtt{T} : \mathcal{X} \rightsquigarrow \mathcal{X}$ and 2-morphisms $\nabla : \mathtt{T} \otimes \mathtt{T} \to \mathtt{T}$ and $\eta : \mathds{I}_{\mathcal{X}} \to \mathtt{T}$, for which the following diagrams are commutative:
$$
\begin{tikzcd}
    \mathtt{T} \otimes \mathtt{T} \otimes \mathtt{T} \ar[r, "\mathrm{id}_{\mathtt{T}} \otimes \nabla"] \ar[d, "\nabla \otimes \mathrm{id}_{\mathtt{T}}"] & \mathtt{T} \otimes \mathtt{T} \ar[d, "\nabla"] & T \otimes \mathds{I}_{\mathcal{X}} \ar[d, "\mathrm{id}_{\mathtt{T}\otimes\eta}"] \ar[r, "\cong"] & \mathtt{T}  \ar[d, equal] & \mathds{I}_{\mathcal{X}} \otimes \mathtt{T} \ar[d, "\eta \otimes \mathrm{id}_{\mathtt{T}}"] \ar[l, "\cong"']\\
    \mathtt{T} \otimes \mathtt{T}  \ar[r, "\nabla"] & \mathtt{T} & \mathtt{T} \otimes \mathtt{T} \ar[r, "\nabla"] & \mathtt{T} & \mathtt{T} \otimes \mathtt{T} \ar[l, "\nabla"']
\end{tikzcd}
$$
For monads $\mathtt{T}$ and $\mathtt{S}$ on $\mathcal{X}$, a monad map $\xi: \mathtt{T} \to \mathtt{S}$ is a map of $(\mathcal{X}\times\mathcal{X})$-indexed $\mathds{k}$-modules which preserves the monad structures.

\begin{notation} \label{202410232034}
    Let $X,Y,Z \in \mathcal{X}$.
    For $f \in \mathtt{T} (Z,Y)$ and $g \in \mathtt{T}(Y,X)$, we use the {\it composition} notation for the monad operation:
    $$
    f \circ g {:=} \nabla (f \otimes g) , \quad \mathrm{where~~}f \otimes g \in \mathtt{T} (Z,Y) \otimes \mathtt{T}(Y,X) \subset \left( \mathtt{T} \otimes \mathtt{T} \right) (Z,X) .
    $$
    For $X \in \mathcal{X}$, we also denote by $1_X \in \mathtt{T} (X,X)$ the element corresponding to the unit $1_{\mathds{k}} \in \mathds{k}$ via the monad unit $\eta : \mathds{I}_{\mathcal{X}} \to \mathtt{T}$.
\end{notation}

The definition of a monad implies that $( f \circ g) \circ h= f \circ (g \circ h)$ and $f \circ 1_X = f = 1_Y \circ f$.
Hence, the monad $\mathtt{T}$ induces a (small) $\mathds{k}$-linear category whose object set is $\mathcal{X}$ and $\mathtt{T}(Y,X)$ is the morphism set from $X$ to $Y$.
Conversely, every monad in $\mathsf{Mat}_{\mathds{k}}$ on $\mathcal{X}$ arises in this way.
This observation naturally leads to the following equivalence of categories:
\begin{align} \label{202401231608}
    \mathsf{Lin}^{\mathcal{X}}_{\mathds{k}} \to \mathsf{Mnd}^{\mathcal{X}},
\end{align}
where $\mathsf{Mnd}^{\mathcal{X}}$ is the category of monads in $\mathsf{Mat}_{\mathds{k}}$ on $\mathcal{X}$ and monad maps;
and $\mathsf{Lin}^{\mathcal{X}}_{\mathds{k}}$ the category of $\mathds{k}$-linear categories with object set $\mathcal{X}$ and $\mathds{k}$-linear functors which are the identity on objects.
In particular, this assigns a monad map to a functor between $\mathds{k}$-linear categories.

In the following, we present examples of monads in $\mathsf{Mat}_{\mathds{k}}$.

\begin{Example} \label{202404111056}
    If $\mathcal{X} = \ast$ is the one-point set, then the category $\mathsf{Lin}^{\mathcal{X}}_{\mathds{k}}$ is equivalent to the category of $\mathds{k}$-algebras.
    Hence, monads in $\mathsf{Mat}_{\mathds{k}}$ on $\ast$ coincide with unital $\mathds{k}$-algebras.
\end{Example}

\begin{Example} \label{202401271351}
    Let $\mathcal{C}$ be a category with object set $\mathcal{X}$.
    The $\mathds{k}$-linearization of $\mathcal{C}$ gives a $\mathds{k}$-linear category, which induces a monad on $\mathcal{X}$ via (\ref{202401231608}).
    We denote by $\mathtt{L}_{\mathcal{C}}$ the induced monad.
    In particular, $$\mathtt{L}_{\mathcal{C}} (Y,X) = \mathds{k} \mathcal{C}(X,Y).$$
\end{Example}

\begin{Example} \label{202508152025}
    Let $\mathtt{T}$ be a monad on $\mathcal{X}$.
    The transposition introduced in Definition \ref{202410131009} induces a monad $\mathtt{T}^{\mathrm{t}}$ on $\mathcal{X}$ by the formula $\nabla_{\mathtt{T}^{\mathrm{t}}} (f \otimes g) {:=} \nabla_{\mathtt{T}} ( g \otimes f)$.
\end{Example}

\subsection{Modules in $\mathsf{Mat}_{\mathds{k}}$} \label{202402281443}

In this section, we study modules over monads in the bicategory $\mathsf{Mat}_{\mathds{k}}$.
Consider monads $\mathtt{T},\mathtt{S}$ in $\mathsf{Mat}_{\mathds{k}}$ on objects $\mathcal{X}$ and $\mathcal{Y}$ respectively.
A {\it $(\mathtt{T},\mathtt{S})$-bimodule} is a morphism $\mathtt{M}: \mathcal{Y} \rightsquigarrow \mathcal{X}$ endowed with 2-morphisms $\rhd : \mathtt{T} \otimes \mathtt{M} \to \mathtt{M}$ and $\lhd : \mathtt{M} \otimes \mathtt{S} \to \mathtt{M}$.
It is required that $\rhd$ gives a left $\mathtt{T}$-action, $\lhd$ gives a right $\mathtt{S}$-action and that the following diagram commutes:
$$
\begin{tikzcd}
    \mathtt{T} \otimes \mathtt{M} \otimes \mathtt{S} \ar[r, "\rhd \otimes \mathrm{id}_{\mathtt{S}}"] \ar[d, "\mathrm{id}_{\mathtt{T}} \otimes \lhd"] & \mathtt{M} \otimes \mathtt{S} \ar[d, "\lhd"] \\
    \mathtt{T} \otimes \mathtt{M} \ar[r, "\rhd"] & \mathtt{M} 
\end{tikzcd} .
$$
For $(\mathtt{T},\mathtt{S})$-bimodules $\mathtt{M},\mathtt{M}^\prime$, a {\it $(\mathtt{T},\mathtt{S})$-homomorphism} $\xi : \mathtt{M} \to \mathtt{M}^\prime$ is a 2-morphism which is compatible with the actions of the monads $\mathtt{T}$ and $\mathtt{S}$.

It is useful to understand the bimodule structure at the level of its elements.
Let $\mathtt{M}$ be a $(\mathtt{T},\mathtt{S})$-bimodule as above.
For $X, X^\prime \in \mathcal{X}$ and $Y \in \mathcal{Y}$, the left $\mathtt{T}$-action $\rhd : \mathtt{T} \otimes \mathtt{M} \to \mathtt{M}$ gives a map
$$\mathtt{T} (X^\prime,X) \otimes \mathtt{M} (X, Y) \to \mathtt{M} (X^\prime ,Y) ; \quad f \otimes v \mapsto f \rhd v . $$
By the axioms of the left action, we have:
\begin{align*}
    g \rhd ( f \rhd v) = (g \circ f) \rhd v, \qquad 1_{X} \rhd v = v .
\end{align*}
Analogously, the right action $\lhd$ satisfies $(v \lhd f^\prime) \lhd g^\prime = v \lhd (f^\prime \circ g^\prime)$ and $v \lhd 1_Y = v$ for appropriate $f^\prime,g^\prime$.
Furthermore, the commutativity of the left and right actions implies 
$$(f \rhd v) \lhd f^\prime = f \rhd ( v \lhd f^\prime) . $$

Let $\mathtt{U} = \{\mathtt{U}(X,Y)\}_{(X,Y) \in \mathcal{X} \times \mathcal{X}}$ be a family of subsets $\mathtt{U}(X,Y) \subset \mathtt{T}(X,Y)$, and $\mathtt{N} = \{\mathtt{N}(X,Y)\}_{(X,Y) \in \mathcal{X} \times \mathcal{Y}}$ a family of subsets $\mathtt{N}(X,Y) \subset \mathtt{M}(X,Y)$.
We denote by $\mathtt{U} \rhd \mathtt{N}$ the submodule of $\mathtt{M}$ whose component at $(X,Y) \in \mathcal{X}\times \mathcal{Y}$ is the $\mathds{k}$-submodule of $\mathtt{M}(X,Y)$ generated by the elements $f \rhd g$ with $Z \in \mathcal{X}$, $f \in \mathtt{U}(X,Z)$ and $g \in \mathtt{N}(Z,Y)$.
For $h \in \mathtt{M}(X,Y)$, we write $\mathtt{U} \rhd h := \mathtt{U} \rhd \mathtt{N}$, where $\mathtt{N}(X',Y') = \{h\}$ if $(X',Y')=(X,Y)$ and $\mathtt{N}(X',Y') = \varnothing$ otherwise.

\begin{Defn}
    A monad $\mathtt{T}$ has a $(\mathtt{T},\mathtt{T})$-bimodule structure arising from the regular $\mathtt{T}$-actions.
    A {\it two-sided ideal} of $\mathtt{T}$ is a submodule $\mathtt{J} \subset \mathtt{T}$ such that $\mathtt{T} \rhd \mathtt{J} \subset \mathtt{J}$ and $\mathtt{J} \lhd \mathtt{T} \subset \mathtt{J}$ where $\rhd,\lhd$ denote the regular $\mathtt{T}$-actions.
\end{Defn}

\begin{remark}
    Based on the equivalence (\ref{202401231608}), this notion corresponds to that of a two-sided ideal of an additive category introduced in \cite[Section 3]{Mitchell1972}.
\end{remark}

\begin{Defn} \label{202410121620}
    We define $\mathtt{T}\mbox{-}\mathsf{Mod}\mbox{-}\mathtt{S}$ as the category of $(\mathtt{T},\mathtt{S})$-bimodules and $(\mathtt{T},\mathtt{S})$-homomorphisms.
\end{Defn}

The category $\mathtt{T}\mbox{-}\mathsf{Mod}\mbox{-}\mathtt{S}$ is a bicomplete abelian category.
In particular, (co)kernel objects of $(\mathtt{T},\mathtt{S})$-homomorphisms consist of (co)kernel modules index-wise.

\begin{remark} \label{202510211626}
    Let $\mathcal{T}$ and $\mathcal{S}$ be the $\mathds{k}$-linear categories corresponding to $\mathtt{T}$ and $\mathtt{S}$ by the equivalence in (\ref{202401231608}).
    Then the discussion above yields a $\mathds{k}$-bilinear functor from the cartesian product of $\mathcal{T}$ and the opposite category of $\mathcal{S}$ to $\mathds{k}\mbox{-}\mathsf{Mod}$:
    \begin{align*}
        \mathcal{T} \times \mathcal{S}^\mathsf{o} \to \mathds{k}\mbox{-}\mathsf{Mod}; \quad (X,Y) \mapsto \mathtt{M}(X,Y) .
    \end{align*}
    Conversely, every $\mathds{k}$-bilinear functor $\mathcal{T} \times \mathcal{S}^\mathsf{o} \to \mathds{k}\mbox{-}\mathsf{Mod}$ arises in this way.
    More precisely, we have an equivalence between the category of $\mathds{k}$-bilinear functors $\mathcal{T} \times \mathcal{S}^\mathsf{o} \to \mathds{k}\mbox{-}\mathsf{Mod}$ and the category $\mathtt{T}\mbox{-}\mathsf{Mod}\mbox{-}\mathtt{S}$.
    This is established in \cite[Section 3]{betti1983variation} in a general setting.
    In particular, it gives a correspondence between $(\mathtt{T},\mathtt{S})$-homomorphisms and natural transformations.
\end{remark}

\begin{Defn} \label{202510021306}
    Let $\ast$ be a one-point set.
    Consider an object $\mathcal{X}$ of $\mathsf{Mat}_{\mathds{k}}$ and a monad $\mathtt{T}$ on $\mathcal{X}$.
    In this paper, a $(\mathtt{T}, \mathds{I}_{\ast})$-bimodule is called a {\it left $\mathtt{T}$-module}.
    Denote by $\mathtt{T}\mbox{-}\mathsf{Mod} {:=} \mathtt{T}\mbox{-}\mathsf{Mod}\mbox{-}\mathds{I}_{\ast}$.
    Similarly, a {\it right $\mathtt{T}$-module} is defined to be a $( \mathds{I}_{\ast}, \mathtt{T})$-bimodule.
    Denote by $\mathsf{Mod}\mbox{-}\mathtt{T} {:=} \mathds{I}_{\ast}\mbox{-}\mathsf{Mod}\mbox{-}\mathtt{T}$.
\end{Defn}

\begin{notation}
    By the definition, a left $\mathtt{T}$-module $\mathtt{M}$ consists of a $\mathds{k}$-module $\mathtt{M} ( X, \ast )$ parametrized by $X \in \mathcal{X}$.
    In this paper, we use the following simplified notation:
    \begin{align*}
        \mathtt{M} (X) {:=} \mathtt{M} (X, \ast) .
    \end{align*}
    In particular, the $(\mathcal{X},\ast)$-indexed module $\mathtt{M}$ is regarded as an $\mathcal{X}$-indexed module.
\end{notation}

\begin{prop} \label{202409271434}
    The category of left $\mathtt{L}_{\mathcal{C}}$-modules is equivalent to that of functors $\mathcal{C} \to \mathds{k}\mbox{-}\mathsf{Mod}$.
    In particular, it gives a correspondence between $\mathtt{L}_{\mathcal{C}}$-homomorphisms and natural transformations.
\end{prop}
\begin{proof}
    By the correspondence in Example \ref{202401271351}, Remark \ref{202510211626} implies that the category of left $\mathtt{L}_{\mathcal{C}}$-modules is equivalent to that of $\mathds{k}$-linear functors $\mathds{k}\mathcal{C} \to \mathds{k}\mbox{-}\mathsf{Mod}$.
    Moreover, the $\mathds{k}$-linearization gives a one-to-one correspondence between $\mathds{k}$-linear functors $\mathds{k}\mathcal{C} \to \mathds{k}\mbox{-}\mathsf{Mod}$ and functors $\mathcal{C} \to \mathds{k}\mbox{-}\mathsf{Mod}$.
\end{proof}

\section{Two-sided ideals, Adjunctions and Properties}
\label{202510131439}

In this section, we investigate the relationships among two-sided ideals
of a monad, adjunctions of module categories and isomorphism-invariant
properties of modules. These provide a unified way of treating properties
of modules throughout this paper. The central result is Theorem \ref{202609281358}, which lets us, in Section \ref{202404111426}, construct the
polynomiality ideal for an arbitrary monoidal category with zero unit.

The framework introduced in this section is of fundamental importance for our series of papers.
Some remarks are given in Appendix \ref{202510142110}.

\subsection{Adjunctions induced by two-sided ideals}
\label{202510121347}

In this section, we show that a two-sided ideal of a monad in $\mathsf{Mat}_{\mathds{k}}$ induces several canonical adjunctions between module categories.
This is made precise in Propositions \ref{202509281759} and \ref{202509301238}.
Let $\mathtt{T}$ be a monad in $\mathsf{Mat}_{\mathds{k}}$ on an object $\mathcal{X}$, and $\mathtt{J}$ be a two-sided ideal of $\mathtt{T}$.

\begin{Defn} \label{202510021423}
    For a left $\mathtt{T}$-module $\mathtt{M}$, we define the {\it $\mathtt{J}$-vanishing module} of $\mathtt{M}$, which we denote by $\mathrm{V}(\mathtt{M};\mathtt{J})$, to be the $\mathcal{X}$-indexed $\mathds{k}$-module
    $$
    \left( \mathrm{V}(\mathtt{M};\mathtt{J}) \right) (X) {:=} \{ x \in \mathtt{M}(X)\mid \mathtt{J} \rhd x \cong 0 \} , \quad X \in \mathcal{X}.
    $$    
\end{Defn}

One can check directly that the $\mathtt{T}$-module structure on $\mathtt{M}$ induces a left $\mathtt{T}/\mathtt{J}$-module structure on $\mathrm{V}(\mathtt{M};\mathtt{J})$.
Hence, we obtain an additive functor:
\begin{align*}
    \mathrm{V}(-;\mathtt{J}) : \mathtt{T}\mbox{-}\mathsf{Mod} \to \mathtt{T}/\mathtt{J}\mbox{-}\mathsf{Mod} .
\end{align*}

\begin{prop} \label{202509281759}
    We have an adjunction:
    $$
    \begin{tikzcd}
            \iota: \mathtt{T}/ \mathtt{J} \mbox{-}\mathsf{Mod} \arrow[r, shift right=1ex, ""{name=G}, hookrightarrow] & \mathtt{T}\mbox{-}\mathsf{Mod} : \mathrm{V}(-;\mathtt{J}) \arrow[l, shift right=1ex, ""{name=F}]
            \arrow[phantom, from=G, to=F, , "\scriptscriptstyle\boldsymbol{\top}"],
    \end{tikzcd}
    $$
    where $\iota$ is the fully faithful functor induced by the quotient monad map $\mathtt{T} \to \mathtt{T}/\mathtt{J}$.
\end{prop}
\begin{proof}
    The counit is given by the inclusion $\iota ( \mathrm{V} (\mathtt{M} ; \mathtt{J})) \to \mathtt{M}$ for $\mathtt{M} \in \mathrm{Obj} (\mathtt{T}\mbox{-}\mathsf{Mod})$.
    The unit is given by the identity $\mathtt{N} \to \mathrm{V}( \iota ( \mathtt{N});\mathtt{J}) = \mathtt{N}$.
\end{proof}

The examples of the primary interest here are given by polynomial functors which are considered in the next section.
Here, we give two easy examples that are also exploited later.
\begin{Example} \label{202509021729}
    Let $\mathcal{C}$ be a category with a zero object $\ast$.
    A functor $F : \mathcal{C} \to \mathds{k}\mbox{-}\mathsf{Mod}$ is {\it reduced} if $F(\ast) \cong 0$.
    To treat this in our framework, we introduce a two-sided ideal $\mathtt{Z}_{\mathcal{C}} \subset \mathtt{L}_{\mathcal{C}}$ generated by zero morphisms $(\ast : X \to Y) \in \mathcal{C}(X,Y) \subset \mathtt{L}_{\mathcal{C}}(Y,X)$ for $X,Y \in \mathcal{X}$.
    The adjunction in Proposition \ref{202509281759} associated with the ideal $\mathtt{Z}_{\mathcal{C}}$ gives
    $$
    \begin{tikzcd}
            \iota : \mathtt{L}_{\mathcal{C}}/ \mathtt{Z}_{\mathcal{C}} \mbox{-}\mathsf{Mod} \arrow[r, shift right=1ex, ""{name=G}, hookrightarrow] & \mathtt{L}_{\mathcal{C}}\mbox{-}\mathsf{Mod} : \mathrm{V}(-;\mathtt{Z}_{\mathcal{C}}) \arrow[l, shift right=1ex, ""{name=F}]
            \arrow[phantom, from=G, to=F, , "\scriptscriptstyle\boldsymbol{\top}"].
    \end{tikzcd}
    $$
    This yields an equivalence between $\mathtt{L}_{\mathcal{C}}/ \mathtt{Z}_{\mathcal{C}}$-modules and reduced $\mathtt{L}_{\mathcal{C}}$-modules.
\end{Example}

This example is a special case of the following:
\begin{Example} \label{202510261038}
    Let $\mathcal{C}, \mathcal{D}$ be categories with the same object set $\mathcal{X}$.
    Consider a full functor $\mathcal{C} \to \mathcal{D}$ which is the identity on objects.
    It induces a monad epimorphism $\mathtt{L}_{\mathcal{C}} \to \mathtt{L}_{\mathcal{D}}$, which we denote by $\xi$.
    Let $\mathtt{K}_{\xi} \subset \mathtt{L}_{\mathcal{C}}$ be the kernel of $\xi$.
    Clearly, it is a two-sided ideal of the monad $\mathtt{L}_{\mathcal{C}}$.
    Note that $\xi$ is an epimorphism since the functor $\mathcal{C} \to \mathcal{D}$ is full, so we have a monad isomorphism $\mathtt{L}_{\mathcal{C}}/\mathtt{K}_{\xi} \cong \mathtt{L}_{\mathcal{D}}$.
    The adjunction in Proposition \ref{202509281759} gives
    $$
    \begin{tikzcd}
            \xi^\ast : \mathtt{L}_{\mathcal{D}}\mbox{-}\mathsf{Mod} \arrow[r, shift right=1ex, ""{name=G}] & \mathtt{L}_{\mathcal{C}}\mbox{-}\mathsf{Mod} : \mathrm{V}(-;\mathtt{K}_{\xi}) \arrow[l, shift right=1ex, ""{name=F}]
            \arrow[phantom, from=G, to=F, , "\scriptscriptstyle\boldsymbol{\top}"].
    \end{tikzcd}
    $$
    This yields an equivalence between the category $\mathtt{L}_{\mathcal{D}}\mbox{-}\mathsf{Mod}$ and the full subcategory of $\mathtt{M} \in \mathtt{L}_{\mathcal{C}}\mbox{-}\mathsf{Mod}$ such that $\mathrm{V}(\mathtt{M};\mathtt{K}_{\xi}) = \mathtt{M}$.
\end{Example}

Before we close this section, we give a dual framework of the above results.

\begin{Defn} \label{202510021430}
    For a left $\mathtt{T}$-module $\mathtt{M}$, we define the {\it  $\mathtt{J}$-covanishing module} of $\mathtt{M}$ as the quotient $\mathcal{X}$-indexed $\mathds{k}$-module $\mathtt{M} / (\mathtt{J} \rhd \mathtt{M})$.
    It inherits a left $\mathtt{T}/\mathtt{J}$-module structure from the $\mathtt{T}$-module structure $\mathtt{M}$.
    We denote this by $\Lambda (\mathtt{M};\mathtt{J})$.
\end{Defn}
In analogy to Proposition \ref{202509281759}, we obtain the following.
In particular, the functor $\mathtt{T}/ \mathtt{J}\mbox{-}\mathsf{Mod} \hookrightarrow \mathtt{T}\mbox{-}\mathsf{Mod}$ admits left and right adjoints:
\begin{prop} \label{202509301238}
    We have an adjunction:
    $$
    \begin{tikzcd}
           \Lambda ( - ; \mathtt{J}) :  \mathtt{T}\mbox{-}\mathsf{Mod}\arrow[r, shift right=1ex, ""{name=G}] & \mathtt{T}/ \mathtt{J}\mbox{-}\mathsf{Mod} : \iota \arrow[l, hookrightarrow, shift right=1ex, ""{name=F}]
            \arrow[phantom, from=G, to=F, , "\scriptscriptstyle\boldsymbol{\top}"].
    \end{tikzcd}
    $$
\end{prop}

\subsection{$\mathtt{J}$-vanishingly generated modules} \label{202512241936}

In this section, for a two-sided ideal $\mathtt{J}$ of $\mathtt{T}$, we introduce the notion of a $\mathtt{J}$-vanishingly generated $\mathtt{T}$-module, and prove some basic arguments.

\begin{Defn}
    A left $\mathtt{T}$-module $\mathtt{M}$ is {\it $\mathtt{J}$-vanishingly generated} if $\mathtt{J} \rhd \mathtt{M} \cong 0$.
    Equivalently, $\mathtt{M}$ is $\mathtt{J}$-vanishingly generated if 
    \begin{align} \label{202510261307}
        \mathrm{V}(\mathtt{M};\mathtt{J})= \mathtt{M}.
    \end{align}
    Let $\mathtt{T}\mbox{-}\mathsf{Mod}^{\mathtt{J}\mathrm{-vg}} \subset \mathtt{T}\mbox{-}\mathsf{Mod}$ denote the full subcategory of $\mathtt{J}$-vanishingly generated $\mathtt{T}$-modules.
    Here, $\mathtt{J}\mathrm{-vg}$ stands for $\mathtt{J}$-vanishingly generated.
\end{Defn}

Clearly, the adjunction given in Proposition \ref{202509281759} gives an equivalence of categories:
\begin{align} \label{202512241357}
    \mathtt{T}/\mathtt{J}\mbox{-}\mathsf{Mod} \simeq \mathtt{T}\mbox{-}\mathsf{Mod}^{\mathtt{J}\mathrm{-vg}} .
\end{align}
In particular, the category $\mathtt{T}\mbox{-}\mathsf{Mod}^{\mathtt{J}\mathrm{-vg}}$ is an abelian category admitting arbitrary limits and colimits.

\begin{Defn} \label{202510121351}
    Let $\mathtt{T}\mbox{-}\mathcal{M}\mathsf{od}$ denote the class of isomorphism classes of left $\mathtt{T}$-modules.
    We denote by $\mathcal{V}_{\mathtt{T}}(\mathtt{J})$, or simply by $\mathcal{V} (\mathtt{J})$, the subclass of $\mathtt{T}\mbox{-}\mathcal{M}\mathsf{od}$ consisting of isomorphism classes of $\mathtt{J}$-vanishingly generated $\mathtt{T}$-modules.
\end{Defn}

\begin{Lemma} \label{202509291501}
    Let $\mathtt{T},\mathtt{S}$ be monads on an object $\mathcal{X}$ of $\mathsf{Mat}_{\mathds{k}}$, and $\pi : \mathtt{T} \to \mathtt{S}$ be a monad epimorphism.
    Let $\pi^\ast : \mathtt{S}\mbox{-}\mathcal{M}\mathsf{od} \to \mathtt{T}\mbox{-}\mathcal{M}\mathsf{od}$ be the class map induced by $\pi$.
    \begin{enumerate}
        \item The map $\pi^\ast$ is injective.
        \item
        For two-sided ideals $\mathtt{J} \subset \mathtt{T}$ and $\mathtt{I} \subset \mathtt{S}$ such that $\pi (\mathtt{J}) \subset \mathtt{I}$, we have $\pi^\ast ( \mathcal{V}_{\mathtt{S}}(\mathtt{I}) ) \subset \mathcal{V}_{\mathtt{T}}(\mathtt{J})$.
        \item 
        For a two-sided ideal $\mathtt{J} \subset \mathtt{T}$ with $\mathtt{I} = \pi (\mathtt{J})$, we have $\mathcal{V}_{\mathtt{S}} (\mathtt{I}) = (\pi^\ast)^{-1} (\mathcal{V}_{\mathtt{T}} (\mathtt{J}))$. 
        Moreover, $\pi^\ast ( \mathcal{V}_{\mathtt{S}}(\mathtt{I}) ) = \mathcal{V}_{\mathtt{T}}(\mathtt{J})$ if and only if $\pi$ induces an isomorphism $\mathtt{T}/\mathtt{J} \to \mathtt{S}/\mathtt{I}$.
    \end{enumerate}
\end{Lemma}
\begin{proof}
    We start by proving that $\pi^\ast$ is injective.
    Let $\mathtt{M},\mathtt{N}$ be $\mathtt{S}$-modules with a $\mathtt{T}$-isomorphism $h : \pi^\ast (\mathtt{M}) \to \pi^\ast (\mathtt{N})$.
    It satisfies $\pi (f) \rhd h(x) = h ( \pi (f) \rhd x)$ for $x \in \mathtt{M} (X)$ and $f \in \mathtt{T} (Y,X)$ where $X,Y\in\mathcal{X}$.
    Since $\pi$ is an epimorphism, this condition implies that $h$ gives an $\mathtt{S}$-isomorphism.
    
    (2) is proved as follows:
    Let $\mathtt{M}$ be an $\mathtt{I}$-vanishingly generated $\mathtt{S}$-module.
    Then $\mathtt{J} \rhd \pi^\ast (\mathtt{M}) = \pi ( \mathtt{J}) \rhd \mathtt{M} \subset \mathtt{I} \rhd \mathtt{M} = 0$, hence $\pi^\ast (\mathtt{M})$ is $\mathtt{J}$-vanishingly generated.
    
    Next, we prove the first statement of (3).
    Note that for a left $\mathtt{S}$-module $\mathtt{M}$, $\mathtt{I} \rhd \mathtt{M}  = \pi ( \mathtt{J}) \rhd \mathtt{M} = \mathtt{J} \rhd \pi^\ast (\mathtt{M})$.
    Hence, $\mathtt{M}$ is $\mathtt{I}$-vanishingly generated if and only if $\pi^\ast (\mathtt{M})$ is $\mathtt{J}$-vanishingly generated.

    We now prove the rest of (3).
    By the bijection $\mathcal{V}_{\mathtt{T}}(\mathtt{J}) \cong \mathtt{T}/\mathtt{J}\mbox{-}\mathcal{M}\mathsf{od}$ (see \eqref{202512241357}), if the map $\mathtt{T}/\mathtt{J} \to \mathtt{S}/\mathtt{I}$ is an isomorphism, then $\pi^\ast ( \mathcal{V}_{\mathtt{S}}(\mathtt{I}) ) = \mathcal{V}_{\mathtt{T}}(\mathtt{J})$.
    Conversely, we assume that $\pi^\ast ( \mathcal{V}_{\mathtt{S}}(\mathtt{I}) ) = \mathcal{V}_{\mathtt{T}}(\mathtt{J})$.
    It suffices to show that $\mathtt{T}/\mathtt{J} \to \mathtt{S}/\mathtt{I}$ is a monomorphism, since it is an epimorphism by the assumption.
    Since $\pi(\mathtt{J}) = \mathtt{I}$, one may check that the kernel of $\mathtt{T}/\mathtt{J} \to \mathtt{S}/\mathtt{I}$ is $(\ker (\pi) + \mathtt{J} )/\mathtt{J} \subset \mathtt{T}/\mathtt{J}$.
    On the other hand, by the hypothesis, for $X \in \mathcal{X}$, there exists a unique $[\mathtt{M}_X] \in \mathcal{V}_{\mathtt{S}}(\mathtt{I})$ such that $\pi^\ast (\mathtt{M}_X) \cong (\mathtt{T}/\mathtt{J}) (-,X)$ as $\mathtt{T}$-modules, since $(\mathtt{T}/\mathtt{J})(-,X)$ is $\mathtt{J}$-vanishingly generated.
    Therefore, we have
    $$\ker ( \pi) \rhd (\mathtt{T}/\mathtt{J})(-,X) = \ker ( \pi) \rhd  \pi^\ast (\mathtt{M}_X) = \pi^\ast \left(  \pi ( \ker ( \pi) ) \rhd \mathtt{M}_X \right) \cong 0 .$$
    In particular, $\ker ( \pi) \rhd 1_X \subset \mathtt{J} (-,X)$.
    By the equality $(\ker (\pi)) (-,X) = \ker ( \pi) \rhd 1_X$, we obtain $(\ker (\pi)) (-,X)  \subset \mathtt{J} (-,X)$.
    This holds for any $X \in \mathcal{X}$, so we obtain $\ker (\pi) \subset \mathtt{J}$.
    Hence, the kernel of $\mathtt{T}/\mathtt{J} \to \mathtt{S}/\mathtt{I}$ is trivial.
\end{proof}

\subsection{Ideals from module properties}
\label{202510121348}

In this section, we construct a retract of the map $\mathcal{V}$ introduced in Definition \ref{202510121351}.
The motivation is to find a two-sided ideal that induces a given subclass of $\mathtt{T}\mbox{-}\mathcal{M}\mathsf{od}$, in the sense of Definition \ref{202512281236} below.
The framework developed in the present section will be used in one of the fundamental theorems of this paper, namely Theorem \ref{202508250938}.
\begin{Defn} \label{202512281236}
    Let $\mathcal{S}$ be a subclass of $\mathtt{T}\mbox{-}\mathcal{M}\mathsf{od}$.
    A two-sided ideal $\mathtt{J}$ of $\mathtt{T}$ is a {\it core $\mathtt{T}$-internalizer} of $\mathcal{S}$ if we have $\mathcal{S} = \mathcal{V} (\mathtt{J})$.
\end{Defn}
Given a subclass $\mathcal{S}$, if one can find a core $\mathtt{T}$-internalizer, then it yields a natural adjunction as in Propositions \ref{202509281759} and \ref{202509301238}.
This respects the class $\mathcal{S}$ in the sense of the equivalence (\ref{202512241357}).

We now introduce a method to construct a left ideal from a given subclass $\mathcal{S} \subset \mathtt{T}\mbox{-}\mathcal{M}\mathsf{od}$, which provides a candidate for a core $\mathtt{T}$-internalizer.

\begin{Defn}
    Let $\mathtt{M}$ be a left $\mathtt{T}$-module.
    The {\it annihilator} of $\mathtt{M}$, denoted by $\mathrm{Ann} (\mathtt{M})$, is defined to be the two-sided ideal $\mathrm{Ann} (\mathtt{M}) \subset \mathtt{T}$,
    $$
    \left( \mathrm{Ann} (\mathtt{M}) \right) (Y,X) {:=} \{ f \in \mathtt{T} (Y,X) ~|~ f \rhd \mathtt{M} (X) \cong 0  \}, \quad X,Y \in \mathcal{X} .
    $$
\end{Defn}

\begin{Defn} \label{202409271755}
    For a subclass $\mathcal{S} \subset \mathtt{T}\mbox{-}\mathcal{M}\mathsf{od}$, we define the two-sided ideal $\mathrm{Ann} (\mathcal{S})$ of $\mathtt{T}$ to be the intersection of the annihilator of $\mathtt{M}$ for every $\mathtt{T}$-module $\mathtt{M}$ whose isomorphism class lies in $\mathcal{S}$:
    \begin{align*}
        \mathrm{Ann} ( \mathcal{S} ) {:=} \bigcap_{[\mathtt{M}] \in \mathcal{S}} \mathrm{Ann} (\mathtt{M}) .
    \end{align*}
\end{Defn}

\begin{prop} \label{202510142114}
    Let $\mathtt{J}$ be a two-sided ideal of $\mathtt{T}$.
    We have $\mathrm{Ann}(\mathcal{V}(\mathtt{J})) = \mathtt{J}$.
\end{prop}
\begin{proof}
    By the definitions, we have 
    $$\mathrm{Ann} ( \mathcal{V} (\mathtt{J}) ) \subset \bigcap_{X\in\mathcal{X}} \mathrm{Ann} ( (\mathtt{T}/ \mathtt{J}) (-,X)) \subset \mathtt{J} .$$
    On the other hand, for a $\mathtt{J}$-vanishingly generated $\mathtt{M}$, we have $\mathrm{V} (\mathtt{M} ; \mathtt{J}) = \mathtt{M}$, so $\mathrm{Ann} ( \mathtt{M})$ contains $\mathtt{J}$.
    Thus, we obtain $\mathtt{J} \subset \mathrm{Ann} ( \mathcal{V} ( \mathtt{J})  )$ which leads to the desired statement.
\end{proof}

In the following, we give a characterization of a core internalizer.
\begin{Corollary} \label{202409271802} 
	If a subclass $\mathcal{S} \subset \mathtt{T}\mbox{-}\mathcal{M}\mathsf{od}$ admits a core $\mathtt{T}$-internalizer, then it is $\mathrm{Ann} ( \mathcal{S} )$.
    In particular, $\mathcal{S}$ has at most one core $\mathtt{T}$-internalizer.
\end{Corollary}
\begin{proof}
	Let $\mathtt{J}$ be a core internalizer of $\mathcal{S}$.
	By Proposition \ref{202510142114}, we have $\mathrm{Ann} ( \mathcal{S}) = \mathrm{Ann}(\mathcal{V}(\mathtt{J})) = \mathtt{J}$.
\end{proof}

\begin{Corollary} \label{202509291354}
    Let $\mathtt{I},\mathtt{J}$ be two-sided ideals of $\mathtt{T}$.
    The following conditions are equivalent:
    \begin{enumerate}
        \item $\mathcal{V}(\mathtt{J}) \subset \mathcal{V} (\mathtt{I})$.
        \item $\mathtt{I} \subset \mathtt{J}$.
        \item $[(\mathtt{T}/\mathtt{J}) (-,X)] \in \mathcal{V} (\mathtt{I})$ for $X \in \mathcal{X}$.
    \end{enumerate}
\end{Corollary}
\begin{proof}
    We first prove (1) starting from (2).
    For a $\mathtt{J}$-vanishingly generated $\mathtt{T}$-module $\mathtt{M}$, $\mathtt{I}\rhd \mathtt{M} \subset \mathtt{J} \rhd \mathtt{M} = 0$, so $\mathcal{V}(\mathtt{J}) \subset \mathcal{V} (\mathtt{I})$.
    Conversely, if $\mathcal{V}(\mathtt{J}) \subset \mathcal{V} (\mathtt{I})$, then by Proposition \ref{202510142114} $\mathtt{I} = \mathrm{Ann} (\mathcal{V} (\mathtt{I})) \subset \mathrm{Ann} (\mathcal{V} (\mathtt{J})) = \mathtt{J}$.
    Therefore, (1) and (2) are equivalent.
    
    We now show that (3) implies (2).
    By Proposition \ref{202510142114}, for $X \in \mathcal{X}$, we have
    $$
    \mathtt{I} (-,X) = \left( \mathrm{Ann} (\mathcal{V}(\mathtt{I})) \right) (-,X)\subset \left( \mathrm{Ann} ((\mathtt{T}/\mathtt{J})(-,X)) \right) (-,X) \subset \mathtt{J}(-,X) .
    $$
    Conversely, it is straightforward to show that (1) implies (3).
\end{proof}

In the rest of this section, we give an equivalent condition for $\mathcal{S}\subset\mathtt{T}\mbox{-}\mathcal{M}\mathsf{od}$ to have a core $\mathtt{T}$-internalizer.

\begin{Defn} \label{202609281518}
    For a subclass $\mathcal{S}\subset\mathtt{T}\mbox{-}\mathcal{M}\mathsf{od}$, let $\tilde{\mathcal{S}}$ denote the full subcategory of $\mathtt{T}\mbox{-}\mathsf{Mod}$ whose objects are $\mathtt{M}$ such that $[\mathtt{M}] \in \mathcal{S}$.
    A subclass $\mathcal{S}\subset\mathtt{T}\mbox{-}\mathcal{M}\mathsf{od}$ is {\it closed} if the full subcategory $\tilde{\mathcal{S}}\subset\mathtt{T}\mbox{-}\mathsf{Mod}$ is closed under subobjects, quotient objects, arbitrary products and arbitrary direct sums, all formed in $\mathtt{T}\mbox{-}\mathsf{Mod}$.
\end{Defn}

Since kernels, cokernels, limits and colimits in $\mathtt{T}\mbox{-}\mathsf{Mod}$ are computed index-wise, a closed class is in particular closed under all limits and colimits.

\begin{Lemma}\label{202609281400}
    For every two-sided ideal $\mathtt{J}\subset\mathtt{T}$, the class $\mathcal{V}(\mathtt{J})$ is closed.
\end{Lemma}
\begin{proof}
    The condition $\mathtt{J}\rhd\mathtt{M}=0$ is inherited by submodules and quotient modules, and it is checked index-wise, so it is inherited by direct products and direct sums.
\end{proof}

\begin{Lemma}\label{202609281319}
    Let $\mathcal{S}$ be closed and set $\mathtt{K} = \mathrm{Ann}(\mathcal{S})$.
    For $X\in\mathcal{X}$, $\mathtt{T}(-,X)/\mathtt{K}(-,X)$ belongs to $\tilde{\mathcal{S}}$.
\end{Lemma}
\begin{proof}
    For a left $\mathtt{T}$-module $\mathtt{M}$ satisfying $[\mathtt{M}] \in \mathcal{S}$ and $v \in \mathtt{M}(X)$, let $$\mathtt{E}_{\mathtt{M},v} (Y,X) {:=} \{ f \in \mathtt{T}(Y,X) \mid f \rhd v = 0 \} , \quad Y \in \mathcal{X}.$$
    Since $\mathtt{T}(-,X) \to \mathtt{M}; f \mapsto f \rhd v$ factors into a monomorphism $(\mathtt{T}/\mathtt{E}_{\mathtt{M},v})(-,X) \to \mathtt{M}$ whose target belongs to $\tilde{\mathcal{S}}$, we obtain $(\mathtt{T}/\mathtt{E}_{\mathtt{M},v})(-,X) \in \tilde{\mathcal{S}}$.
    Let $\mathcal{A}$ be the set of submodules of $\mathtt{T}(-,X)$ of the form $\mathtt{E}_{\mathtt{M},v}(-,X)$, where $[\mathtt{M}] \in \mathcal{S}$ and $v \in \mathtt{M}(X)$.
    By the argument above, $(\mathtt{T}/\mathtt{E})(-,X) \in \tilde{\mathcal{S}}$ for every $\mathtt{E} \in \mathcal{A}$.
    By definition,
    $$
    \mathrm{Ann}(\mathcal{S})(-,X) = \bigcap_{\mathtt{E} \in \mathcal{A}} \mathtt{E} .
    $$
    The canonical map
    $$
    \mathtt{T}(-,X) \Big/ \bigcap_{\mathtt{E} \in \mathcal{A}} \mathtt{E} \longrightarrow \prod_{\mathtt{E} \in \mathcal{A}} \mathtt{T}(-,X)/\mathtt{E}
    $$
    is a monomorphism, and its target lies in $\tilde{\mathcal{S}}$ since $\mathcal{S}$ is closed under products. 
    Hence its source $(\mathtt{T}/\mathrm{Ann}(\mathcal{S}))(-,X)$ lies in $\tilde{\mathcal{S}}$
    by closure under subobjects.
\end{proof}

\begin{theorem}\label{202609281358}
    Let $\mathcal{S}\subset\mathtt{T}\mbox{-}\mathcal{M}\mathsf{od}$ be a subclass.
    The following are equivalent:
    \begin{enumerate}
    \item $\mathcal{S}$ admits a core $\mathtt{T}$-internalizer;
    \item $\mathcal{S}$ is closed.
    \end{enumerate}
    In this case, $\mathrm{Ann}(\mathcal{S})$ is the unique core $\mathtt{T}$-internalizer of $\mathcal{S}$.
\end{theorem}
\begin{proof}
    (1)$\Rightarrow$(2) is Lemma \ref{202609281400}.
    
    (2)$\Rightarrow$(1): Put $\mathtt{K}:=\mathrm{Ann}(\mathcal{S})$, a two-sided ideal of $\mathtt{T}$.
    It suffices to show that $\mathcal{V}(\mathtt{K})=\mathcal{S}$.
    The inclusion $\mathcal{S}\subset\mathcal{V}(\mathtt{K})$ is immediate, since
    $\mathtt{K}\subset\mathrm{Ann}(\mathtt{M})$ whenever $[\mathtt{M}]\in\mathcal{S}$.
    To prove the converse, let $\mathtt{M}$ be a left $\mathtt{T}$-module with
    $\mathtt{K}\rhd\mathtt{M}=0$.
    Then each component of the canonical epimorphism
    $$
    \bigoplus_{X\in\mathcal{X}}\ \bigoplus_{v\in\mathtt{M}(X)}\mathtt{T}(-,X)\longrightarrow\mathtt{M},\qquad g\longmapsto g\rhd v\ \text{ on the $(X,v)$-summand},
    $$
    vanishes on $\mathtt{K}(-,X)$.
    Hence the epimorphism factors through $\bigoplus_{X\in\mathcal{X}}\ \bigoplus_{v\in\mathtt{M}(X)}\mathtt{T}(-,X)/\mathtt{K}(-,X)$, which lies in $\tilde{\mathcal{S}}$ by Lemma \ref{202609281319} and closure under direct sums.
    Thus $\mathtt{M}$ is a quotient of an object of $\tilde{\mathcal{S}}$, so $\mathtt{M}\in\tilde{\mathcal{S}}$.
\end{proof}

\subsection{Filtration of internalizers}
\label{202510121349}

Let $\mathtt{T}$ be a monad in $\mathsf{Mat}_{\mathds{k}}$.
For a descending filtration by two-sided ideals of $\mathtt{T}$
$$\mathtt{J}^{(\bullet)} : \quad  \cdots \subset \mathtt{J}^{(d+1)} \subset \mathtt{J}^{(d)} \subset \cdots \subset \mathtt{J}^{(0)} \subset \mathtt{T} ,$$
let $\mathcal{S}_{\infty} \subset \mathtt{T}\mbox{-}\mathcal{M}\mathsf{od}$ be the subclass of isomorphism classes of the left $\mathtt{T}$-module $\mathtt{M}$ which satisfy $$\varinjlim_{d} \mathrm{V} ( \mathtt{M} ; \mathtt{J}^{(d)}) \cong \mathtt{M}.$$
In this section, we give some conditions that ensure that $\mathcal{S}_{\infty}$ admits a core internalizer.
The results given here will be applied in Section \ref{202509301437} to study analytic functors.

\begin{Defn} \label{202508270931}
    The filtration $\mathtt{J}^{(\bullet)}$ {\it left-stabilizes} if, for every $X \in \mathcal{X}$, there exists $d_0 \in \mathds{N}$ such that the following filtration of $\mathtt{T}(-,X)$ stabilizes after $d_0$:
    $$\cdots = \mathtt{J}^{(d_0+1)} (-,X) = \mathtt{J}^{(d_0)}(-,X) \subset \mathtt{T} (-,X).$$ 
\end{Defn}

\begin{remark}
    If a filtration $\mathtt{J}^{(\bullet)}$ does not left-stabilize, then, in particular, it does not stabilize: there does not exist $d_0 \in \mathds{N}$ such that $\mathtt{J}^{(d_0)} = \mathtt{J}^{(d_0+1)} =\cdots$.
\end{remark}

\begin{prop} \label{202405291746}
    The following conditions are equivalent:
    \begin{enumerate}
        \item 
        There exists a core $\mathtt{T}$-internalizer of $\mathcal{S}_{\infty}$.
        \item
        The filtration $\mathtt{J}^{(\bullet)}$ left-stabilizes.
        \item $\mathcal{S}_{\infty}$ is closed.
    \end{enumerate}
    In this case, the core $\mathtt{T}$-internalizer of $\mathcal{S}_{\infty}$, say $\mathtt{J}$, is given by $\mathtt{J}(-,X)= \mathtt{J}^{(d_0)}(-,X)$ where $X \in \mathcal{X}$ and $d_0$ is given in Definition \ref{202508270931}. 
\end{prop}
\begin{proof}
    (1) $\Leftrightarrow$ (3) follows from Theorem \ref{202609281358}.

    Next, we prove (3) $\Rightarrow$ (2).
    Let $\tilde{\mathcal{S}}_{\infty}$ denote the full subcategory of $\mathtt{T}\mbox{-}\mathsf{Mod}$ whose objects are $\mathtt{M}$ such that $[\mathtt{M}] \in \mathcal{S}_{\infty}$.
    Fix $X \in \mathcal{X}$.
    By definition of $\mathcal{S}_{\infty}$, $(\mathtt{T}/\mathtt{J}^{(d)}) (-,X) \in \tilde{\mathcal{S}}_{\infty}$, where $d \in \mathds{N}$.
    Since $\mathcal{S}_{\infty}$ is closed, $\mathtt{P} = \prod_{d \in \mathds{N}} (\mathtt{T}/\mathtt{J}^{(d)}) (-,X) \in \tilde{\mathcal{S}}_{\infty}$.
    Consider the element $f = \left( f_d \right)_{d\in\mathds{N}} \in \mathtt{P}(X) = \prod_{d \in \mathds{N}} (\mathtt{T}/\mathtt{J}^{(d)}) (X,X)$ given by $f_d = 1_X$ for every $d \in \mathds{N}$.
    There exists $d_0 \in \mathds{N}$ such that $\mathtt{J}^{(d_0)} \rhd f =0$, since $\varinjlim_{d} \mathrm{V} ( \mathtt{P} ; \mathtt{J}^{(d)}) \cong \mathtt{P}$.
    This condition is equivalent to $\mathtt{J}^{(d_0)} (-,X)= \mathtt{J}^{(d)} (-,X)$ where $d \geq d_0$.
    Thus, the proof of (2) is complete.

    We now prove (1) from (2).
    Let $\mathtt{J} \subset \mathtt{T}$ be a submodule defined by $\mathtt{J}(-,X){:=} \mathtt{J}^{(d_0)}(-,X)$ for $X$ and $d_0$ in Definition \ref{202508270931}.
    The condition (2) implies that $\mathtt{J} = \bigcap_{d \in \mathds{N}} \mathtt{J}^{(d)}$, so that $\mathtt{J}$ is a two-sided ideal of $\mathtt{T}$, since so $\mathtt{J}^{(d)}$ is.
    It suffices to prove that $\mathcal{S}_{\infty} = \mathcal{V}(\mathtt{J})$.
    If a left $\mathtt{T}$-module $\mathtt{M}$ satisfies $\varinjlim_{d} \mathrm{V} ( \mathtt{M} ;\mathtt{J}^{(d)}) \cong \mathtt{M}$, then, for $X \in \mathcal{X}$ and $v \in \mathtt{M}(X)$, there exists $d \in \mathds{N}$ such that $\mathtt{J}^{(d)} \rhd v = 0$.
    Hence, we obtain $\mathtt{J} \rhd v = 0$ by $\mathtt{J} = \bigcap_{d \in \mathds{N}} \mathtt{J}^{(d)}$.
    This implies $\mathrm{V} ( \mathtt{M} ; \mathtt{J}) = \mathtt{M}$, so $\mathcal{S}_{\infty} \subset \mathcal{V}(\mathtt{J})$.
    On the other hand, if $\mathtt{M}$ is $\mathtt{J}$-vanishingly generated, i.e. $\mathrm{V} ( \mathtt{M} ; \mathtt{J}) = \mathtt{M}$, then one can show $\varinjlim_{d} \mathrm{V} ( \mathtt{M} ;\mathtt{J}^{(d)}) \cong \mathtt{M}$ from the definition of $\mathtt{J}$.
    This implies $\mathcal{V}(\mathtt{J}) \subset \mathcal{S}_{\infty}$.
\end{proof}

\section{Polynomiality ideal}
\label{202408041517}
Applying the framework of Section \ref{202510131439}, we introduce a two-sided ideal of $\mathtt{L}_{\mathcal{C}}$ which encodes the property of being of finite polynomial degree, for a small monoidal category $\mathcal{C}$ with zero unit (recall the monad $\mathtt{L}_{\mathcal{C}}$ from Example \ref{202401271351}).
Using this, we also study the sets $\Gamma(\xi)$ and $\mathrm{D}(\mathcal{C})$ introduced in the introduction.
When $\mathcal{C}$ has finite products, or finite coproducts, we give in Sections \ref{202401241445} and \ref{202508152158} respectively a useful description of this ideal, which is need in Section \ref{202510161122}.

\subsection{The polynomiality ideal}
\label{202404111426}

In this section, we revisit the theory of polynomial functors using the framework of Section \ref{202510131439}.
We show, in Theorem \ref{202508250938}, that the class of isomorphism classes of $\mathtt{L}_{\mathcal{C}}$-modules of degree at most $d$
admits a unique core $\mathtt{L}_{\mathcal{C}}$-internalizer, called the {\it polynomiality ideal} $\mathtt{I}^{(d)}_{\mathcal{C}}$,
and use it to prove Theorem \ref{202510011337}.
We also give an explicit form of $\mathtt{I}^{(d)}_{\mathcal{C}}$ (Proposition \ref{202609282123}).

Let $\mathcal{C}$ be a small monoidal category with zero unit.
Denote by $\mathcal{X}$ the object set of $\mathcal{C}$.
Recall the monad $\mathtt{L}_{\mathcal{C}}$ in $\mathsf{Mat}_{\mathds{k}}$ on $\mathcal{X}$ defined in Example \ref{202401271351}.
By Proposition \ref{202409271434}, we identify functors $\mathcal{C} \to \mathds{k}\mbox{-}\mathsf{Mod}$ with left $\mathtt{L}_{\mathcal{C}}$-modules in $\mathsf{Mat}_{\mathds{k}}$.
A left $\mathtt{L}_{\mathcal{C}}$-module $\mathtt{M}$ is {\it of polynomial degree at most $d$}, i.e. $\deg \mathtt{M} \leq d$, if so is the corresponding $\mathcal{C}$-module.

\begin{Defn} \label{202512241512}
    Let $\mathtt{L}_{\mathcal{C}}\mbox{-} \mathsf{Mod}^{\leq d}$ be the full subcategory of $\mathtt{L}_{\mathcal{C}}\mbox{-}\mathsf{Mod}$ consisting of $\mathtt{L}_{\mathcal{C}}$-modules of $\deg \leq d$.
    Let $\mathcal{S}^{d}_{\mathcal{C}}$ be the class of isomorphism classes of $\mathtt{L}_{\mathcal{C}}$-modules of $\deg \leq d$.
\end{Defn}

Recall that $\mathtt{L}_{\mathcal{C}}\mbox{-}\mathcal{M}\mathsf{od}$ denotes the class of isomorphism classes of left $\mathtt{L}_{\mathcal{C}}$-modules.
By Proposition \ref{202609281508}, we have an increasing filtration of $\mathtt{L}_{\mathcal{C}}\mbox{-}\mathcal{M}\mathsf{od}$:
\begin{align} \label{202509291620}
    \mathcal{S}^{0}_{\mathcal{C}} \subset \cdots \subset \mathcal{S}^{d}_{\mathcal{C}} \subset \mathcal{S}^{d+1}_{\mathcal{C}} \subset \cdots \subset \mathtt{L}_{\mathcal{C}}\mbox{-}\mathcal{M}\mathsf{od} .
\end{align}

Recall from Section \ref{202509301224} that the $d$-th polynomial approximation is right adjoint to the inclusion.
With a slight abuse of notation, let $P_d : \mathtt{L}_{\mathcal{C}}\mbox{-}\mathsf{Mod} \to \mathtt{L}_{\mathcal{C}}\mbox{-}\mathsf{Mod}^{\leq d}$ be the induced functor, so that we have an adjunction:
\begin{equation} \label{202407041637}
\begin{tikzcd}
         \mathtt{L}_{\mathcal{C}}\mbox{-} \mathsf{Mod}^{\leq d}  \arrow[r, shift right=1ex, ""{name=G}, hookrightarrow] & \mathtt{L}_{\mathcal{C}}\mbox{-} \mathsf{Mod} : P_{d} \arrow[l, shift right=1ex, ""{name=F}]
        \arrow[phantom, from=G, to=F, , "\scriptscriptstyle\boldsymbol{\top}"].
\end{tikzcd}
\end{equation}

The following is the main theorem of this section.

\begin{theorem} \label{202508250938}
    For $d \in \mathds{N}$, the class $\mathcal{S}^{d}_{\mathcal{C}}$ admits a unique core $\mathtt{L}_{\mathcal{C}}$-internalizer, given by $\mathtt{I}^{(d)}_{\mathcal{C}} {:=} \mathrm{Ann} ( \mathcal{S}^{d}_{\mathcal{C}} )$.
    We call $\mathtt{I}^{(d)}_{\mathcal{C}}$ the $d$-th polynomiality ideal for $\mathcal{C}$.
    Furthermore, this satisfies the following:
    \begin{enumerate}
        \item The adjunction associated with $\mathtt{I}^{(d)}_{\mathcal{C}}$, given by Proposition \ref{202509281759}, induces
        $$
        \mathtt{L}_{\mathcal{C}} /\mathtt{I}^{(d)}_{\mathcal{C}}\mbox{-} \mathsf{Mod} \simeq \mathtt{L}_{\mathcal{C}}\mbox{-} \mathsf{Mod}^{\leq d}  \quad \mathrm{with} \quad P_d \cong \mathrm{V}(-;\mathtt{I}^{(d)}_{\mathcal{C}}) .
        $$
        \item The two-sided ideals $\mathtt{I}^{(d)}_{\mathcal{C}}$ for $d \in \mathds{N}$ form a descending filtration of $\mathtt{L}_{\mathcal{C}}$:
        \begin{align*}
            \cdots \subset \mathtt{I}^{(d+1)}_{\mathcal{C}} \subset \mathtt{I}^{(d)}_{\mathcal{C}} \subset \cdots \subset \mathtt{L}_{\mathcal{C}} .
        \end{align*}
        \item Let $\mathcal{D}$ be another monoidal category with zero unit.
        For a full strong monoidal functor $\xi : \mathcal{C} \to \mathcal{D}$ that is the identity on objects, we have $$\xi ( \mathtt{I}^{(d)}_{\mathcal{C}}) = \mathtt{I}^{(d)}_{\mathcal{D}}.$$
    \end{enumerate}
\end{theorem}

To prove this theorem, we need the following lemma.
\begin{Lemma} \label{202609281638}
    Let $\mathcal{D}$ be another monoidal category with zero unit.
    For a strong monoidal functor $\xi : \mathcal{C} \to \mathcal{D}$ that is the identity on objects, we have $(\xi^\ast)^{-1}(\mathcal{S}^d_{\mathcal{C}})=\mathcal{S}^d_{\mathcal{D}}$ where $\xi^\ast:\mathtt{L}_{\mathcal{D}}\mbox{-}\mathcal{M}\mathsf{od}\to\mathtt{L}_{\mathcal{C}}\mbox{-}\mathcal{M}\mathsf{od}$ is the induced class map.
\end{Lemma}
\begin{proof}
    Let $F$ be a $\mathcal{D}$-module.
    Recall the morphisms $\epsilon_{X_1,\cdots,X_{d+1};k}$ and $\eta_{X_1,\cdots,X_{d+1};k}$ from Section \ref{202509301224}.
    Since $\xi$ is strong monoidal and the identity on objects, the naturality of its monoidal structure maps
    identifies $\xi(\epsilon_{X_1,\cdots,X_{d+1};k})$ and $\xi(\eta_{X_1,\cdots,X_{d+1};k})$ with the corresponding
    morphisms of $\mathcal{D}$.
    Hence there is an isomorphism
    $\mathrm{cr}_{d+1}(\xi^\ast F)(X_1,\cdots,X_{d+1})\cong\mathrm{cr}_{d+1}(F)(X_1,\cdots,X_{d+1})$
    for all objects $X_1,\cdots,X_{d+1}$.
    Thus $\deg F\le d$ if and only if $\deg\xi^\ast F\le d$, that is,
    $(\xi^\ast)^{-1}(\mathcal{S}^d_{\mathcal{C}})=\mathcal{S}^d_{\mathcal{D}}$.
\end{proof}

\begin{proof}[Proof of Theorem \ref{202508250938}]
    By Lemma \ref{202609281530}, the subclass $\mathcal{S}^{d}_{\mathcal{C}} \subset \mathtt{L}_{\mathcal{C}}\mbox{-}\mathcal{M}\mathsf{od}$ is closed in the sense of Definition \ref{202609281518}.
    Theorem \ref{202609281358} shows that there exists a unique core $\mathtt{L}_{\mathcal{C}}$-internalizer of $\mathcal{S}^{d}_{\mathcal{C}}$, given by  $\mathtt{I}^{(d)}_{\mathcal{C}} = \mathrm{Ann} ( \mathcal{S}^{d}_{\mathcal{C}} )$, i.e. $\mathcal{S}^{d}_{\mathcal{C}} = \mathcal{V} ( \mathtt{I}^{(d)}_{\mathcal{C}} )$.
    Therefore, a left $\mathtt{L}_{\mathcal{C}}$-module is of degree at most $d$ if and only if it is $\mathtt{I}^{(d)}_{\mathcal{C}}$-vanishingly generated.
    Hence, the adjunction of Proposition \ref{202509281759} applied to $\mathtt{J} = \mathtt{I}^{(d)}_{\mathcal{C}}$ yields the equivalence of (1).
    Under this equivalence, the adjunction \ref{202407041637} coincides with the adjunction associated with $\mathtt{I}^{(d)}_{\mathcal{C}}$, given by Proposition \ref{202509281759}.
    By uniqueness of right adjoint functors, we obtain $ P_d \cong \mathrm{V}(-;\mathtt{I}^{(d)}_{\mathcal{C}})$.
    
    (2) is proved by applying Corollary \ref{202509291354} to the filtration \eqref{202509291620}.

    We now prove (3).
    Since $\xi$ is full, it induces a monad epimorphism $\xi:\mathtt{L}_{\mathcal{C}}\to\mathtt{L}_{\mathcal{D}}$;
    in particular, $\xi(\mathtt{I}^{(d)}_{\mathcal{C}})$ is a two-sided ideal of $\mathtt{L}_{\mathcal{D}}$,
    since $\xi$ is an epimorphism.
    By Corollary \ref{202509291354}, it suffices to prove
    $\mathcal{V}(\xi(\mathtt{I}^{(d)}_{\mathcal{C}}))=\mathcal{V}(\mathtt{I}^{(d)}_{\mathcal{D}})=\mathcal{S}^d_{\mathcal{D}}$.
    By (3) of Lemma \ref{202509291501},
    $$
    \mathcal{V}(\xi(\mathtt{I}^{(d)}_{\mathcal{C}}))
    =(\xi^\ast)^{-1}(\mathcal{V}(\mathtt{I}^{(d)}_{\mathcal{C}}))
    =(\xi^\ast)^{-1}(\mathcal{S}^d_{\mathcal{C}}),
    $$
    where $\xi^\ast:\mathtt{L}_{\mathcal{D}}\mbox{-}\mathcal{M}\mathsf{od}\to\mathtt{L}_{\mathcal{C}}\mbox{-}\mathcal{M}\mathsf{od}$
    is the induced class map.
    By Lemma \ref{202609281638}, $(\xi^\ast)^{-1}(\mathcal{S}^d_{\mathcal{C}})=\mathcal{S}^d_{\mathcal{D}}$.
\end{proof}

\begin{remark}
    By (1) of Theorem \ref{202508250938}, $P_d(\mathtt{M})=\mathrm{V}(\mathtt{M};\mathtt{I}^{(d)}_{\mathcal{C}})$, so (2) of the theorem
    recovers the filtration (\ref{202609281811}) in Section \ref{202509301224}: the inclusions $P_d(\mathtt{M})\subset P_{d+1}(\mathtt{M})$ correspond to
    $\mathtt{I}^{(d+1)}_{\mathcal{C}}\subset\mathtt{I}^{(d)}_{\mathcal{C}}$.
\end{remark}

\begin{remark} \label{202509301226}
    When $\mathcal{C}$ is a category with finite coproducts and a zero object, and $\mathds{k} = \mathds{Z}$, the $\mathds{k}$-linear category corresponding to the quotient monad $\mathtt{L}_{\mathcal{C}} / \mathtt{I}^{(d)}_{\mathcal{C}}$ can essentially be found in \cite[Proposition 6.25]{hartl2015polynomial} (see Remark \ref{202409261755} for details).
\end{remark}

The following statement is well-known for additive categories $\mathcal{C}$.
Recall Definition \ref{202510021430}.

\begin{Corollary} \label{202508151849}
    The inclusion $\mathtt{L}_{\mathcal{C}}\mbox{-}\mathsf{Mod}^{\leq d} \hookrightarrow \mathtt{L}_{\mathcal{C}}\mbox{-}\mathsf{Mod}$ admits a left adjoint, which sends a left $\mathtt{L}_{\mathcal{C}}$-module $\mathtt{M}$ to 
    $Q_d(\mathtt{M}) {:=} \Lambda(\mathtt{M};\mathtt{I}^{(d)}_{\mathcal{C}})$,
    which is the maximal quotient of $\mathtt{M}$ of degree at most $d$.
\end{Corollary}
\begin{proof}
    By (1) of Theorem \ref{202508250938}, a quotient $\mathtt{M}\twoheadrightarrow\mathtt{Q}$ has $\deg\mathtt{Q}\le d$ if and only if $\mathtt{I}^{(d)}_{\mathcal{C}}\rhd\mathtt{Q}=0$.
    This holds if and only if the quotient map factors through
    $\mathtt{M}\twoheadrightarrow\mathtt{M}/(\mathtt{I}^{(d)}_{\mathcal{C}}\rhd\mathtt{M})=\Lambda(\mathtt{M};\mathtt{I}^{(d)}_{\mathcal{C}})$.
    Hence, $\Lambda(\mathtt{M};\mathtt{I}^{(d)}_{\mathcal{C}})$ is the maximal such quotient, and Proposition
    \ref{202509301238} applied to $\mathtt{J}=\mathtt{I}^{(d)}_{\mathcal{C}}$ gives the left adjoint in the assertion.
\end{proof}

We now prove Theorem \ref{202510011337}, which shows that the ideals $\mathtt{I}^{(d)}_{\mathcal{C}}$ yield an effective criterion for answering question~(Q1) raised in the introduction.

\begin{proof}[Proof of Theorem \ref{202510011337}]
    We first establish that for $d \in \mathds{N}$, the following two conditions are equivalent:
    \begin{enumerate}
        \item $\mathtt{I}^{(d+1)}_{\mathcal{C}} \subsetneq \mathtt{I}^{(d)}_{\mathcal{C}}$.
        \item $d +1 \in \mathrm{D}(\mathcal{C})$.
    \end{enumerate}
    By definition, $d+1 \in \mathrm{D}(\mathcal{C})$ if and only if $\mathcal{S}^{d}_{\mathcal{C}} \subsetneq \mathcal{S}^{d+1}_{\mathcal{C}}$.
    Corollary \ref{202509291354} applied to the filtration given in Theorem \ref{202508250938} shows that $\mathcal{S}^{d}_{\mathcal{C}} \subsetneq \mathcal{S}^{d+1}_{\mathcal{C}}$ if and only if $\mathtt{I}^{(d+1)}_{\mathcal{C}} \subsetneq \mathtt{I}^{(d)}_{\mathcal{C}}$, since $\mathcal{V} ( \mathtt{I}^{(d)}_{\mathcal{C}}) = \mathcal{S}^{d}_{\mathcal{C}}$ by Theorem \ref{202508250938}.
    
    The above assertion shows that $\mathtt{I}^{(d+1)}_{\mathcal{C}} \subsetneq \mathtt{I}^{(d)}_{\mathcal{C}}$ for every $d \in \mathds{N}$ if and only if $\mathrm{D}(\mathcal{C}) = \mathds{N}$; and $\mathtt{I}^{(0)}_{\mathcal{C}} = \mathtt{I}^{(d)}_{\mathcal{C}}$ for every $d \in \mathds{N}$ if and only if $\mathrm{D}(\mathcal{C}) = \{ 0\}$.
    By Theorem \ref{202609142120}, exactly one of these cases occurs.
    Hence, the former is equivalent to $1 \in \mathrm{D}(\mathcal{C})$ which is in turn equivalent to $\mathtt{I}^{(1)}_{\mathcal{C}} \subsetneq \mathtt{I}^{(0)}_{\mathcal{C}}$, and the latter is equivalent to $\mathtt{I}^{(1)}_{\mathcal{C}} = \mathtt{I}^{(0)}_{\mathcal{C}}$ by taking negations.
\end{proof}

We conclude this section giving an explicit description of the polynomiality ideal $\mathtt{I}^{(d)}_{\mathcal{C}}$.
Recall the notation considered in Remark \ref{202609282124}.

\begin{prop}\label{202609282123}
    Let $\mathcal{C}$ be a small monoidal category with zero unit.
    For $X_1,\cdots,X_{d+1}\in\mathcal{X}$ and $W=\bigvee_{i=1}^{d+1} X_i$, set
    $$
    \alpha_{X_1,\cdots,X_{d+1}}:=\sum_{\epsilon \in\{0,1\}^{d+1}}(-1)^{|\epsilon|}\,e_{X_1}^{\epsilon_1}\vee\cdots\vee e_{X_{d+1}}^{\epsilon_{d+1}}\in\mathtt{L}_{\mathcal{C}}(W,W)
    $$
    where $|\epsilon| = \sum_{i=1}^{d+1} \epsilon_i$.
    Then $\mathtt{I}^{(d)}_{\mathcal{C}}$ is the two-sided ideal of $\mathtt{L}_{\mathcal{C}}$ generated by
    $\alpha_{X_1,\cdots,X_{d+1}}$ for all $X_1,\cdots,X_{d+1}$.
\end{prop}
\begin{proof}
    Let $\mathtt{K}$ be the two-sided ideal generated by these elements.
    By the observation in Remark \ref{202609282124}, a left $\mathtt{L}_{\mathcal{C}}$-module $\mathtt{M}$ has $\deg\mathtt{M}\le d$ if and only if
    $\alpha_{X_1,\cdots,X_{d+1}}\rhd\mathtt{M}=0$ for all $X_1,\cdots,X_{d+1}$.
    Since $\mathtt{K}(Y,X)$ is spanned by elements $g\circ \alpha_{X_1,\cdots,X_{d+1}} \circ h$ for $g \in \mathtt{L}_{\mathcal{C}} (Y,W)$ and $h \in \mathtt{L}_{\mathcal{C}} (W,X)$, this is equivalent to $\mathtt{K}\rhd\mathtt{M}=0$.
    Hence $\mathcal{V}(\mathtt{K})=\mathcal{S}^d_{\mathcal{C}}$, i.e.\ $\mathtt{K}$ is a core internalizer of $\mathcal{S}^d_{\mathcal{C}}$,
    and $\mathtt{K}=\mathtt{I}^{(d)}_{\mathcal{C}}$ by Corollary \ref{202409271802} and Theorem \ref{202508250938}.
\end{proof}

\begin{remark}
    By Proposition \ref{202609282123}, for any monoidal category $\mathcal{C}$ with zero unit, $\mathtt{I}^{(d)}_{\mathcal{C}}$ is generated as a two-sided ideal by the elements $\alpha_{X_1,\cdots,X_{d+1}}$, one for each $(d+1)$-tuple of objects.
    In this paper, the main applications concern categories $\mathcal{C}$ with finite products and a zero object, regarded as monoidal categories with zero unit.
    When $\mathcal{C}$ has finite products, Section \ref{202401241445} shows that it is the left ideal generated by a single element $\pi^d_X$ for each object $X$.
    The finite product assumption on $\mathcal{C}$ is used to obtain the explicit generators $\pi^{d}_{X}$, on which the computations of Sections \ref{202510211820}--\ref{202509301437} rest.
\end{remark}

Recall the transposition from Definition \ref{202410131009}.

\begin{Corollary} \label{202609282136}
    Under $\mathtt{L}_{\mathcal{C}^{\mathsf{o}}}=\mathtt{L}_{\mathcal{C}}^{\mathrm{t}}$, we have
    $\mathtt{I}^{(d)}_{\mathcal{C}^{\mathsf{o}}}=\left(\mathtt{I}^{(d)}_{\mathcal{C}}\right)^{\mathrm{t}}$.
\end{Corollary}
\begin{proof}
    Since $\alpha_{X_1,\cdots,X_{d+1}}$ is expressed by $\vee$ and $e_X$'s, its counterpart in
    $\mathtt{L}_{\mathcal{C}^{\mathsf{o}}}(W,W)=\mathds{k}\mathcal{C}(W,W)$ is the same element.
    By Proposition \ref{202609282123}, $\mathtt{I}^{(d)}_{\mathcal{C}^{\mathsf{o}}}$ is the two-sided ideal of $\mathtt{L}_{\mathcal{C}^{\mathsf{o}}}$
    generated by these elements.
    Transposition exchanges left and right multiplication, so it carries the two-sided ideal generated by a family
    to the one generated by the same family, in $\mathtt{L}_{\mathcal{C}^{\mathsf{o}}}=\mathtt{L}_{\mathcal{C}}^{\mathrm{t}}$.
    By Proposition \ref{202609282123} again, this is $(\mathtt{I}^{(d)}_{\mathcal{C}})^{\mathrm{t}}$.
\end{proof}

\subsection{Polynomiality ideals and $\Gamma(\xi)$} \label{202512241413}

In this section, we give several alternative descriptions of the set $\Gamma(\xi)$ introduced in Section \ref{202510211607} (see Definition \ref{202609211701} below).
Let $\mathcal{C}$ and $\mathcal{D}$ be small monoidal categories with zero unit and the same object set $\mathcal{X}$.
Let $\xi : \mathcal{C} \to \mathcal{D}$ be a full strong monoidal functor that is the identity on objects.
The set $\Gamma(\xi)$ is recalled as follows:
\begin{Defn} \label{202609211701}
    $\Gamma(\xi) {:=} \{ d \in \mathds{N}\mid \xi^\ast : \mathtt{L}_{\mathcal{D}}\mbox{-}\mathsf{Mod}^{\leq d} \stackrel{\simeq}{\to} \mathtt{L}_{\mathcal{C}}\mbox{-}\mathsf{Mod}^{\leq d} \}$.
\end{Defn}

\begin{Defn}
    Using the ideals given in Theorem \ref{202508250938}, we define
    $$
    \mathrm{B}(\xi) {:=} \{  d \in \mathds{N} \mid \xi:  \mathtt{I}^{(d-1)}_{\mathcal{C}}/\mathtt{I}^{(d)}_{\mathcal{C}} \stackrel{\cong}{\longrightarrow}  \mathtt{I}^{(d-1)}_{\mathcal{D}}/\mathtt{I}^{(d)}_{\mathcal{D}} \}, \quad \mathrm{where~} \mathtt{I}^{(-1)}_{\mathcal{C}} {:=} \mathtt{L}_{\mathcal{C}} .
    $$
\end{Defn}

\begin{notation} \label{202510221957}
    For a subset $S \subset \mathds{N}$, we denote by $S^\spadesuit$ the maximal subset of consecutive integers starting at $0$ contained in $S$.
    For instance, $\{0,1,2,5,8\}^{\spadesuit} = \{0,1,2\}$.
\end{notation}

Recall the ideal $\mathtt{K}_{\xi}$ associated with the functor $\xi$, introduced in Example \ref{202510261038}.
\begin{theorem} \label{202509291857}
    For $d\in \mathds{N}$, the following conditions are equivalent:
    \begin{enumerate}
        \item $d \in \Gamma(\xi)$.
        \item The map $\xi^\ast :\mathcal{S}^{d}_{\mathcal{D}} \to \mathcal{S}^{d}_{\mathcal{C}}$ induced by $\xi$ is bijective.
        \item The map $\xi: \mathtt{L}_{\mathcal{C}}/\mathtt{I}^{(d)}_{\mathcal{C}} \to \mathtt{L}_{\mathcal{D}}/\mathtt{I}^{(d)}_{\mathcal{D}}$ is an isomorphism of $(\mathcal{X} \times \mathcal{X})$-indexed $\mathds{k}$-modules.
        \item $\mathtt{K}_{\xi} \subset \mathtt{I}^{(d)}_{\mathcal{C}}$.
    \end{enumerate}
    Furthermore, we have $\Gamma (\xi) = \mathrm{B}(\xi)^\spadesuit$.
    In particular, $\Gamma(\xi)$ consists of consecutive nonnegative integers starting at $0$.
\end{theorem}
\begin{proof}
    We begin by proving (1) $\Leftrightarrow$ (2).
    Since $\xi: \mathcal{C} \to \mathcal{D}$ is full, it is immediate that $\xi^\ast : \mathtt{L}_{\mathcal{D}}\mbox{-}\mathsf{Mod}^{\leq d} \to \mathtt{L}_{\mathcal{C}}\mbox{-}\mathsf{Mod}^{\leq d}$ is fully faithful.
    Moreover, the induced monad map $\xi : \mathtt{L}_{\mathcal{C}} \to \mathtt{L}_{\mathcal{D}}$ is an epimorphism, so the map $\xi^\ast : \mathcal{S}^{d}_{\mathcal{D}} \to  \mathcal{S}^{d}_{\mathcal{C}}$ is injective by (1) of Lemma \ref{202509291501}.
    Hence, $\xi^\ast : \mathtt{L}_{\mathcal{D}}\mbox{-}\mathsf{Mod}^{\leq d} \to \mathtt{L}_{\mathcal{C}}\mbox{-}\mathsf{Mod}^{\leq d}$ is essentially surjective if and only if $\xi^\ast :\mathcal{S}^{d}_{\mathcal{D}} \to \mathcal{S}^{d}_{\mathcal{C}}$ is surjective.
    These results prove (1) $\Leftrightarrow$ (2).

    Next, we prove (2) $\Leftrightarrow$ (3).
    By (3) of Theorem \ref{202508250938}, the epimorphism $\xi: \mathtt{L}_{\mathcal{C}} \to \mathtt{L}_{\mathcal{D}}$ satisfies $\xi (\mathtt{I}^{(d)}_{\mathcal{C}}) = \mathtt{I}^{(d)}_{\mathcal{D}}$.
    Thus, we obtain (2) $\Leftrightarrow$ (3) from (3) of Lemma \ref{202509291501}, since, by Theorem \ref{202508250938}, $\mathtt{I}^{(d)}_{\mathcal{C}}$ is a core internalizer for $\mathcal{S}^{d}_{\mathcal{C}}$.

    We prove (1) $\Leftrightarrow$ (4).
    From Theorem \ref{202508250938} and the discussion in Example \ref{202510261038}, it follows that an $\mathtt{L}_{\mathcal{C}}$-module $\mathtt{M}$ induces an $\mathtt{L}_{\mathcal{D}}$-module of $\deg \leq d$ if and only if $\mathtt{K}_{\xi} \rhd \mathtt{M} = 0 = \mathtt{I}^{(d)}_{\mathcal{C}} \rhd \mathtt{M}$.
    By definition, this condition is equivalent to
    $[\mathtt{M}] \in \mathcal{V} ( \mathtt{K}_{\xi} + \mathtt{I}^{(d)}_{\mathcal{C}})$.
    Therefore, for $d \in \mathds{N}$, $d \in \Gamma(\xi)$ if and only if $$\mathcal{V} ( \mathtt{K}_{\xi} + \mathtt{I}^{(d)}_{\mathcal{C}}) = \mathcal{V} (\mathtt{I}^{(d)}_{\mathcal{C}}).$$
    By Corollary \ref{202509291354}, this is equivalent to $\mathtt{K}_{\xi} + \mathtt{I}^{(d)}_{\mathcal{C}} = \mathtt{I}^{(d)}_{\mathcal{C}}$, i.e. $\mathtt{K}_{\xi} \subset \mathtt{I}^{(d)}_{\mathcal{C}}$.

    The set $\Gamma (\xi)$ is a set of consecutive integers starting at $0$.
    In fact, the above results imply that it equals the set of $d \in \mathds{N}$ such that $\mathtt{K}_{\xi} \subset \mathtt{I}^{(d)}_{\mathcal{C}}$.
    Therefore, if $d \in \Gamma (\xi)$, then $\mathtt{K}_{\xi} \subset \mathtt{I}^{(d)}_{\mathcal{C}}$ so $\mathtt{K}_{\xi} \subset \mathtt{I}^{(d')}_{\mathcal{C}}$ for $0 \leq d' \leq d$ by (2) of Theorem \ref{202508250938}.
    
    Based on the above results, we prove $\Gamma (\xi) \subset \mathrm{B}(\xi)$.
    Let $d \in \Gamma(\xi)$.
    Then we have $[d] \subset \Gamma(\xi)$, since $\Gamma (\xi)$ is an interval.
    Using (1) $\Leftrightarrow$ (3) and the following usual commutative diagram (in the category of $(\mathcal{X}\times\mathcal{X})$-indexed $\mathds{k}$-modules), one can recursively prove that $\mathtt{I}^{(d-1)}_{\mathcal{C}}/\mathtt{I}^{(d)}_{\mathcal{C}} \stackrel{\cong}{\longrightarrow}  \mathtt{I}^{(d-1)}_{\mathcal{D}}/\mathtt{I}^{(d)}_{\mathcal{D}}$ is an isomorphism.
    Hence, $d \in \mathrm{B}(\xi)$.
    $$
    \begin{tikzcd}
        0 \ar[r]\ar[d, equal]  &
        \mathtt{I}^{(d-1)}_{\mathcal{C}} / \mathtt{I}^{(d)}_{\mathcal{C}} \ar[r] \ar[d, "\xi"] & \mathtt{L}_{\mathcal{C}} / \mathtt{I}^{(d)}_{\mathcal{C}}  \ar[r] \ar[d, "\xi"] & \mathtt{L}_{\mathcal{C}} / \mathtt{I}^{(d-1)}_{\mathcal{C}}  \ar[r] \ar[d, "\xi"] & 0\ar[d, equal]  \\
        0 \ar[r] &
        \mathtt{I}^{(d-1)}_{\mathcal{D}} / \mathtt{I}^{(d)}_{\mathcal{D}}  \ar[r] & \mathtt{L}_{\mathcal{D}} / \mathtt{I}^{(d)}_{\mathcal{D}}  \ar[r] & \mathtt{L}_{\mathcal{D}} / \mathtt{I}^{(d-1)}_{\mathcal{D}}  \ar[r] & 0
    \end{tikzcd}
    $$
    
    If $\Gamma(\xi)$ is infinite, then $\Gamma(\xi) = \mathds{N}$, since $\Gamma(\xi)$ is an interval of $\mathds{N}$ containing $0$.
    By the above result, $\mathrm{B}(\xi) = \mathds{N}$.
    
    We now assume that $\Gamma(\xi)$ is finite, and prove the last assertion.
    Let $d_0$ be the maximum of $\Gamma(\xi)$, and it suffices to show that $d_0 +1 \not\in \mathrm{B}(\xi)$.
    If $d_0 + 1 \in \mathrm{B}(\xi)$, then, considering the above commutative diagram for $d = d_0 +1$, the five-lemma implies that the middle map $\xi: \mathtt{L}_{\mathcal{C}} / \mathtt{I}^{(d_0+1)}_{\mathcal{C}} \cong \mathtt{L}_{\mathcal{D}} / \mathtt{I}^{(d_0+1)}_{\mathcal{D}}$ is an isomorphism.
    In fact, by the hypothesis $d_0 + 1 \in \mathrm{B}(\xi)$, the left map is an isomorphism, and, by $d_0 \in \Gamma(\xi)$, the right map is an isomorphism.
    Hence, we obtain $d_0 +1 \in \Gamma(\xi)$; this contradicts the maximality of $d_0 \in \Gamma(\xi)$.
\end{proof}

\begin{Corollary} \label{202510141158}
    $\Gamma (\xi: \mathcal{C}\to \mathcal{D}) = \Gamma(\xi^{\mathsf{o}}: \mathcal{C}^{\mathsf{o}}\to \mathcal{D}^{\mathsf{o}})$.
\end{Corollary}
\begin{proof}
    By Theorem \ref{202509291857}, it is sufficient to prove $\mathrm{B}(\xi)^{\spadesuit} = \mathrm{B}(\xi^{\mathsf{o}})^{\spadesuit}$.
    Note that these sets are defined by using the polynomiality ideals associated with $\mathcal{C}$, $\mathcal{D}$ and their opposites.
    One can directly show the assertion by using Corollary \ref{202609282136}.
\end{proof}

\subsection{Polynomiality ideals for a category with finite products}
\label{202401241445}

In this section, we present a useful form of the polynomiality ideal $\mathtt{I}^{(d)}_{\mathcal{C}} \subset \mathtt{L}_{\mathcal{C}}$ (see Theorem \ref{202609281707}) when $\mathcal{C}$ is a category having finite products and a zero object, regarded as a monoidal category with zero unit.

\begin{notation} \label{202510281503}
    For $a_i,b_j \in \mathds{k}$ and morphisms $f_i \in \mathcal{C}(X,Y) ,g_j \in \mathcal{C}(X^\prime,Y^\prime)$, we write
    $$
    \sum_{i} a_i f_i \times \sum_{j} b_j g_j {:=} \sum_{i,j} a_ib_j (f_i \times g_j) \in \mathtt{L}_{\mathcal{C}} ( Y \times Y^\prime , X \times X^\prime).
    $$
    More generally, for $d \in \mathds{N}$ and $f \in \mathtt{L}_{\mathcal{C}}(Y,X)$, we denote by $f^{\times d}$ the product $f \times f \times \cdots \times f$ where $f$ appears $d$-times.
\end{notation}

\begin{Defn} \label{202403041329}
    Let $d \in \mathds{N}$ and $X \in \mathcal{X}$.
    For $k \in [d] = \{ 0,1, \cdots, d\}$, let $p_{k+1} \in \mathcal{C} ( X^{\times (d+1)} , X)$ be the $(k+1)$-th projection.
    For $T \subset [d]$, we take $D^{T}_{X} \in \mathcal{C} (X, X^{\times (d+1)})$ that is characterized by
    \begin{align} \label{202407311833}
        p_{k+1} \circ D^{T}_{X} =
        \begin{cases}
            \mathrm{id}_{X} & (k \in T) , \\
            e_{X} & (\mathrm{otherwise}) .
        \end{cases}
    \end{align}
    We then define $\pi^{\mathcal{C},d}_{X}$ as
    \begin{align*}
        \pi^{\mathcal{C},d}_{X} {:=} \sum_{T \subset [d]} (-1)^{d+1-|T|} \cdot D^{T}_X \in \mathtt{L}_{\mathcal{C}} ( X^{\times (d+1)} , X ) .
    \end{align*}
    Using Notation \ref{202510281503}, we give an equivalent definition in closed form:
    $$
    \pi^{\mathcal{C},d}_{X} {:=} (\mathrm{id}_X - e_{X})^{\times (d+1)} \circ \Delta^{(d+1)}_X 
    $$
    where $\Delta^{(d+1)}_X \in \mathcal{C}( X , X^{\times (d+1)})$ is the $(d+1)$-fold diagonal map on $X$.
    If $\mathcal{C}$ is understood, then we write $\pi^d_{X} {:=} \pi^{\mathcal{C},d}_{X}$.
    For instance, $\pi^{\mathcal{C},0}_{X} = \mathrm{id}_{X} - e_X$.
\end{Defn}

\begin{theorem} \label{202609281707}
    The polynomiality ideal $\mathtt{I}^{(d)}_{\mathcal{C}}$ equals the left ideal of $\mathtt{L}_{\mathcal{C}}$ generated by $\pi^{d}_{X}$ for $X \in \mathcal{X}$.
\end{theorem}

\begin{remark} 
    (1) of Theorem \ref{202508250938} gives $P_d \cong \mathrm{V}(-;\mathtt{I}^{(d)}_{\mathcal{C}})$.
    Hence, it follows from Theorem \ref{202609281707} that $P_d(\mathtt{M})(X)$ is the kernel of $\pi^d_{X} \rhd (-) : \mathtt{M}(X) \to \mathtt{M}(X^{\times (d+1)})$.
    This description was already obtained from \cite[Section 3.3]{kim2024analytic}.
\end{remark}

In the rest of this section, we prove Theorem \ref{202609281707}.
Fix $d \in \mathds{N}$.
We set $\mathtt{J}^d$ as the left ideal of $\mathtt{L}_{\mathcal{C}}$ generated by $\pi^{d}_{X}$ for $X \in \mathcal{X}$, and we shall prove $\mathtt{J}^d = \mathtt{I}^{(d)}_{\mathcal{C}}$.

\begin{Lemma} \label{202401231807}
    For $f \in \mathcal{C}(X,Y) \subset \mathtt{L}_{\mathcal{C}} (Y,X)$, we have $\pi^{d}_Y \circ f = (f^{\times (d+1)} ) \circ \pi^{d}_X$.
\end{Lemma}
\begin{proof}
    The universality of the product gives $\Delta^{(d+1)}_Y \circ f = (f^{\times (d+1)}) \circ \Delta^{(d+1)}_X$, and that of the zero object $\ast$ gives $e_Y \circ f = f \circ e_{X}$.
    Hence, the statement follows from the definition of $\pi^{d}_{X}$.
\end{proof}

\begin{Lemma} \label{202401241347}
    The left ideal $\mathtt{J}^d$ is a two-sided ideal of the monad $\mathtt{L}_{\mathcal{C}}$.
\end{Lemma}
\begin{proof}
    Let $X,Y \in \mathcal{X}$.
    Any element in $\mathtt{L}_{\mathcal{C}} (Y,X)$ is given by a linear combination of some elements of $\mathcal{C} ( X,Y)$, say $f$'s.
    By Lemma \ref{202401231807}, we obtain $\pi^{d}_Y \circ f = (f^{\times (d+1)} ) \circ \pi^{d}_X \in \mathtt{J}^d ( Y^{\times (d+1)} ,X )$.
    Therefore, $\mathtt{J}^d$ is a right ideal of $\mathtt{L}_{\mathcal{C}}$, as required.
\end{proof}

We now prove Theorem \ref{202609281707}.
\begin{proof}[Proof of Theorem \ref{202609281707}]
    We prove $\mathtt{J}^d=\mathtt{I}^{(d)}_{\mathcal{C}}$.
    Since $\alpha_{X,\cdots,X}=(\mathrm{id}_X-e_X)^{\times(d+1)}$, we have
    $$\pi^d_X=\alpha_{X,\cdots,X}\circ\Delta^{(d+1)}_X,$$
    which lies in $\mathtt{I}^{(d)}_{\mathcal{C}}$ by Proposition \ref{202609282123}.
    Hence $\mathtt{J}^d\subset\mathtt{I}^{(d)}_{\mathcal{C}}$.
    Conversely, let $X_1,\cdots,X_{d+1}\in\mathcal{X}$, $W=\prod_{i=1}^{d+1} X_i$, and let $\mathrm{pr}_i \in \mathcal{C} (W, X_i)$ be the $i$-th projection.
    Put $q {:=} \prod_{i=1}^{d+1} \mathrm{pr}_i \in \mathcal{C}( W^{\times(d+1)}, W)$.
    For $T\subset[d]$, the $i$-th component of $q\circ D^T_W$ is $\mathrm{pr}_i$ if $i-1\in T$, and factors through $\ast$ otherwise, so $q\circ D^T_W=\prod_{i=1}^{d+1} e_{X_i}^{\epsilon_i}$ where $\epsilon \in \{0,1\}^{d+1}$ such that $\epsilon_i=1$ if and only if $i-1\notin T$.
    Since $\sum_{i=1}^{d+1} \epsilon_i =d+1-|T|$, we obtain $$q\circ\pi^d_W=\alpha_{X_1,\cdots,X_{d+1}}.$$
    This proves $\alpha_{X_1,\cdots,X_{d+1}}\in\mathtt{J}^d (W,W)$.
    Since $\mathtt{J}^d$ is a two-sided ideal by Lemma \ref{202401241347}, Proposition \ref{202609282123} gives $\mathtt{I}^{(d)}_{\mathcal{C}}\subset\mathtt{J}^d$.
\end{proof}

\subsection{Polynomiality ideals for a category with finite coproducts}
\label{202508152158}

Let $\mathcal{D}$ be a small category having finite coproducts and a zero object $\ast$, the setting dual to that of Section \ref{202401241445}.
The purpose of this section is to clarify Remark \ref{202509301226}.
For an object $X$, let $\nabla^{(d+1)}_{X} \in \mathcal{D} (X^{\amalg (d+1)}, X)$ denote the $(d+1)$-fold folding map on $X$.
The polynomiality ideal $\mathtt{I}^{(d)}_{\mathcal{D}}$ admits a generator description, as in Theorem \ref{202609281707}:
\begin{Corollary} \label{202609291033}
    The polynomiality ideal $\mathtt{I}^{(d)}_{\mathcal{D}}$ equals the right ideal of $\mathtt{L}_{\mathcal{D}}$ generated by $\nabla^{(d+1)}_{X} \circ (\mathrm{id}_{X}- e_{X})^{\amalg (d+1)}$ for $X \in \mathcal{X}$.
\end{Corollary}
\begin{proof}
    Note that $\mathcal{D}^{\mathsf{o}}$ admits finite products and a zero object $\ast$.
    Hence, the results in Section \ref{202401241445} apply to $\mathcal{D}^{\mathsf{o}}$.
    The desired statement follows from Theorem \ref{202609281707} applied to $\mathcal{D}^{\mathsf{o}}$, since $\mathtt{I}^{(d)}_{\mathcal{D}} = \left( \mathtt{I}^{(d)}_{\mathcal{D}^{\mathsf{o}}} \right)^{\mathrm{t}}$ by Corollary \ref{202609282136}.
\end{proof}

We now give several remarks concerning the results of Hartl-Pirashvili-Vespa \cite{hartl2015polynomial}.

\begin{remark}
    Let $\mathcal{D}$ be a small category having finite coproducts and a zero object $\ast$.
    Recall from Corollary \ref{202508151849} that $Q_d(\mathtt{M})=\Lambda(\mathtt{M};\mathtt{I}^{(d)}_{\mathcal{D}})$ is the maximal
    quotient of $\mathtt{M}$ of degree $\leq d$.
    We remark that $Q_d(\mathtt{M})$ coincides with the $d$-Taylorization of $\mathtt{M}$
    introduced in \cite[Definition 3.16]{hartl2015polynomial}: by
    Corollary \ref{202609291033}, $\mathtt{I}^{(d)}_{\mathcal{D}}(X,-)$ is spanned by right multiples of $f_X=\nabla^{(d+1)}_X\circ(\mathrm{id}_X-e_X)^{\amalg(d+1)}$, so
    $$
    \Lambda(\mathtt{M};\mathtt{I}^{(d)}_{\mathcal{D}})(X)=\mathrm{coker}\big(f_X\rhd(-):\mathtt{M}(X^{\amalg(d+1)})\to\mathtt{M}(X)\big).
    $$
\end{remark}

\begin{remark} \label{202409261755}
    We consider the $\mathds{Z}$-linear category $T_d\bar{\mathds{Z}}[\mathcal{D}]$ introduced in \cite[Proposition 6.25]{hartl2015polynomial}; we refer the reader there for the definition.
    By abuse of notation, we also write $T_d\bar{\mathds{Z}}[\mathcal{D}]$ for the corresponding monad in $\mathsf{Mat}_{\mathds{Z}}$ on $\mathcal{X} = \mathrm{Obj} (\mathcal{D})$ via the equivalence (\ref{202401231608}).
    In \cite{hartl2015polynomial}, the main objects of study are {\it reduced} functors from $\mathcal{D}$.
    To compare their results, we recall the ideal $\mathtt{Z}_{\mathcal{D}}$ introduced in Example \ref{202509021729}.
    By combining the results in Example \ref{202509021729} and Corollary \ref{202508151849}, we obtain an equivalence of categories:
    \begin{align*}
        q^\ast : \mathtt{L}_{\mathcal{D}} / (\mathtt{I}^{(d)}_{\mathcal{D}} + \mathtt{Z}_{\mathcal{D}} )\mbox{-} \mathsf{Mod} \stackrel{\simeq}{\to} \{ \mathtt{M} \in \mathtt{L}_{\mathcal{D}}\mbox{-} \mathsf{Mod}^{\leq d}\mid \mathtt{Z}_{\mathcal{D}} \rhd \mathtt{M} \cong 0 \} 
    \end{align*}
    where the right hand side denotes the full subcategory corresponding to the condition.
    By applying Corollary \ref{202508151849} to $\mathtt{M}= \mathtt{L}_{\mathcal{D}} / \mathtt{Z}_{\mathcal{D}}$, we obtain an isomorphism
    \begin{align*}
        \mathtt{L}_{\mathcal{D}} / (\mathtt{I}^{(d)}_{\mathcal{D}} + \mathtt{Z}_{\mathcal{D}} ) \cong Q_d ( \mathtt{L}_{\mathcal{D}} / \mathtt{Z}_{\mathcal{D}}) .
    \end{align*}
    In particular, the $(\mathcal{X}\times\mathcal{X})$-indexed module in the right hand side inherits the monad structure from the quotient monad $ \mathtt{L}_{\mathcal{D}} / (\mathtt{I}^{(d)}_{\mathcal{D}} + \mathtt{Z}_{\mathcal{D}} )$.
    By definitions, for $\mathds{k}= \mathds{Z}$, we have $$Q_d ( \mathtt{L}_{\mathcal{D}} / \mathtt{Z}_{\mathcal{D}}) = T_d\bar{\mathds{Z}}[\mathcal{D}] .$$
    Under this identification, the above equivalence of categories leads to a proof of \cite[Corollary 6.29]{hartl2015polynomial} for a general ground ring $\mathds{k}$.
\end{remark}

\section{Application to Lawvere theories}
\label{202510211820}

An algebraic theory in the sense of Lawvere \cite{Lawvere1963}, a {\it Lawvere theory} for short, is a small category with finite products, whose set of objects is $\mathds{N}$, with the product on objects given by the addition of natural numbers.
In particular, $0$ is a terminal object.
In this section, we study the polynomiality ideal $\mathtt{I}^{(d)}_{\mathcal{C}}$ associated with a Lawvere theory $\mathcal{C}$.
In particular, we give, in Sections \ref{202510161122} and \ref{202404121704}, the explicit statement of Theorem \ref{202509241101} with the proof.

\subsection{Some conditions on Lawvere theory} \label{202509031755}

In this section, we introduce some conditions of a Lawvere theory that serve as the main source of examples in this paper.
Let $\mathcal{C}$ be a Lawvere theory.

\begin{notation} \label{202512241120}
    For $n,m \in \mathds{N}$, let $\mathcal{C}(n,m)$ denote the set of morphisms from $n$ to $m$.
    We set $$\mathcal{C}_n {:=} \mathcal{C} (n,1).$$
\end{notation}
The universality of products gives a bijection:
\begin{align*}
    \mathcal{C}_n^{\times m} \stackrel{\cong}{\to} \mathcal{C}(n,m) ; \quad (f_1, \cdots, f_m) \mapsto (f_1 \times \cdots \times f_m) \circ  \Delta^{(m)}_{n}
\end{align*}
Based on this identification, throughout this paper, for $f_1, \cdots, f_m \in \mathcal{C}_n$, we denote by $(f_1, \cdots, f_m)$ the corresponding morphism from $n$ to $m$.

In the following, we introduce some conditions on $\mathcal{C}$, called (Z), (M) and (M*).
\begin{itemize}
    \item[{\bf (Z)}] $\mathcal{C}$ has a zero object, i.e. an object that is initial and terminal.
\end{itemize}
This condition, by the uniqueness of a terminal object, implies that the object $0$ is a zero object.

\begin{itemize}
    \item[{\bf (M)}] The object $1$ of $\mathcal{C}$ is equipped with a monoid object structure.
    We denote by $\nabla : 2= (1+1) \to 1$ the multiplication and $\eta: 0 \to 1$ the unit.
\end{itemize}
By this condition, the set $\mathcal{C}_n$ obtains a monoid structure: 
for $f,g \in \mathcal{C}_n$, we define
$$
f \star g {:=} \nabla \circ ( f \times g ) \circ \Delta_{n} 
$$
where $\Delta_{n} : n \to (n+n) = 2n$ is the diagonal map.
Then the morphism $e : n \to 1$ defined by the composition $n \stackrel{z}{\to} 0 \stackrel{\eta}{\to} 1$ gives the unit of the monoid $\mathcal{C}_n$ where $z$ is the morphism to the terminal object.

\begin{remark}
    The assignment of the set $\mathcal{C}_n = \mathcal{C}(n,1)$ to $n \in \mathds{N}$ induces a functor $\mathcal{C}^{\mathsf{o}} \to \mathsf{Set}$ where $\mathsf{Set}$ denotes the category of sets.
    One can check that the condition (M) holds if and only if the functor $\mathcal{C}_{(-)}$ is endowed with a factorization:
    $$
    \begin{tikzcd}
        & \mathsf{Mon} \ar[d,"U"] \\
        \mathcal{C}^{\mathsf{o}} \ar[r, "\mathcal{C}_{(-)}"'] \ar[ur, "\exists"] & \mathsf{Set}
    \end{tikzcd}
    $$
    where $\mathsf{Mon}$ is the category of monoids, and $U$ denotes the forgetful functor.
\end{remark}

We introduce a condition stronger than (M) as follows:
\begin{itemize}
    \item[{\bf (M*)}] 
        Let $p_{k} : n \to 1$ be the $k$-th projection from $n = \overbrace{1 + \cdots + 1}^{n}$ to $1$.
        The Lawvere theory $\mathcal{C}$ satisfies (M), and, for $n\in \mathds{N}$, the union $\bigcup^{n}_{k=1} p_{k}^\ast ( \mathcal{C}_1 )$ generates the monoid $\mathcal{C}_n$.
\end{itemize}

We simply write $\mathrm{(ZM^*)} = \mathrm{(Z)} + \mathrm{(M^*)}$.
In what follows, we give some examples of $\mathcal{C}$ satisfying condition (ZM*).
We present further examples in Section \ref{202510071642} after introducing some general notions.

\begin{Example} \label{202512281213}
    Let $R$ be a unital (not necessarily commutative) ring.
    Let $\mathbf{M}_{R}$ be the category of free left $R$-modules of finite rank.
    It is a Lawvere theory satisfying condition (ZM*).
\end{Example}

\begin{Example} \label{202509301105}
    Let $\mathbf{W}$ be the category of finitely generated free monoids and monoid maps.
    The opposite category $\mathbf{W}^{\mathsf{o}}$ is a Lawvere theory satisfying (ZM*), and the associated monoid $\mathcal{C}_n$ is the free monoid generated by $n$ elements.
    We recall that an idempotent monoid is a monoid in which any element is an idempotent.
    There are several variants of $\mathbf{W}$ whose opposites satisfy (ZM*):
    \begin{itemize}
        \item $\mathbf{W}_{\mathrm{id}}$: the category of finitely generated free idempotent monoids.
        \item $\mathbf{W}_{\mathrm{c}}$: the category of finitely generated free commutative monoids.
        \item $\mathbf{W}_{\mathrm{c,id}}$: the category of finitely generated free commutative idempotent monoids.
    \end{itemize}
\end{Example}

\subsection{Polynomiality ideals and augmentation ideals}
\label{202510161122}

In this section, we investigate the polynomiality ideals associated with a Lawvere theory $\mathcal{C}$.
We show that, if $\mathcal{C}$ satisfies condition (ZM*), then these ideals can be identified with the augmentation ideals of certain monoids naturally associated with $\mathcal{C}$.

\begin{Defn}
    For a monoid $M$, let $\mathds{k} [M]$ denote the {\it monoid algebra} associated with $M$, whose multiplication is, in particular, induced by the monoid structure of $M$.
    We define the {\it augmentation ideal} $\Aug (M)$ to be the two-sided ideal of $\mathds{k} [M]$ generated by $(g-h)$ for $g,h \in M$.
    Denote by $\Aug (M)^d$ the $d$-th power of the ideal $\Aug (M)$.
\end{Defn}

\begin{notation} \label{202609141203}
    In this paper, $\Aug (M)$ denotes the augmentation ideal over a general ground ring $\mathds{k}$.
    When we consider a specific ground ring such as $\mathds{Z}$ or $\mathds{Q}$, we use $\Aug_{\mathds{Z}}$ or $\Aug_{\mathds{Q}}$.
\end{notation}

We fix $n,m,d \in \mathds{N}$.
Using the notation and arguments in Section \ref{202509031755}, we have an isomorphism $\mathds{k}[\mathcal{C}_n^{\times m}] \cong \mathds{k}[\mathcal{C}(n,m)] =  \mathtt{L}_{\mathcal{C}} (m,n)$.

\begin{theorem} \label{202408011134}
    Under the isomorphism $\mathds{k}[\mathcal{C}_n^{\times m}] \cong \mathtt{L}_{\mathcal{C}} (m,n)$, we have
    $$
    \Aug \left( \mathcal{C}_n^{\times m} \right)^{d+1} \cong \mathtt{I}^{(d)}_{\mathcal{C}} (m,n)  .
    $$
\end{theorem}
The proof is deferred to Section \ref{202404121704}.
Note that the augmentation ideal is defined using the monoid structure induced by condition (M), whereas the polynomiality ideal depends only on the category $\mathcal{C}$.

\begin{remark}
    To the best of our knowledge, the ideals $\mathtt{I}^{(d)}_{\mathcal{C}}$ of the category $\mathcal{C}$ have not been considered in the literature.
    For some specific categories, however, the left-hand side of Theorem~\ref{202408011134}, namely the powers of the augmentation ideal, has already appeared in the study of polynomial functors: see \cite[Proposition 2.7]{djament2016cohomologie} for $\mathcal{C}= \G^{\mathsf{o}}$, and the construction of the ringoid $P^{d} \mathcal{C}$ in \cite[Section 3]{Piras1993} for additive categories $\mathcal{C}$. 
\end{remark}

In what follows, we give two important consequences of the theorem.
The first is related to the set $\mathrm{D}(\mathcal{C})$ introduced in Definition \ref{202511251137}; we will investigate this consequence in Section \ref{202510171330} where some of our main theorems are proved.
The second is the following corollary which states that the polynomial property of $\mathtt{L}_{\mathcal{C}}$-modules can be understood from the classical notion due to Passi.
We refer to \cite[Chapter V]{MR537126} for the degree of a map from a monoid.

\begin{Corollary} \label{202509251010}
    A left $\mathtt{L}_{\mathcal{C}}$-module $\mathtt{M}$ has polynomial degree $\leq d$ if and only if, for $n,m \in \mathds{N}$, the following map has polynomial degree $\leq d$ (in the sense of Passi) with respect to the monoid $\mathcal{C}_n^{\times m}$:
    \begin{align*}
        \mathcal{C}_n^{\times m} \hookrightarrow  \mathds{k}\mathcal{C}(n,m) = \mathtt{L}_{\mathcal{C}} (m,n)\stackrel{\mathtt{M}}{\to} \hom_{\mathds{k}} ( \mathtt{M} (n) , \mathtt{M} (m) ) .
    \end{align*}
\end{Corollary}
\begin{proof}
    By Theorems \ref{202508250938} and \ref{202408011134}, $\deg \mathtt{M} \leq d$ if and only if $I(\mathcal{C}_n^{\times m})^{d+1} \cong \mathtt{I}^{(d)}_{\mathcal{C}}(m,n)$ acts trivially on $\mathtt{M}(n)$ for all $n,m$, which is the assertion.
\end{proof}

\begin{remark}
    This result is well known for an additive category.
    For instance, see \cite[Example 1.3, Chapter V]{MR537126}.
\end{remark}

\begin{remark}
    The corollary resembles the application of \cite[Theorem 6.27]{hartl2015polynomial} to $\mathcal{C}^{\mathsf{o}}$; however, Corollary \ref{202509251010} is a stronger statement which can be described by using the monoid structures, especially the assumption (M).
    To clarify this, we restate the theorem in our context: it asserts that a left $\mathtt{L}_{\mathcal{C}^{\mathsf{o}}}$-module $\mathtt{M}$ has polynomial degree $\leq d$ if and only if $\mathcal{C} (n,m) \to \hom_{\mathds{k}} ( \mathtt{M} (m) , \mathtt{M} (n) )$ is polynomial of degree $\leq d$ in the sense of \cite[Definition 6.23]{hartl2015polynomial}.
    The latter is equivalent to the condition that $\mathtt{L}_{\mathcal{C}}(m,n)  \to \hom_{\mathds{k}} ( \mathtt{M} (m) , \mathtt{M} (n) )$ factors through $\mathtt{L}_{\mathcal{C}}(m,n) \to (\mathtt{L}_{\mathcal{C}}/\mathtt{I}^{(d)}_{\mathcal{C}})(m,n)$, which does not imply the previous corollary.
\end{remark}

\subsection{A computable criterion for $\mathrm{D}(\mathcal{C})$}
\label{202510171330}

In this section, we apply Theorem \ref{202408011134} to study the set $\mathrm{D}(\mathcal{C})$ introduced in Definition \ref{202511251137} when $\mathcal{C}$ is a Lawvere theory satisfying condition (ZM*).
The following theorem is the main result of this section.

\begin{notation}
    For a monoid $M$, let $\mathfrak{Q}_1 (M) {:=} \Aug_{\mathds{Z}} (M)/\Aug_{\mathds{Z}} (M)^2$.
\end{notation}

\begin{theorem} \label{202512241009}
    Let $\mathcal{C}$ be a Lawvere theory satisfying condition (ZM*).
    The following conditions are equivalent:
    \begin{enumerate}
        \item $\mathrm{D}(\mathcal{C}) = \{ 0 \}$; in other words, every polynomial $\mathcal{C}$-module is constant.
        \item $\mathds{k} \otimes_{\mathds{Z}} \mathfrak{Q}_1 (\mathcal{C}_1) \cong 0$.
    \end{enumerate}
    Otherwise, $\mathrm{D}(\mathcal{C}) = \mathds{N}$.
\end{theorem}

\begin{remark}
    This theorem gives a {\it non-additive} generalization of 
    \cite[Proposition 2.13]{DTV2023}.
\end{remark}

\begin{Example} \label{202509241443}
    Examples satisfying $\mathrm{D}(\mathcal{C}) = \{ 0 \}$ include $\mathbf{M}_{R}$ for a ring $R$ such that $\mathds{k} \otimes_{\mathds{Z}} R \cong 0$.
    Moreover, the category $\mathcal{C} \in \{ \mathbf{W}_{\mathrm{id}}^{\mathsf{o}} , \mathbf{W}_{\mathrm{c,id}}^{\mathsf{o}}\}$ given in Example \ref{202509301105} satisfies $\mathrm{D}(\mathcal{C}) = \{ 0\}$.
    In fact, $\mathcal{C}_1$ is the idempotent monoid with one generator, so $\mathfrak{Q}_1 (\mathcal{C}_1)$ is trivial.
\end{Example}

\begin{Example} \label{202509301321}
    More examples arise from some results of the following sections:
    see Theorem \ref{202512231657}, and the case $p \neq q$ in Theorem \ref{202510171046}.
    Moreover, if $\mathds{k} = \mathds{Q}$, then the opposite category of free Burnside groups with a given exponent (see Example \ref{202609221743}) provides an example.
\end{Example}

In the rest of this section, we give the proof of Theorem \ref{202512241009}.
One may observe that, for a monoid $M$, $\mathfrak{Q}_1 (M)$ is the universal abelian group associated with $M$.
This can be verified directly; alternatively, one may invoke \cite[Lemma 6.5]{hartl2015polynomial} for a different proof, since, for a group $G$, $\mathfrak{Q}_1 (G)$ is the abelianization of $G$.

\begin{Lemma} \label{202512251348}
    For a monoid $M$, we have a natural isomorphism $\mathds{k} \otimes_{\mathds{Z}} \mathfrak{Q}_1 (M) \cong \Aug (M)/ \Aug (M)^{2}$.
\end{Lemma}
\begin{proof}
    The map $\Aug_{\mathds{Z}} (M) \to \Aug (M)$ induces a natural surjection $\mathds{k} \otimes_{\mathds{Z}} \mathfrak{Q}_1 (M) \to \Aug (M)/ \Aug (M)^{2}$.
    On the other hand, $\Aug(M)$ is a free $\mathds{k}$-module generated by $(g-e)$ for $g \in M \backslash \{e\}$.
    This observation yields a map $\Aug (M) \to \mathds{k} \otimes_{\mathds{Z}} \mathfrak{Q}_1 (M)$ which induces the inverse of the above map.
\end{proof}

\vspace{1mm}
\begin{proof}[Proof of Theorem \ref{202512241009}]
    By Theorems \ref{202510011337} and \ref{202408011134}, the condition (1) of the theorem is equivalent to 
    $$
    \Aug (\mathcal{C}_n^{\times m})^2 = \Aug (\mathcal{C}_n^{\times m}) , \quad n,m \in \mathds{N} .
    $$
    By Lemma \ref{202512251348} applied to the monoid $M = \mathcal{C}_n^{\times m}$, this is equivalent to
    \begin{align} \label{202512251451}
        \mathds{k} \otimes_{\mathds{Z}} \mathfrak{Q}_1 (\mathcal{C}_n^{\times m}) \cong 0 , \quad n,m \in \mathds{N}.
    \end{align}
    Hence, all that remain is to derive (\ref{202512251451}) from the condition $\mathds{k} \otimes_{\mathds{Z}} \mathfrak{Q}_1 (\mathcal{C}_1) \cong 0$. 
    By the natural isomorphism $\mathfrak{Q}_1 (M \times N)\cong \mathfrak{Q}_1 (M) \oplus \mathfrak{Q}_1 (N)$ for monoids $M,N$, it suffices to show $\mathds{k} \otimes_{\mathds{Z}} \mathfrak{Q}_1 (\mathcal{C}_n) \cong 0,~ n \in \mathds{N}$ in place of (\ref{202512251451}).
    To this end, we assume that $\mathds{k} \otimes_{\mathds{Z}} \mathfrak{Q}_1 (\mathcal{C}_1) \cong 0$.
    Note that the $\mathds{k}$-module $\mathds{k} \otimes_{\mathds{Z}} \mathfrak{Q}_1 (\mathcal{C}_n)$ is generated by the images of the induced maps
    $$
    p_{k}^\ast : \mathds{k} \otimes_{\mathds{Z}} \mathfrak{Q}_1 (\mathcal{C}_1) \to 
    \mathds{k} \otimes_{\mathds{Z}} \mathfrak{Q}_1 (\mathcal{C}_n)
    $$
    for $1 \leq k \leq n$, since the Lawvere theory $\mathcal{C}$ fulfills the condition (M*).
    Thus, by the hypothesis, we obtain $\mathds{k} \otimes_{\mathds{Z}} \mathfrak{Q}_1 (\mathcal{C}_n) \cong 0$.

    The last assertion follows from Theorem \ref{202609142120}.
\end{proof}

\subsection{Proof of Theorem \ref{202408011134}}
\label{202404121704}

In this section, we prove Theorem \ref{202408011134}.

\begin{notation}
    Let $A$ be a $\mathds{k}$-algebra with unit $1_A$.
    For $n \in \mathds{N}$, $a_0,a_1, \cdots , a_n \in A$ and $T \subset [n]$, we introduce $\prod_{i \in T} a_i {:=} a_{i_1} a_{i_2} \cdots a_{i_r} \in A$ where $i_1 < i_2 < \cdots < i_r$ are elements of $T \neq \emptyset$; and $\prod_{i \in T} a_i {:=} 1_A$ if $T = \emptyset$.   
\end{notation}

\begin{Lemma} \label{202508192133}
    Let $N \in \mathds{N}$.
    For monoids $G,H$, consider two families of monoid homomorphisms $\{ l_i : G \to H\mid i \in [N] \}$ and $\{ F_T : H \to G \mid T \subset [N] \}$.
    We assume the following:
    \begin{enumerate}
        \item $\bigcup^{N}_{i=0} l_i (G)$ generates the monoid $H$.
        \item $F_T \circ l_i = \mathrm{id}_{G}$ if $i \in T$; and $F_T \circ l_i = e$ otherwise.
    \end{enumerate}
    Then $\Aug(G)^{N+1}$ coincides with the image of the map
    $$\sum_{T \subset [N]} (-1)^{N+1-|T|} F_T : \mathds{k}[H] \to \mathds{k}[G].$$
\end{Lemma}
\begin{proof}
    Let $V \subset \mathds{k}[G]$ be the image of $\sum_{T \subset [N]} (-1)^{N+1-|T|} F_T$.
    Since $\Aug(G)^{N+1}$ is generated by $\prod_{i\in[N]} (g_i -e)$ as a $\mathds{k}$-module for $g_i \in G$, we have $\Aug(G)^{N+1} \subset V$.
    In fact, (2) implies the following:
    $$
    \prod_{i\in[N]} (g_i -e) = \sum_{T \subset [N]} (-1)^{N+1-|T|} \prod_{i \in T}g_i = \sum_{T \subset [N]} (-1)^{N+1-|T|} F_T(l_0 (g_0) \cdots l_N (g_N)) .
    $$
    
    We now prove that $V \subset \Aug(G)^{N+1}$.
    By (1), every element of $H$ is given by $\prod_{k\in [r]} \prod_{i\in[N]} l_i (g_i^{(k)})$ for some $r\in\mathds{N}$ and a family $g_i^{(k)} \in G,~k\in [r],i\in [N]$.
    For such an element, we put $x_i^{(k)} {:=} g_i^{(k)} -1 \in \mathds{k}[G]$ where $1$ denotes the unit of $G$.
    By (2), for $T \subset [N]$, we have:
    \begin{align*}
        F_T \left( \prod_{k\in[r]} \prod_{i\in[N]} l_i (g_i^{(k)}) \right) = \prod_{k\in[r]} \prod_{i\in T} g_i^{(k)}  = \prod_{k\in[r]} \prod_{i\in T} (1 + x_i^{(k)}) ,\\
        = \prod_{k\in[r]} \left( \sum_{S \subset T} \left(\prod_{i\in S} x_i^{(k)} \right) \right)
        =  \sum_{S_0, \cdots,S_r \subset T}  \prod_{k\in[r]} \prod_{i\in S_k} x_i^{(k)} .
    \end{align*}
    This yields the following:
    \begin{align*}
        \sum_{T \subset [N]} (-1)^{N+1-|T|} F_T \left( \prod_{k\in[r]} \prod_{i\in[N]} l_i (g_i^{(k)}) \right) &= \sum_{T \subset [N]}\left( \sum_{S_0, \cdots,S_r \subset T} (-1)^{N+1-|T|} \prod_{k\in[r]} \prod_{i\in S_k} x_i^{(k)} \right)  , \\
        &= \sum_{S_0,\cdots, S_r \subset [N]} \left( \sum_{\bigcup^{r}_{k=0} S_k \subset T\subset [N]} (-1)^{N+1-|T|} \prod_{k\in[r]} \prod_{i\in S_k} x_i^{(k)} \right), \\
        &= \sum_{S_0,\cdots, S_r \subset [N]} \left( \sum_{\bigcup^{r}_{k=0} S_k \subset T \subset [N]} (-1)^{N+1-|T|} \right)\left( \prod_{k\in[r]} \prod_{i\in S_k} x_i^{(k)}  \right) .
    \end{align*}
    Note that for $S_0,\cdots, S_r \subset [N]$,
    \begin{align*}
        \sum_{\bigcup^{r}_{k=0} S_k \subset T \subset [N]} (-1)^{N+1-|T|}
        =
        \begin{cases}
            1 , & \mathrm{if~} \bigcup^{r}_{k=1} S_k = [N] , \\
            0 , & \mathrm{otherwise.}
        \end{cases}
    \end{align*}
    By the above results,
    $$
    \sum_{T \subset [N]} (-1)^{N+1-|T|} F_T \left( \prod_{k\in[r]} \prod_{i\in[N]} l_i (g_i^{(k)})  \right) = \sum_{\substack{S_0,\cdots, S_r \subset [N]\\ \bigcup^{r}_{k=0} S_k =[N]}} \prod_{k\in[r]} \prod_{i\in S_k} x_i^{(k)} ,
    $$
    and the right hand side lies in $\Aug (G)^{N+1}$,
    since $x_i^{(k)} \in \Aug(G)$ and $\sum^{r}_{k=0}|S_k| \geq N+1$.
\end{proof}

\begin{proof}[Proof of Theorem \ref{202408011134}]
Let $n,m\in \mathds{N}$.
For the proof, we recall the morphisms introduced in Definition \ref{202403041329} for $X = n \in \mathds{N}$.
We then have monoid maps $(D_n^T)^\ast : \mathcal{C}_{n(d+1)}^{\times m} \to \mathcal{C}_n^{\times m}$ for $T \subset [d]$; and $(p_{k+1})^\ast : \mathcal{C}_n^{\times m} \to \mathcal{C}_{n(d+1)}^{\times m}$ for $k \in [d]$.
By Theorem \ref{202609281707}, $\mathtt{I}^{(d)}_{\mathcal{C}}(m,n)$ equals the image of the map
$$(\pi^{d}_{n})^\ast = \sum_{T \subset [d]} (-1)^{d+1-|T|} (D^T_n)^\ast : \mathds{k}[\mathcal{C}_{n(d+1)}^{\times m} ] \to \mathds{k} [\mathcal{C}_{n}^{\times m}],$$
where $(-)^\ast$ denotes the precomposition as morphisms in $\mathcal{C}$.
On the other hand, we can apply Lemma \ref{202508192133} to 
$$
N=d,~G= \mathcal{C}_n^{\times m},~ H =  \mathcal{C}_{n(d+1)}^{\times m},~ l_i = (p_{i+1})^\ast,~ F_T = (D^T_n)^\ast .
$$
Indeed, the assumptions of the lemma are checked, as the first one follows from the assumption (M*); and the second one follows from (\ref{202407311833}).
Therefore, the image of $(\pi^{d}_{n})^\ast : \mathds{k}[\mathcal{C}_{n(d+1)}^{\times m} ] \to \mathds{k} [\mathcal{C}_{n}^{\times m}]$ coincides with $\Aug(\mathcal{C}_n^{\times m})^{d+1}$.
\end{proof}

\section{Lawvere theory of free $\mathcal{R}$-semisimple groups}
\label{202510071642}

In this section, we introduce the notion of a radical functor for groups, and the associated Lawvere theory.
This framework yields numerous examples to which the results of Section \ref{202510161122} apply.
Examples include the categories considered in the main results.

\subsection{Radical functors for groups} \label{202510211610}

In this section, we introduce the notion of a radical functor for groups.
For $\mathcal{R}$ a radical functor for groups, we also introduce the notion of an $\mathcal{R}$-semisimple group.
We will prove that the opposite category of free $\mathcal{R}$-semisimple groups is a Lawvere theory satisfying condition (ZM*) introduced in Section \ref{202509031755}.

\begin{Defn} \label{202509041252}
    Let $\mathsf{Gr}$ denote the category of all groups.
    A functor $\mathcal{R} : \mathsf{Gr}\to \mathsf{Gr}$ is a {\it radical functor for groups} if the following conditions hold:
    \begin{enumerate}
        \item It is a subfunctor of the identity functor on $\mathsf{Gr}$. In other words, we have $\mathcal{R}(G) \subset G$ and, for a group homomorphism $\rho : G \to H$, we have $\rho ( \mathcal{R} (G)) \subset \mathcal{R} (H)$.
        \item $\mathcal{R}(G) \subset G$ is a normal subgroup.
        \item For a group $G$, we have $\mathcal{R} (G/\mathcal{R}(G)) = 1$.
    \end{enumerate}
\end{Defn}

We will give numerous examples in the next section.

\begin{remark}
    We adopt this terminology from \cite[Section 1]{casacuberta1999singly}, and this concept originates in the systematic study of radicals in group theory.
\end{remark} 

\begin{Defn}
    Let $\mathcal{R}$ be a radical functor for groups.
    A group $G$ is {\it $\mathcal{R}$-semisimple} if $\mathcal{R}(G) = 1$.
\end{Defn}

\begin{Lemma} \label{202609221120}
    Let $\mathcal{R}$ be a radical functor for groups and let $H$ be an $\mathcal{R}$-semisimple group. 
    Every homomorphism $G \to H$ factors uniquely through $G \to G/\mathcal{R}(G)$.
    Consequently, $G \mapsto G/\mathcal{R}(G)$ is left adjoint to the inclusion of $\mathcal{R}$-semisimple groups into $\mathsf{Gr}$, and hence preserves coproducts.
\end{Lemma}
\begin{proof}
    For $\rho : G \to H$, condition (1) gives $\rho(\mathcal{R}(G)) \subset \mathcal{R}(H) = 1$. Condition (3) ensures that $G/\mathcal{R}(G)$ is $\mathcal{R}$-semisimple.
\end{proof}

By the lemma, for an $\mathcal{R}$-semisimple group $H$ we have $\mathsf{Gr}(\mathsf{F}_n/\mathcal{R}(\mathsf{F}_n), H) \cong H^{\times n}$.
Accordingly, $\mathsf{F}_n/\mathcal{R}(\mathsf{F}_n)$ is the free $\mathcal{R}$-semisimple group of rank $n$.

The following compatibility with products will be used in Section \ref{202510211825}.
\begin{Lemma} \label{202510131400}
    Let $\mathcal{R}$ be a radical functor for groups.
    We have $\mathcal{R}(G \times H) = \mathcal{R}(G) \times \mathcal{R}(H)$.
\end{Lemma}
\begin{proof}
     Via the functoriality of $\mathcal{R}$, the projection $p_1 : G \times H \to G$ gives $p_1 ( \mathcal{R}(G \times H)) \subset \mathcal{R}(G)$.
     Similarly, for the other projection $p_2 : G \times H \to H$, we obtain $p_2 ( \mathcal{R}(G \times H)) \subset \mathcal{R}(H)$.
     Therefore, we obtain $\mathcal{R}(G \times H) \subset \mathcal{R}(G) \times \mathcal{R}(H)$.
     On the other hand, the injections $G \to G \times H$ and $H \to G \times H$ induce $\mathcal{R} (G) \times 1 \subset \mathcal{R}(G \times H)$ and $1 \times \mathcal{R} (H)  \subset \mathcal{R}(G \times H)$.
     Since $\mathcal{R}(G \times H)$ is a subgroup of $G \times H$, we obtain $\mathcal{R} (G) \times \mathcal{R}(H) \subset \mathcal{R}(G \times H)$.
\end{proof}

\begin{Defn}
    We denote by $\Rac$ the category of finitely generated free $\mathcal{R}$-semisimple groups and group homomorphisms.
    In this paper, we fix a skeleton of this category by identifying $n \in \mathds{N}$ with the free $\mathcal{R}$-semisimple group of rank $n$, i.e. $\mathsf{F}_n / \mathcal{R} (\mathsf{F}_n)$.
\end{Defn}

\begin{prop}
    Let $\mathcal{R}$ be a radical functor for groups.
    The opposite category $\Rac^{\mathsf{o}}$ is a Lawvere theory that satisfies condition (ZM*).
\end{prop}
\begin{proof}
    We first show that $\Rac$ has finite coproducts, with $n \amalg m = n+m$ and zero object $0$.
    By Lemma~\ref{202609221120}, the coproduct of $\mathsf{F}_n/\mathcal{R}(\mathsf{F}_n)$ and $\mathsf{F}_m/\mathcal{R}(\mathsf{F}_m)$ in the category of $\mathcal{R}$-semisimple groups is
    $$
    (\mathsf{F}_n \ast \mathsf{F}_m)/\mathcal{R}(\mathsf{F}_n \ast \mathsf{F}_m) \cong \mathsf{F}_{n+m}/\mathcal{R}(\mathsf{F}_{n+m}),
    $$
    which is again an object of $\Rac$.
    Since $\Rac$ is a full subcategory of the category of $\mathcal{R}$-semisimple groups, this is also the coproduct in $\Rac$.
    Likewise, the trivial group $\mathsf{F}_0/\mathcal{R}(\mathsf{F}_0)$ is initial and terminal in $\Rac$.
    Hence, the opposite $\Rac^{\mathsf{o}}$ is a Lawvere theory satisfying (Z).

    Next, the quotients $\mathsf{F}_n \to \mathsf{F}_n/\mathcal{R}(\mathsf{F}_n)$ induce a functor $\mathbf{F} \to \Rac$ which is bijective on objects and, by the first paragraph, preserves finite coproducts.
    It is full: a homomorphism $\mathsf{F}_n/\mathcal{R}(\mathsf{F}_n) \to \mathsf{F}_m/\mathcal{R}(\mathsf{F}_m)$ lifts to $\mathsf{F}_n \to \mathsf{F}_m$ since $\mathsf{F}_n$ is free and $\mathsf{F}_m \to \mathsf{F}_m/\mathcal{R}(\mathsf{F}_m)$ is surjective.
    Thus $\mathbf{F}^{\mathsf{o}} \to \Rac^{\mathsf{o}}$ is a full functor that is the identity on objects and preserves finite products.

    Finally, $\mathbf{F}^{\mathsf{o}}$ satisfies (M) and (M*).
    Indeed, the diagonal map $\mathsf{F}_1 \to \mathsf{F}_1 \ast \mathsf{F}_1$ induces a morphism $2 \to 1$ in $\G^{\mathsf{o}}$, which gives a monoid object structure on $1$.
    Furthermore, for $\mathcal{C} = \G^{\mathsf{o}}$, the monoid $\mathcal{C}_n$ coincides with the underlying monoid of the group $\mathsf{F}_n$, which is generated by $\bigcup_{k=1}^{n} \iota_k(\mathsf{F}_1)$, where $\iota_k$ is the inclusion induced by $x_1 \mapsto x_k$; so (M*) holds.
    Via the functor $\G^{\mathsf{o}} \to \Rac^{\mathsf{o}}$ of the previous paragraph, $\Rac^{\mathsf{o}}$ inherits (M) and (M*).
\end{proof}

The following allows us to apply the discussion in Section \ref{202512241413} to such a functor $\mathbf{F}_{\mathcal{R}_1}^{\mathsf{o}} \to \mathbf{F}_{\mathcal{R}_2}^{\mathsf{o}}$:
\begin{prop}
    Let $\mathcal{R}_1,\mathcal{R}_2$ be radical functors for groups such that $\mathcal{R}_1 \subset \mathcal{R}_2$.
    The quotient functor $\mathbf{F}_{\mathcal{R}_1}^{\mathsf{o}} \to \mathbf{F}_{\mathcal{R}_2}^{\mathsf{o}}$ is full and preserves finite products.
\end{prop}
\begin{proof}
    By definitions, for $\mathcal{C} = \Rac^{\mathsf{o}}$, we have $\mathcal{C}(n,m) = \mathcal{C}_n^{\times m} \cong (\mathsf{F}_n /\mathcal{R}(\mathsf{F}_n))^{\times m}$.
    The fullness follows from the surjectivity of the map $(\mathsf{F}_n /\mathcal{R}_1(\mathsf{F}_n))^{\times m} \to (\mathsf{F}_n /\mathcal{R}_2(\mathsf{F}_n))^{\times m}$.
    The rest follows from the definition of finite products.
\end{proof}

\subsection{Examples}

In this section, we present concrete examples of radical functors together with machinery for producing them.
Most examples of our interest satisfy the following:
\begin{Defn} \label{202510091828}
    A radical functor for groups $\mathcal{R}$ is {\it surjection-preserving} if, for a surjective group homomorphism $\rho : G \to H$, we have $\rho ( \mathcal{R} (G)) = \mathcal{R} (H)$.
\end{Defn}

\begin{notation}
    For a group $G$ and $g,h\in G$, we write $[g,h] {:=} g^{-1}h^{-1}gh$.
    For subgroups $H_1,H_2 \subset G$, we denote by $[H_1,H_2]$ the {\it commutator subgroup} that is generated by $[a,b]$ for $a \in H_1$ and $b \in H_2$.
\end{notation}

\begin{Example} \label{202609221743}
    Let $F$ be a free group and $w \in F$.
    For a group $G$, the verbal subgroup of $G$ induced by $w$ is defined as 
    $$
    \mathcal{R} (G) = \langle \alpha (w) \in G \mid \alpha : F \to G \rangle
    $$
    where $\alpha$ runs over all group homomorphisms.
    One may show that it gives a surjection-preserving radical functor for groups.
    Examples include the following:
    \begin{enumerate}
        \item Consider $w = x_1^{r} \in \mathsf{F}_1$. 
        The associated radical functor is given by $\mathcal{R} (G) = G^r$, so the associated category $\Rac$ consists of finitely generated free Burnside groups with respect to the exponent $r$.
        \item One may extend the notion of a verbal subgroup by allowing a collection of $w$'s.
        For instance, given $r \in \mathds{N}^\ast$, $\mathcal{R} (G) = \langle g^r , [g,h] \mid g ,h \in G\rangle$ gives a radical functor for groups.
        In this case, a free $\mathcal{R}$-semisimple group is given by a product of $\mathds{Z}/r$, so the associated category $\Rac$ coincides with the category $\mathbf{M}_{\mathds{Z}/r}$ considered in Example \ref{202512281213}.
    \end{enumerate}
\end{Example}

\begin{prop} \label{202510111100}
    Let $\mathcal{R}_1, \mathcal{R}_2$ be radical functors for groups that are surjection-preserving.
    The following yields a surjection-preserving radical functor for groups:
    $$
    [\mathcal{R}_1, \mathcal{R}_2] (G) {:=} [\mathcal{R}_1(G), \mathcal{R}_2 (G)] .
    $$
\end{prop}
\begin{proof}
    Since the commutator subgroup of normal subgroups is normal, so is $[\mathcal{R}_1, \mathcal{R}_2] (G) \subset G$.
    We show that $[\mathcal{R}_1,\mathcal{R}_2] (G/N)$ is trivial where $N= [\mathcal{R}_1,\mathcal{R}_2](G)$.
    By definition, $[\mathcal{R}_1,\mathcal{R}_2](G/N)$ is generated by $[g_1N, g_2N]$ for $g_1N \in \mathcal{R}_1 (G/N)$ and $g_2N \in \mathcal{R}_2 (G/N)$.
    Since $\mathcal{R}_1$ and $\mathcal{R}_2$ are surjection-preserving, we can choose $g_i$ to be $g_i \in \mathcal{R}_i (G)$.
    We then have $[g_1,g_2] \in N$.
    Hence, $[g_1N, g_2N] = [g_1,g_2]N =N$, so $[\mathcal{R}_1,\mathcal{R}_2](G/N)$ is trivial.
    Consequently, $[\mathcal{R}_1, \mathcal{R}_2]$ is a radical functor for groups.
    Furthermore, $[\mathcal{R}_1, \mathcal{R}_2]$ is surjection-preserving, since for a surjective group homomorphism $\rho : G \to H$, we have
    $$\rho (  [\mathcal{R}_1, \mathcal{R}_2] (G) ) = \rho ([\mathcal{R}_1(G), \mathcal{R}_2 (G)] ) = [\rho (\mathcal{R}_1(G)), \rho (\mathcal{R}_2 (G))]  = [\mathcal{R}_1 (H), \mathcal{R}_2(H) ] = [\mathcal{R}_1, \mathcal{R}_2] (H).$$
\end{proof}

\begin{Example} \label{202510121401}
    Obviously, the identity functor $\mathrm{Id}_{\mathsf{Gr}}$ on $\mathsf{Gr}$ is a surjection-preserving radical functor for groups.
    By Proposition \ref{202510111100}, every functor obtained by recursive commutators of $\mathrm{Id}_{\mathsf{Gr}}$ is a surjection-preserving radical functor for groups.
    There are two typical examples: 
    \begin{enumerate}
        \item The $c$-th term of the lower central series $\gamma_c$ that is recursively defined by $\gamma_1(G) = G$ and $\gamma_{c+1}(G) = [\gamma_{c}(G), G]$.
        The category $\Rac$ for $\mathcal{R} = \gamma_{c+1}$ consists of finitely generated free nilpotent groups of class $\leq c$.
        \item The $c$-th derived subgroup $\gamma_{2}^{c} = \overbrace{\gamma_{2}\circ \cdots \circ \gamma_{2}}^{c}$.
        This can be recursively defined by $\gamma_{2}^{1} (G)= [G,G]$ and $\gamma_{2}^{c+1} (G)= [\gamma_2^{c} (G) , \gamma_2^{c} (G)] $.
        In this case, $\Rac$ is the category of finitely generated free solvable groups of derived length $\leq c$.
    \end{enumerate}
\end{Example}

\begin{notation}
    In this paper, we set $\mathbf{N}_{c} {:=} \G_{\gamma_{c+1}}$.
    By definition, $\mathbf{N}_{c}$ is the category of finitely generated free nilpotent groups of class at most $c$.
\end{notation}

For a group $G$ and a family of normal subgroups $H_i \subset G,~i \in I$, let $\prod_{i\in I} H_i$ denote the smallest subgroup containing all the $H_i$.
Similarly, one may prove the following:
\begin{prop} \label{202510161130}
    Let $\mathcal{R}_i,~i\in I$ be a family of radical functors for groups that are surjection-preserving.
    Then $(\prod_{i\in I} \mathcal{R}_i)(G) {:=} \prod_{i \in I} \mathcal{R}_i(G)$ yields a surjection-preserving radical functor for groups.
\end{prop}

\vspace{1mm}
We also present another method which does not necessarily produce surjection-preserving radical functors.
Let $F$ be a free group and let $f \in \mathds{Z}[F]$.
For a unital commutative ring $R$ and a group $G$, the {\it polynomial ideal} \cite{MR537126}, denoted by $A_{f,R}(G)$, is defined as the two-sided ideal of $R[G]$ generated by $\alpha_\ast (f)\in R[G]$ for group homomorphisms $\alpha : F \to G$.
    
\begin{prop}
    The assignment $\mathcal{R} (G) = G \cap ( 1 + A_{f,R}(G) )$ to a group $G$ gives a radical functor for groups.
\end{prop}
\begin{proof}
    The only nontrivial part is proving the condition that $\mathcal{R}(G/\mathcal{R}(G)) \cong 1$.
    Denote by $\pi : G \to G/ \mathcal{R}(G)$ the projection.
    Let $g \in G$ such that $\pi(g) \in \mathcal{R}(G/\mathcal{R}(G))$.
    Then $\pi(g)-1 \in A_{f,R} (G/\mathcal{R}(G))$.
    By definition, there exist $a_j,b_j \in R[G]$, group homomorphisms $\alpha_j : F \to G/\mathcal{R}(G)$ such that $\pi (g)-1 = \sum_{j} \pi_\ast (a_j) (\alpha_j)_\ast (f) \pi_\ast (b_j )$.
    Since $F$ is free, $\alpha_j$ lifts to a homomorphism $\beta_j : F \to G$.
    In particular, $\pi (g) -1 = \pi_\ast ( \sum_j a_j (\beta_j)_\ast (f) b_j)$, so there exist $z_i \in G$ and $u_i \in \mathcal{R}(G)$ such that
    $$
    g- 1 = \sum_j a_j (\beta_j)_\ast (f) b_j + \sum_i z_i (u_i -1) .
    $$
    Since $u_i \in \mathcal{R}(G)$, $(u_i -1) \in A_{f,R}(G)$ which leads to $g-1 \in A_{f,R}(G)$.
    Hence, $g \in \mathcal{R}(G)$, so $\pi(g) =1$.
\end{proof}

\begin{Example} \label{202510121414}
    The radical functor associated with $f = (x_1-1) \cdots (x_c-1) \in \mathds{Z}[\mathsf{F}_{c}]$ in the proposition yields the $c$-th {\it dimension subgroup} functor, which we denote by $D_{c,R}$.
    It is well known that, if $R$ is a field, then $D_{c,R}$ is determined by the characteristic $p$ of $R$ (see \cite{MR537126} or \cite[Section 4]{massuyeau2007finite}). 
    We set 
    \begin{align*}
        D_{c,p} {:=} D_{c,\mathds{Q}},~\mathrm{if~}p=0; \quad \mathrm{and~~}D_{c,p} {:=} D_{c,\mathds{F}_{p}},~\mathrm{if~}p>0 .
    \end{align*}
    The following can be found in the literature:
    \begin{itemize}
        \item If $p > 0$, then
        $$
        D_{c,p} (G) = \prod_{ip^j \geq c} \gamma_{i}(G)^{p^j} .
        $$
        Since $\gamma_i$'s are surjection-preserving, so $D_{c,p}$ is by Proposition \ref{202510161130}.
        \item If $p = 0$, then
        $$
        D_{c,0} (G) = \sqrt{\gamma_{c}(G)} .
        $$
        For a group $G$ and a normal subgroup $H \subset G$, $\sqrt{H}$ denotes the subset of $G$ consisting of $g \in G$ such that $g^n \in H$ for some $n \in \mathds{N}^\ast$.
        If $\gamma_{c+1}(G) \subset H$ for some $c$, then $\sqrt{H}$ is a normal subgroup of $G$.
        It follows from the fact that, in a nilpotent group, the set of finite order elements yields a subgroup.
        For $c \geq 2$, $D_{c,0}$ is not surjection-preserving, as $\sqrt{\gamma_c (\mathds{Z})} = 0$ while $\sqrt{\gamma_c (\mathds{Z}/n)} = \mathds{Z}/n$.
    \end{itemize}
\end{Example}

\begin{notation} \label{202512281222}
    For a prime $p$ or $p = 0$, let $\Dim_{c,p} {:=} \Rac$ where $\mathcal{R} = D_{c+1,p}$.
\end{notation}

\begin{remark} \label{202510221959}
    For a free group $F$, we have $D_{c,0}(F) = \sqrt{\gamma_{c}(F)} = \gamma_{c}(F)$.
    Hence, $\Dim_{c,0} = \mathbf{N}_{c}$.
\end{remark}

\subsection{A computable criterion for $\mathrm{D} (\Rac^{\mathsf{o}})$}
\label{202510271047}

In this section, we investigate the strictness of polynomial degree filtration given in (\ref{202509291620}) applied to $\mathcal{C} = \Rac^{\mathsf{o}}$.
In application, we prove Theorems \ref{202510171045} and \ref{202510171046}.

\begin{Defn}
    Let $\mathcal{R}$ be a radical functor for groups.
    We define $|\mathcal{R}| \in \mathds{N}$ as $r \in \mathds{N}$ such that $$\mathcal{R}(\mathds{Z}) = r \mathds{Z}.$$
\end{Defn}

\begin{Example}
    We have $|\gamma_{c}| = 0$ for $c>1$ and $|\gamma_{1}| = 1$.
    By the arguments in Example \ref{202510121414}, if $p$ is a prime, then $|D_{c,p}| = p^j$ where $j$ is the least integer such that $c \leq p^j$.
\end{Example}

\begin{theorem} \label{202512231657}
    Suppose that $\mathds{k} \neq 0$.
    For $\mathcal{R}$ a radical functor for groups and $r = |\mathcal{R}|$, the following are equivalent:
    \begin{enumerate}
        \item $r \cdot 1_\mathds{k} \in \mathds{k}$ is invertible.
        \item $\mathrm{D} ( \Rac^{\mathsf{o}}) =\{ 0\}$.
    \end{enumerate}
    Otherwise, $\mathrm{D} ( \Rac^{\mathsf{o}}) = \mathds{N}$.
\end{theorem}
\begin{proof}
    Let $\mathcal{C} = \Rac^{\mathsf{o}}$.
    By definition, $\mathcal{C}_1 = \mathds{Z}/r$ (see Notation \ref{202512241120}), so $\mathfrak{Q}_1 (\mathcal{C}_1) \cong \mathds{Z}/r$.
    Theorem \ref{202512241009} shows that $\mathrm{D}(\mathcal{C}) = \{ 0\}$ if and only if
    $\mathds{k}\otimes_{\mathds{Z}} \mathds{Z}/r \cong \mathds{k}/ r~ \mathds{k} \cong 0$.
    The latter condition is equivalent to saying that $r \cdot 1_\mathds{k}$ is invertible.
\end{proof}
    
\begin{remark}
    This theorem implies that, if $|\mathcal{R}|=0$, then the inclusion functor  $\mathtt{L}_{\Rac^{\mathsf{o}}}\mbox{-}\mathsf{Mod}^{\leq d} \hookrightarrow \mathtt{L}_{\Rac^{\mathsf{o}}}\mbox{-}\mathsf{Mod}^{\leq (d+1)}$ is not essentially surjective.
    In our sequel paper \cite{kimprop}, we study the Serre quotient of these abelian categories.
\end{remark}

Applying the above theorem, we now prove Theorem \ref{202510171045}:
\begin{proof}[Proof of Theorem \ref{202510171045}]
    Let $c \in \mathds{N}^\ast$.
    Theorem \ref{202512231657} applied to $\mathcal{R} = \gamma_{c+1}$ proves Theorem \ref{202510171045}, since $|\gamma_{c+1}|=0$.
    We obtain the same result for $\mathrm{D} (\mathbf{N}_{c})$ from Proposition \ref{202512231721}.
\end{proof}

Likewise, we give the proof of Theorem \ref{202510171046}.

\begin{proof}[Proof of Theorem \ref{202510171046}]

To prove the statement, we apply Theorem \ref{202512231657} to the case $\mathcal{R} = D_{c+1,q}$.
Note that $|D_{c+1,q}|$ is a power of $q$.
If $p \neq q$, then $|D_{c+1,q}|$ is invertible in $\mathds{k}$.
Hence, we obtain $\mathrm{D} (\Dim_{c,q}^{\mathsf{o}}) = \{0\}$.
Likewise, if $p = q$, then Theorem \ref{202512231657} proves $\mathrm{D} (\Dim_{c,q}^{\mathsf{o}}) = \mathds{N}$.
Finally, Proposition \ref{202512231721} yields the same result for $\mathrm{D} (\Dim_{c,q})$.
\end{proof}

\subsection{Radical functors on quotient groups}

In this section, we describe the quotients $\mathcal{R}_2(G/H)/\mathcal{R}_1(G/H)$ of radical functors $\mathcal{R}_1 \subset \mathcal{R}_2$ evaluated on a quotient group $G/H$ in terms of $G$ and $H$.
These results will be used to prove Theorems \ref{202509021231} and \ref{202510021437}.

\begin{Lemma} \label{202510031417}
    Let $\mathcal{R}_1 \subset \mathcal{R}_2$ be radical functors for groups.
    If $\mathcal{R}_1$ and $\mathcal{R}_2$ are surjection-preserving, then, for a group $G$ and a normal subgroup $H \subset G$, the map $\mathcal{R}_2(G) \to \mathcal{R}_2(G/H)$ induces an isomorphism of groups
    $$
    \mathcal{R}_2 (G)/ \left( \mathcal{R}_1 (G) (\mathcal{R}_2(G)\cap H) \right)  \to \mathcal{R}_2 (G/H)/\mathcal{R}_1 (G/H) .
    $$
\end{Lemma}
\begin{proof}
    It is routine to check that the map in the statement is well-defined.
    Since $\mathcal{R}_2$ is surjection-preserving, $\mathcal{R}_2 (G) \to \mathcal{R}_2 (G/H)$ is surjective.
    Hence, the map in the statement is surjective.
    We prove that it is injective.
    Let $N = \mathcal{R}_1 (G) (\mathcal{R}_2(G)\cap H)$.
    Let $x N \in \mathcal{R}_2 (G) / N$ be an element in the kernel of the map.
    We then have $x H \in \mathcal{R}_1 ( G /H)$.
    Since $\mathcal{R}_1$ is surjection-preserving, there exists $y \in \mathcal{R}_1 (G)$ such that $x H = y H$.
    On the other hand, $y^{-1}x \in \mathcal{R}_1(G) \mathcal{R}_2(G) \subset \mathcal{R}_2(G)$, so $y^{-1}x \in \mathcal{R}_2(G) \cap H$.
    Hence, $x = y (y^{-1}x)  \in \mathcal{R}_1 (G) (\mathcal{R}_2(G)\cap H) = N$.
\end{proof}

Even though $D_{c,0} = \sqrt{\gamma_{c}}$ is not surjection-preserving, as explained in Example \ref{202510121414}, we have an analogous statement for $\sqrt{\gamma_{c}}$:
\begin{Lemma} \label{202510091525}
    Let $G$ be a group with a normal subgroup $H \subset G$.
    Then the natural map $\sqrt{\gamma_{c}(G)H} \to \sqrt{\gamma_{c}(G/H)}$ induces an isomorphism
    $$
    \sqrt{\gamma_{c}(G)H} / \sqrt{\gamma_{c+1}(G)H} \to \sqrt{\gamma_{c}(G/H)} / \sqrt{\gamma_{c+1}(G/H)}.
    $$
\end{Lemma}
\begin{proof}
    We first show that the map $\sqrt{\gamma_{c}(G)H} \to \sqrt{\gamma_{c}(G/H)}$ is surjective.
    Let $xH \in \sqrt{\gamma_{c}(G/H)}$ with $n \in \mathds{N}^\ast$ such that $x^nH \in \gamma_{c}(G/H)$.
    Since $\gamma_{c}$ is surjection-preserving, there exists $y \in \gamma_{c}(G)$ such that $x^nH = yH$.
    In other words, we have $y^{-1}x^n \in H$.
    Hence, $x^n = y(y^{-1}x^n) \in \gamma_{c}(G) H$, so $x \in \sqrt{\gamma_c(G)H}$.

    We now see that the map in the statement is surjective.
    It suffices to show that this is injective.
    Suppose that $x \in \sqrt{\gamma_{c}(G)H}$ satisfies $xH \in \sqrt{\gamma_{c+1}(G/H)}$.
    Then for some $n \in \mathds{N}$, we have $x^n H \in \gamma_{c+1}(G/H)$.
    Since $\gamma_{c+1}$ is surjection-preserving, we can take $y \in \gamma_{c+1}(G)$ such that $x^n H = yH$.
    As above, we then obtain $x^n = y(y^{-1}x^n) \in \gamma_{c+1}(G)H$.
    Hence, $x \in \sqrt{\gamma_{c+1}(G)H}$.
\end{proof}

\section{Recollection of some N-series}
\label{202510231558}

This section provides a recollection of some classical structural theorems on the associated graded algebra of a group.
A general reference on this topic is the book by Passi \cite{MR537126}.
We shall apply these results to prove our main theorems in the following sections.

\subsection{N-series}

Let $G$ be a group.
An {\it N-series} of $G$ is a descending filtration of normal subgroups
$$
\cdots \subset N_{c+1} \subset N_{c} \subset \cdots \subset N_1 = G
$$
such that $[N_i,N_j] \subset N_{i+j}$.
For instance, the lower central series and the dimension subgroups over a ring $R$ give N-series of $G$.
It is clear that the quotient group $N_c /N_{c+1}$ is an abelian group, so $\bigoplus_{c\geq 1} N_{c}/N_{c+1}$ gives a non-negatively graded abelian group.

The commutator on the group $G$ induces a graded Lie $\mathds{Z}$-algebra structure on
$$\mathfrak{L}(G;N_{\bullet}) {:=} \bigoplus_{c\geq 1} N_{c}/N_{c+1} .$$
For a free group $F$, $\mathfrak{L} (F; \gamma_{\bullet})$ is known to be a free Lie $\mathds{Z}$-algebra where $\gamma_{\bullet}= \gamma_{\bullet}(F)$ is the lower central series.
Furthermore, the homogeneous component of $\mathfrak{L} (F; \gamma_{\bullet})$ is a free abelian group whose basis can be described by using basic commutators (see, e.g., \cite[Theorem 11.2.4]{PHall1976}).
This proves the following proposition:
\begin{Lemma} \label{202512301529}
    For a free group $F$ of rank at least $2$, each homogeneous component of $\mathfrak{L} (F; \gamma_{\bullet})$ is non-trivial; in particular $\gamma_{c+1}(F) \subsetneq \gamma_{c}(F), ~ c \in \mathds{N}^\ast$.
\end{Lemma}

\subsection{Sandling–Tahara Theorem}
\label{202510091301}

This section presents a review of the Sandling--Tahara theorem \cite{sandling1979augmentation} which relates the lower central factors to the augmentation ideals.

We start with a general definition.
For an abelian group $A$, let $S(A)$ be the symmetric algebra of $A$.
It is given by the quotient algebra of the tensor $\mathds{Z}$-algebra of $A$ by the two-sided ideal generated by $a \otimes b - b \otimes a$ for $a,b \in A$.
Let $S^k (A) \subset S(A)$ be the $k$-th symmetric power of $A$ with $S^0 (-) {:=} \mathds{Z}$.

If $A$ is a graded abelian group, i.e. $A = \bigoplus_{i \in \mathds{Z}} A_i$, then the symmetric power $S(A)$ inherits a natural grading $S(A) = \bigoplus_{d\in\mathds{Z}} \mathds{S}_d (A)$ with
\begin{align} \label{202510091635}
\mathds{S}_d (A) = \bigoplus \bigotimes_{i\in\mathds{Z}} S^{a_i} ( A_i ) , \quad d \in \mathds{Z}  ,
\end{align}
where the tensor product is over $\mathds{Z}$, and the direct sum is taken over $a_i \in \mathds{Z}$ such that $\sum_{i\in\mathds{Z}} ia_i =d$ with finitely many nonzero $a_i$'s.

Recall $\Aug_{\mathds{Z}} (G)$ from Notation \ref{202609141203}.
\begin{notation} \label{202509021132}
    For a group (more generally, a monoid) $G$, we set
    \begin{align*}
        \mathfrak{Q}_{d} (G) {:=} \Aug_{\mathds{Z}}(G)^{d} / \Aug_{\mathds{Z}}(G)^{d+1} , \mathrm{~and} \quad \mathfrak{G}_{\mathds{Z}}(G) {:=} \bigoplus_{d \in \mathds{N}} \mathfrak{Q}_{d} (G) .
    \end{align*}
\end{notation}

The Sandling--Tahara theorem states that, for a group $G$ whose lower central factors are free-abelian, there exists an isomorphism of graded abelian groups between the associated graded ring $\mathfrak{G}_{\mathds{Z}}(G)$ and the symmetric algebra $S ( \mathfrak{L}(G; \gamma_{\bullet}))$ of the underlying graded abelian group of $\mathfrak{L}(G; \gamma_{\bullet})$:
\begin{theorem}[\cite{sandling1979augmentation}] \label{202509021146}
    Let $G$ be a group whose lower central factor $\gamma_{c}(G)/\gamma_{c+1}(G)$ is free-abelian for any $c \geq 1$.
    We then have $\mathfrak{Q}_d (G) \cong \mathds{S}_d ( \mathfrak{L}(G; \gamma_{\bullet})), \quad d \geq 1$.
\end{theorem}
Note that this isomorphism depends on a choice of a basis of $\gamma_c(G)/\gamma_{c+1}(G)$ for $c \geq 1$.

\subsection{Quillen's theorem}
\label{202510091302}
In this section, we assume that the ground ring $\mathds{k}$ is a field of characteristic $p \geq 0$.
In the following, we give a review of Quillen's theorem \cite{Quillen1968} on the associated graded $\mathds{k}$-algebra of a group $G$:
$$\mathfrak{G} (G) {:=} \bigoplus_{d \in \mathds{N}} \Aug(G)^{d}/\Aug(G)^{d+1} . $$
We refer to \cite[Theorem 5.2, Chapter VIII]{MR537126} for an alternative treatment, and \cite[Section 5]{massuyeau2007finite} for an extension to general N-series.

Let $p = 0$.
For a Lie $\mathds{k}$-algebra $\mathfrak{g}$, we denote by $\mathcal{U}^{0} (\mathfrak{g})$ the universal enveloping algebra.
It is defined by the quotient $\mathds{k}$-algebra of the tensor $\mathds{k}$-algebra by the two-sided ideal generated by $x \otimes y - y\otimes x - [x,y]$ for $x,y \in \mathfrak{g}$.
Below, a $0$-restricted Lie algebra means a usual Lie algebra, and the $0$-restricted universal enveloping algebra is the universal enveloping algebra.

We consider the case $p > 0$.
Let $\mathfrak{g}$ be a $p$-restricted Lie algebra over $\mathds{k}$.
Denote by $x \mapsto x^{[p]}$ the associated {\it $p$-th power} map.
The $p$-restricted universal enveloping algebra of $\mathfrak{g}$ is defined by the quotient $\mathds{k}$-algebra of the tensor $\mathds{k}$-algebra by the two-sided ideal generated by $x \otimes y - y\otimes x - [x,y]$ and $x^{\otimes p} - x^{[p]}$ for $x,y \in \mathfrak{g}$.
We denote this by $\mathcal{U}^{p}(\mathfrak{g})$.

Let $p \geq 0$.
A Lie algebra $\mathfrak{g}$ is {\it non-negatively graded} if $\mathfrak{g}$ is endowed with $\mathfrak{g}=\bigoplus_{n \in \mathds{N}} \mathfrak{g}_n$ such that $[\mathfrak{g}_n,\mathfrak{g}_m]\subset \mathfrak{g}_{n+m}$; and, if $p > 0$, $x^{(p)} \in \mathfrak{g}_{pn}$ for $x \in \mathfrak{g}_n$.
If a $p$-restricted Lie algebra $\mathfrak{g}$ is non-negatively graded, then so is $\mathcal{U}^{p}(\mathfrak{g})$, since, for homogeneous elements $x,y\in\mathfrak{g}$, the generators of the aforementioned ideals are homogeneous.
We denote by $\mathcal{U}^{p}_{c}(\mathfrak{g})$ the $c$-th component which is generated by the product $v_1 v_2 \cdots v_n$ for homogeneous elements $v_j \in \mathfrak{g}$ with degree $n_j$ such that $c = \sum^{n}_{k=1} n_k$.

We recall the notation given in Example \ref{202510121414}.
For $p \geq 0$, the associated graded Lie $\mathds{Z}$-algebra $\mathfrak{L}(G; D_{\bullet,p})$ gives a $p$-restricted Lie algebra.
Quillen established a natural isomorphism of graded $\mathds{k}$-algebras:
\begin{theorem}[\cite{Quillen1968}] \label{202510091518}
    $\mathcal{U}^{p}(\mathds{k} \otimes_{\mathds{Z}} \mathfrak{L}(G; D_{\bullet,p})) \cong \mathfrak{G} (G)$.
\end{theorem}

The following simple observation with respect to the Lie algebras associated with $\gamma_{\bullet}$ and $D_{\bullet,p}$ is probably known, although we could not find an explicit reference.
\begin{prop} \label{202510231737}
    Let $c \in \mathds{N}^\ast$ and $p$ be a prime.
    For a free group $F$, the induced map $\mathds{F}_p \otimes_{\mathds{Z}} \mathfrak{L}(F; \gamma_{\bullet}) \to \mathfrak{L}(F; D_{\bullet,p})$ is injective.
    In particular, if $F$ has rank at least $2$, then $D_{c+1,p}(F) \subsetneq D_{c,p}(F)$.
\end{prop}

For the proof, we formulate the following lemma, which will also be exploited later.

\begin{Lemma} \label{202508291109}
    Let $d_0 \in \mathds{N}$.
    Let $G$ be a group such that $\mathfrak{Q}_d (G)$ is a free abelian group for $ d \leq d_0$.
    Then $\mathds{k} \otimes_{\mathds{Z}} \mathfrak{Q}_{d} (G) \to \Aug ( G)^{d} /\Aug (G)^{d+1}$ is an isomorphism for $d \leq d_0$.
\end{Lemma}
\begin{proof}
    The proof is based on some elementary linear algebra.
    Let $d \leq d_0$.
    The coefficient extension gives a surjection $\mathds{k} \otimes_{\mathds{Z}} \mathfrak{Q}_{d} (G) \to \Aug ( G)^{d} /\Aug (G)^{d+1}$.
    We should show that its kernel, denoted by $A_{d}$, is trivial.
    By the freeness assumption on $\mathfrak{Q}_{d} (G)$, the following exact sequence splits:
    $$
    0 \to \Aug_{\mathds{Z}} (G)^{d+1} \to \Aug_{\mathds{Z}} (G)^{d} \to  \mathfrak{Q}_{d} (G) \to 0 .
    $$
    Hence, taking the tensor product $\mathds{k} \otimes_{\mathds{Z}} (-)$ yields a short exact sequence of $\mathds{k}$-modules.
    Let $K_{d}$ be the kernel of the surjection $\mathds{k} \otimes_{\mathds{Z}} \Aug_{\mathds{Z}} (G)^{d} \to \Aug (G)^{d}$.
    We then obtain the following commutative diagram whose columns and the bottom two rows are exact:
    $$
    \begin{tikzcd}[sep=small]
        & 0 \ar[d] & 0 \ar[d] & 0 \ar[d] & \\
        0 \ar[r] & K_{d+1} \ar[r] \ar[d] & K_{d} \ar[r] \ar[d] & A_{d} \ar[d] \ar[r] &0 \\
        0 \ar[r] & \mathds{k} \otimes_{\mathds{Z}} \Aug_{\mathds{Z}} ( G)^{d+1} \ar[d] \ar[r] & \mathds{k} \otimes_{\mathds{Z}} \Aug_{\mathds{Z}} ( G)^{d}\ar[d] \ar[r] & \mathds{k} \otimes_{\mathds{Z}} \mathfrak{Q}_{d} (G) \ar[d] \ar[r] &0 \\
         0 \ar[r] & \Aug (G)^{d+1}  \ar[r] \ar[d] & \Aug ( G )^{d} \ar[r] \ar[d] & \Aug ( G )^{d}/ \Aug (G)^{d+1} \ar[d] \ar[r] &0 \\
         & 0 & 0 & 0 &
    \end{tikzcd}
    $$
    By the nine-lemma, the first row is exact for any $d \in \mathds{N}$.
    Since $\Aug (G)$ has a basis consisting of $(g-e)$ for $g \in G \backslash \{ e\}$, the map $\mathds{k} \otimes_{\mathds{Z}} \Aug_{\mathds{Z}} (G) \to \Aug (G)$ is an isomorphism, so $K_1 \cong 0$.
    Therefore, it is inductively proved that $K_{d} \cong 0$ and $A_{d} \cong 0$.
\end{proof}

\begin{proof}[Proof of Proposition \ref{202510231737}]
    To prove the statement, we recall a classical result analogous to Theorem \ref{202510091518}.
    Let $\mathcal{U}( \mathfrak{L}(F; \gamma_{\bullet}))$ be the universal enveloping algebra of $\mathfrak{L}(F; \gamma_{\bullet})$ (over $\mathds{Z}$).
    We then have $\mathds{Z}$-algebra isomorphisms between $\mathcal{U}( \mathfrak{L}(F; \gamma_{\bullet}))$, $\mathfrak{G}_{\mathds{Z}}(F)$ and $T ( F/\gamma_2 (F))$ the tensor $\mathds{Z}$-algebra of the abelian group $F/\gamma_2 (F)$ (for instance, see \cite[Theorem 6.2, Chapter VIII]{MR537126}).
    Note that $\mathcal{U}( \mathfrak{L}(F; \gamma_{\bullet}))$ contains $\mathfrak{L}(F; \gamma_{\bullet})$ as a direct summand, since $\mathfrak{L}(F; \gamma_{\bullet}) \hookrightarrow T ( F/\gamma_2 (F))$ can be seen as the subspace consisting of Lie polynomials generated by $F/\gamma_2(F)$ (this follows from the study of $T ( F/\gamma_2 (F))$ using Hall basis \cite{reutenauer2003free}).
    Hence, we obtain the following commutative diagram:
    $$
    \begin{tikzcd}[row sep = small, column sep = small]
        \mathfrak{L}(F; D_{\bullet,p}) \ar[r, hookrightarrow] & \mathcal{U}^{p}(\mathds{F}_{p} \otimes \mathfrak{L}(F; D_{\bullet,p}) ) \ar[r, "\cong"] & \mathfrak{G}_{\mathds{F}_{p}}(F) \\
        \mathds{F}_p \otimes \mathfrak{L}(F; \gamma_{\bullet}) \ar[u] \ar[r, hookrightarrow] & \mathds{F}_{p} \otimes \mathcal{U}( \mathfrak{L}(F; \gamma_{\bullet})) \ar[r, "\cong"] \ar[u] & 
         \mathds{F}_{p} \otimes  \mathfrak{G}_{\mathds{Z}}(F) \ar[u]
    \end{tikzcd}
    $$
    where all the tensor products are defined on $\mathds{Z}$.
    Since every homogeneous subspace of $\mathfrak{G}_{\mathds{Z}}(F)$ is free, by Lemma \ref{202508291109}, the map in the right column is an isomorphism.
    Moreover, the map in the upper left corner is an inclusion by the PBW theorem.
    Therefore, the desired statement follows from the commutativity.
    The final claim follows from Lemma \ref{202512301529} and the above results.
\end{proof}

\section{Polynomial $\mathbf{N}_{c}^{\mathsf{o}}$-modules}
\label{202510211823}

In this section, we study polynomial functors on the category of free nilpotent groups, with a particular focus on their comparison across varying nilpotency levels.
In particular, we give the main argument for the proof of Theorem \ref{202510101544}.

Let $c_0, c_1 \in \mathds{N}^\ast \cup \{\infty\}$ such that $c_0 < c_1$.
Then the full functor $\mathbf{N}_{c_1}^{\mathsf{o}} \to \mathbf{N}_{c_0}^{\mathsf{o}}$ induces a monad epimorphism $\mathtt{L}_{\mathbf{N}_{c_1}^{\mathsf{o}}} \to \mathtt{L}_{\mathbf{N}_{c_0}^{\mathsf{o}}}$.
The following is the main theorem of this section:
\begin{theorem} \label{202509021231}
    Suppose that $\mathds{k} \neq 0$.
    For $c_0, c_1 \in \mathds{N}^\ast \cup \{\infty\}$ such that $c_0 < c_1$, the functor $\mathtt{L}_{\mathbf{N}_{c_0}^{\mathsf{o}}}\mbox{-}\mathsf{Mod}^{\leq d} \to \mathtt{L}_{\mathbf{N}_{c_1}^{\mathsf{o}}}\mbox{-}\mathsf{Mod}^{\leq d}$ gives an equivalence of categories if and only if $0 \leq d \leq c_0$.
    Furthermore, if $d> c_0$, then this is not essentially surjective.
\end{theorem}

Recalling the set $\Gamma$ from Definition \ref{202609211701}, one may obtain $\Gamma(\mathbf{N}_{c_1}^{\mathsf{o}} \to \mathbf{N}_{c_0}^{\mathsf{o}}) = [c_0]$, which leads to Theorem \ref{202510101544} by Corollary \ref{202510141158}.

\subsection{Proof of Theorem \ref{202509021231}}
In this section, we prove Theorem \ref{202509021231}.
We begin with some general observations.
Let $G$ be a product of free nilpotent groups of class $\leq c_0$.
In other words, $G = (F/\gamma_{c_0+1}(F))^{\times m} = F^{\times m}/ \gamma_{c_0+1}(F^{\times m})$ for a free group $F$.
Then we have:
\begin{align} \label{202509011229}
    \gamma_{c} ( G) / \gamma_{c+1} (G) \cong 
    \begin{cases}
        (\gamma_{c} ( F) / \gamma_{c+1} (F))^{\times m} & c \leq c_0 , \\
        1 & \mathrm{otherwise.}
    \end{cases}
\end{align}
This follows from Lemma \ref{202510031417} applied to $\mathcal{R}_2 = \gamma_c$ and $\mathcal{R}_1 = \gamma_{c+1}$.
Hence, the Lie $\mathds{Z}$-algebra $\mathfrak{L} (G; \gamma_{\bullet})$ gives a truncation of $\mathfrak{L} (F^{\times m}; \gamma_{\bullet})$.
In particular, each $\gamma_{c} ( G) / \gamma_{c+1} (G)$ is free abelian, so, by Theorem \ref{202509021146}, $\mathfrak{Q}_d (G)$ is also free-abelian.

Let $A$ be a non-negatively graded free abelian group whose homogeneous component $A_i$ is of finite rank.
Clearly, the abelian group $\mathds{S}_{d}(A)$, introduced in (\ref{202510091635}), is free and has finite rank.
Using this, we give an elementary observation as follows:
\begin{Lemma} \label{202510091657}
    Let $A,B$ be non-negatively graded free abelian groups whose homogeneous components are of finite rank.
    Let $f : A \to B$ be a surjective homomorphism of graded abelian groups.
    For $c \in \mathds{N}$, we have the following:
    \begin{enumerate}
        \item If $f$ is an isomorphism in degrees $\leq c$, then $\mathrm{rank} (\mathds{S}_d (A)) = \mathrm{rank}(\mathds{S}_d (B))$ for $d \leq c$.
        \item If $A_1 \not\cong 0$, and $f$ is not injective at the degree $c+1$, then $\mathrm{rank} (\mathds{S}_d  (A)) > \mathrm{rank} (\mathds{S}_d  (B))$ for $d > c$. 
    \end{enumerate}
\end{Lemma}
\begin{proof}
    Since $A,B$ are non-negatively graded,
    the direct sum in (\ref{202510091635}) for $A,B$ is taken over all non-negative integers $a_1, \cdots, a_n$ such that $\sum^{d}_{i=1} ia_i =d$.
    Hence, by the assumption $d \leq c$, the induced morphism $\mathds{S}_d (A) \to \mathds{S}_d (B)$ is an isomorphism for $d \leq c$.
    It proves (1).

    We now prove (2).
    Let $K$ be the kernel of $f : A_{c+1} \to B_{c+1}$.
    This is not trivial by the assumption since $f$ is not injective at degree $c +1$.
    Since $A_1$ is free and non-trivial, so is $S^{d-(c+1)} ( A_{1})$, hence $S^{d-(c+1)} ( A_{1}) \otimes_{\mathds{Z}} K \not\cong 0$.
    On the one hand, $S^{d-(c+1)} ( A_{1}) \otimes_{\mathds{Z}} K \subset S^{d-(c+1)} ( A_{1}) \otimes_{\mathds{Z}} A_{ c+1} \subset \mathds{S}_{d} (A)$ vanishes under the application of $f$.
    Hence, $\mathds{S}_d  (A) \to \mathds{S}_d  (B)$ is not injective for $d > c$.
    This implies the strict inequality of the ranks in the statement, since $f$ is surjective.
\end{proof}

\begin{notation}
    Let $\mathsf{F}_{n}^{(c)} {:=} \mathsf{F}_n / \gamma_{c+1} (\mathsf{F}_n)$ for $c \in \mathds{N}$.
    We also set $\mathsf{F}_{n}^{(\infty)} {:=} \mathsf{F}_n$.
\end{notation}

Let $c_0,c_1 \in \mathds{N}^\ast \cup \{ \infty \}$ such that $c_0 < c_1$.
The canonical projection $\mathsf{F}_{n}^{(c_1)} \to \mathsf{F}_{n}^{(c_0)}$ gives a surjection $\mathfrak{Q}_d (\mathsf{F}_{n}^{(c_1)}) \to \mathfrak{Q}_d (\mathsf{F}_{n}^{(c_0)})$. 
In what follows, we give a necessary and sufficient condition for this to be an isomorphism:
\begin{Lemma} \label{202509021227}
    Let $d \in \mathds{N}$ and $n,m \in \mathds{N}^\ast$.
    \begin{enumerate}
        \item If $d=0$ or $n=1$, then the map $\mathfrak{Q}_d ((\mathsf{F}_{n}^{(c_1)})^{\times m}) \to \mathfrak{Q}_d ((\mathsf{F}_{n}^{(c_0)})^{\times m})$ is an isomorphism.
        \item If $d \geq 1$ and $n \geq 2$, then $\mathfrak{Q}_d ((\mathsf{F}_{n}^{(c_1)})^{\times m}) \to \mathfrak{Q}_d ((\mathsf{F}_{n}^{(c_0)})^{\times m})$ is an isomorphism if and only if $d \leq c_0$.
    \end{enumerate}
\end{Lemma}
\begin{proof}
    The first assertion is immediate from the fact that $\mathfrak{Q}_0 (G) \cong \mathds{Z}$ for a group $G$; and $\mathsf{F}_{1}^{(c)} \cong \mathsf{F}_1$.

    We now assume that $d \geq 1$ and $n \geq 2$, and prove the second part.
    Since $\mathfrak{Q}_d ((\mathsf{F}_{n}^{(c_1)})^{\times m}) \to \mathfrak{Q}_d ((\mathsf{F}_{n}^{(c_0)})^{\times m})$ is surjective, it suffices to compare the ranks of $\mathfrak{Q}_d ((\mathsf{F}_{n}^{(c_1)})^{\times m})$ and $\mathfrak{Q}_d ((\mathsf{F}_{n}^{(c_0)})^{\times m})$.
    By the Sandling--Tahara theorem, their ranks are computed from the homogeneous component of the symmetric algebra of the associated Lie rings.
    We now apply Lemma \ref{202510091657} to $A = \mathfrak{L} ( (\mathsf{F}_{n}^{(c_1)})^{\times m} ; \gamma_{\bullet})$ and $B = \mathfrak{L} ( (\mathsf{F}_{n}^{(c_0)})^{\times m} ; \gamma_{\bullet})$ with the map $f$ induced by the canonical surjection $\mathsf{F}_{n}^{(c_1)} \to \mathsf{F}_{n}^{(c_0)}$.
    By the isomorphism in (\ref{202509011229}), the map $\mathfrak{L} ( (\mathsf{F}_{n}^{(c_1)})^{\times m} ; \gamma_{\bullet}) \to \mathfrak{L} ( (\mathsf{F}_{n}^{(c_0)})^{\times m} ; \gamma_{\bullet})$ is an isomorphism in degrees $\leq c_0$.
    Furthermore, it is not injective at degree $c_0+1$, since, for $i \in \mathds{N}^\ast$, the quotient $\gamma_{i} (\mathsf{F}_n)/ \gamma_{i+1} (\mathsf{F}_n)$ is nontrivial by the assumption $n \geq 2$.
    It is obvious that the degree $1$ component of $\mathfrak{L} ( (\mathsf{F}_{n}^{(c_1)})^{\times m} ; \gamma_{\bullet})$, i.e. $(\mathsf{F}_{n}/\gamma_2(\mathsf{F}_{n}))^{\times m} \cong \mathds{Z}^{\times nm}$, is nontrivial and free.
\end{proof}

As stated in Theorem \ref{202508250938}, the monad epimorphism $\mathtt{L}_{\mathbf{N}_{c_1}^{\mathsf{o}}} \to \mathtt{L}_{\mathbf{N}_{c_0}^{\mathsf{o}}}$ maps $\mathtt{I}^{(d)}_{\mathbf{N}_{c_1}^\mathsf{o}}$ onto $\mathtt{I}^{(d)}_{\mathbf{N}_{c_0}^\mathsf{o}}$.
In particular, it induces an epimorphism $\mathtt{I}^{(d-1)}_{\mathbf{N}_{c_1}^\mathsf{o}}/ \mathtt{I}^{(d)}_{\mathbf{N}_{c_1}^\mathsf{o}} \to \mathtt{I}^{(d-1)}_{\mathbf{N}_{c_0}^\mathsf{o}}/ \mathtt{I}^{(d)}_{\mathbf{N}_{c_0}^\mathsf{o}}$.

\begin{Lemma} \label{202509021409}
    Let $d \in \mathds{N}$.
    We have $d \leq c_0$ if and only if the induced homomorphism $\mathtt{I}^{(d-1)}_{\mathbf{N}_{c_1}^\mathsf{o}}/ \mathtt{I}^{(d)}_{\mathbf{N}_{c_1}^\mathsf{o}} \to \mathtt{I}^{(d-1)}_{\mathbf{N}_{c_0}^\mathsf{o}}/ \mathtt{I}^{(d)}_{\mathbf{N}_{c_0}^\mathsf{o}}$ is an isomorphism.
\end{Lemma}
\begin{proof}
    By Theorem \ref{202408011134}, it suffices to prove that $0 \leq d \leq c_0$ if and only if, for $n,m \in \mathds{N}$, the map $$\Aug((\mathsf{F}_{n}^{(c_1)})^{\times m})^{d} / \Aug((\mathsf{F}_{n}^{(c_1)})^{\times m})^{d+1} \to \Aug((\mathsf{F}_{n}^{(c_0)})^{\times m})^{d} / \Aug((\mathsf{F}_{n}^{(c_0)})^{\times m})^{d+1}$$ is an isomorphism.
    The proof for the ground ring $\mathds{k}=\mathds{Z}$ follows from Lemma \ref{202509021227}.
    Using Lemma \ref{202508291109}, it extends to general $\mathds{k}$.
    In fact, by Theorem \ref{202509021146}, $\Aug_{\mathds{Z}}(G)^{d}/\Aug_{\mathds{Z}}(G)^{d+1} = \mathfrak{Q}_d(G)$ is free-abelian for a product of free nilpotent groups $G$.
\end{proof}

\begin{proof}[Proof of Theorem \ref{202509021231}]
    The theorem is proved by combining Lemma \ref{202509021409} and Theorem \ref{202509291857}.
    In particular, if $d > c_0$, then the map $\mathcal{S}^{d}_{\mathbf{N}_{c_0}^{\mathsf{o}}} \to \mathcal{S}^{d}_{\mathbf{N}_{c_1}^{\mathsf{o}}}$ is properly injective, so the functor $\mathtt{L}_{\mathbf{N}_{c_0}^{\mathsf{o}}}\mbox{-}\mathsf{Mod}^{\leq d} \to \mathtt{L}_{\mathbf{N}_{c_1}^{\mathsf{o}}}\mbox{-}\mathsf{Mod}^{\leq d}$ is not essentially surjective.
\end{proof}

\subsection{Comparison with $\Rac^{\mathsf{o}}$-modules}
\label{202510261507}
In this section, we give a comparison between polynomial $\Rac^{\mathsf{o}}$-modules and polynomial $\mathbf{N}_{c}^{\mathsf{o}}$-modules, as a generalization of Theorem \ref{202509021231}, where $\mathcal{R}$ is a radical functor for groups.
To this end, we begin with a reinterpretation of Theorem \ref{202509021231}, which follows from Theorem \ref{202509291857}:
\begin{Corollary}
    For $\{c_0 < c_1\} \subset \mathds{N}^\ast \cup \{ \infty\}$ and $d \in \mathds{N}$, we have
    $$
    \mathtt{K}_{\mathbf{N}_{c_1}^{\mathsf{o}}\to \mathbf{N}_{c_0}^{\mathsf{o}}} \subset \mathtt{I}^{(d)}_{\mathbf{N}_{c_1}^{\mathsf{o}}} \quad \mathrm{if~and~only~if} \quad d \leq c_0.
    $$
\end{Corollary}
This consequence can be viewed as an interpretation of Theorem \ref{202509021231} in terms of the intrinsic information of the monad $\mathtt{L}_{\mathbf{N}_{c_1}^{\mathsf{o}}}$.
In what follows, we give a partial generalization of this:

\begin{Lemma} \label{202510261204}
    Let $c_0 \in \mathds{N}^\ast$ and $\mathcal{R}$ be a radical functor for groups such that $\mathcal{R} \subset \gamma_{c_0+1}$.
    Let $\pi : \G^{\mathsf{o}} \to \Rac^{\mathsf{o}}$ be the full functor.
    For $d\in [c_0]$, we have $\mathtt{K}_{\pi}\subset \mathtt{I}^{(d)}_{\G^{\mathsf{o}}}$.
\end{Lemma}
\begin{proof}
    We prove $\mathtt{K}_{\pi} (m,n) \subset \mathtt{I}^{(c_0)}_{\G^{\mathsf{o}}} (m,n)$ for $n,m\in \mathds{N}$.
    For $G = \mathsf{F}_n^{\times m}$, by Theorem \ref{202408011134}, we have $\mathtt{I}^{(c_0)}_{\G^{\mathsf{o}}} (m,n) \cong \Aug ( G)^{c_0+1}$.
    On the other hand, by definition, $\mathtt{K}_{\pi} (m,n)$ can be viewed as the kernel of the monoid algebra map $\mathds{k}[G] \to \mathds{k}[G/\mathcal{R}(G)]$.
    Hence, it is isomorphic to the two-sided ideal of the monoid algebra $\mathds{k}[G]$ generated by $(u-1)$ for $u \in \mathcal{R}(G)$.
    Since $\mathcal{R} \subset \gamma_{c_0+1}$, we see that $(u-1) \in \Aug (G)^{c_0+1}$, so that we obtain $\mathtt{K}_{\pi} (m,n) \subset \mathtt{I}^{(c_0)}_{\G^{\mathsf{o}}} (m,n)$.
\end{proof}

\begin{theorem}  \label{202510281425}
    Let $c_0 \in \mathds{N}^\ast$ and $\mathcal{R}$ be a radical functor for groups such that $\mathcal{R} \subset \gamma_{c_0+1}$.
    Then we have
    $$
    [c_0] \subset \Gamma ( \G^{\mathsf{o}}  \to \Rac^{\mathsf{o}}) \cap \Gamma (\Rac^{\mathsf{o}} \to \mathbf{N}_{c_0}^{\mathsf{o}}).
    $$
\end{theorem}
\begin{proof}
    It suffices to prove that $[c_0] \subset \Gamma ( \G^{\mathsf{o}}  \to \Rac^{\mathsf{o}})$.
    If so, then, by the 2-out-of-3 property of bijections and Theorem \ref{202509021231}, we see that $[c_0] \subset \Gamma (\Rac^{\mathsf{o}} \to \mathbf{N}_{c_0}^{\mathsf{o}})$.

    The proof of $[c_0] \subset \Gamma ( \G^{\mathsf{o}}  \to \Rac^{\mathsf{o}})$ also reduces to that of $c_0 \in \Gamma ( \G^{\mathsf{o}}  \to \Rac^{\mathsf{o}})$ by the last statement of Theorem \ref{202509291857}.
    By the theorem, the condition that $c_0 \in \Gamma ( \G^{\mathsf{o}}  \to \Rac^{\mathsf{o}})$ is equivalent to $\mathtt{K}_{\pi} \subset \mathtt{I}^{(c_0)}_{\G^{\mathsf{o}}}$.
    Thus, the desired statement follows from Lemma \ref{202510261204}.
\end{proof}

\begin{Example} \label{202510221956}
    Let $c \geq c_0$.
    One may consider the $c$-th derived subgroup $\gamma_2^{c}$ as $\mathcal{R}$ in the theorem.
    The theorem yields an equivalence of polynomial functors, of $\deg \leq c_0$, over the category of free solvable groups of derived length $\leq c$ and that of free nilpotent groups of class $\leq c_0$.
    This example is also generalized as follows.
\end{Example}

\begin{Example} \label{202510151133}
    Let $\mathcal{R}_i, ~i\in [c_0]$ be surjection-preserving radical functors for groups in the sense of Definition \ref{202510091828}. 
    By Proposition \ref{202510111100}, the recursive commutator $\mathcal{R} = [\cdots [[\mathcal{R}_0,\mathcal{R}_1], \mathcal{R}_2], \cdots, \mathcal{R}_{c_0}]$ gives a radical functor such that $\mathcal{R} \subset \gamma_{c_0+1}$, to which the theorem applies.
\end{Example}

\section{Polynomial $\Dim_{c,p}^{\mathsf{o}}$-modules}
\label{202510211825}

We assume, only in this section, that the ground ring $\mathds{k}$ is a field with characteristic $p \geq 0$.
This section aims at proving some variations of Theorem \ref{202509021231} by considering the category $\Dim_{c,p}^{\mathsf{o}}$ in place of $\mathbf{N}_{c}^{\mathsf{o}}$.
We recall $\Dim_{c,p}$ from Notation \ref{202512281222}.
As a consequence, we obtain a modular version of the theorem which leads to the main results of this paper, namely Theorems \ref{202510191826} and \ref{202510161105}.

We begin with the main theorem of this section in full generality.

\begin{Defn}
    For a group $G$, we define the subset $\Psi_{\mathcal{R},p} (G) \subset \mathds{N}$ by
    \begin{align*}
        \Psi_{\mathcal{R},p} (G) {:=}
        \begin{cases}
            \{ c \in \mathds{N} \mid D_{c+1,p} (G) (D_{c,p}(G ) \cap \mathcal{R}( G ) )= D_{c,p}(G) \} & p>0  , \\
            \{ c \in \mathds{N} \mid \sqrt{\gamma_{c+1}(G)\mathcal{R}( G)} = \sqrt{\gamma_{c}(G)\mathcal{R}( G )} \} & p = 0 .
        \end{cases}
    \end{align*}
\end{Defn}

Recall the notation $(-)^\spadesuit$ from Notation \ref{202510221957}.

\begin{theorem} \label{202510021437}
    Suppose that the ground ring $\mathds{k}$ is a field with characteristic $p \geq 0$.
    Let $c_0 \in \mathds{N}^\ast$ and $\mathcal{R}$ be a radical functor for groups such that $\mathcal{R} \subset D_{c_0+1,p}$.
    Then we have 
    \begin{align*}
        \Gamma ( \Rac^{\mathsf{o}} \to \Dim_{c_0,p}^{\mathsf{o}})  = \left( [c_0] ~\cup~ \bigcap_{n \in \mathds{N}} \Psi_{\mathcal{R},p} (\mathsf{F}_n) \right)^\spadesuit .
    \end{align*}
    In particular, $[c_0] \subset \Gamma ( \Rac^{\mathsf{o}} \to \Dim_{c_0,p}^{\mathsf{o}})$.
\end{theorem}

We defer the proof to Sections \ref{202510151145} and \ref{202512261437}.
The proof is parallel to that of Theorem \ref{202509021231} except that we need Quillen's theorem instead of the Sandling--Tahara theorem.
Since $\mathcal{R}$ is fixed below, we abbreviate $\Psi_{\mathcal{R},p} (G) = \Psi_{p} (G)$ unless otherwise specified.
    
\begin{remark}
    Note that, in general, the set $\Gamma ( \Rac^{\mathsf{o}} \to \Dim_{c_0,p}^{\mathsf{o}})$ properly contains the set $[c_0]$, as demonstrated in the following example.
\end{remark}

\begin{Example} 
    We assume that $p > 0$.
    Let $\mathcal{R}(G) = \prod_{ip^j = c_0+1} \gamma_i (G)^{p^j}$.
    It yields a radical functor for groups by Proposition \ref{202510161130}.
    By the description of $D_{c,p}$ given in Example \ref{202510121414}, $\mathcal{R} \subset D_{c_0+1,p}$, so we can apply Theorem \ref{202510021437} to this $\mathcal{R}$.
    We then have $$[c_0+1] \subset \Gamma ( \Rac^{\mathsf{o}} \to \Dim_{c_0,p}^{\mathsf{o}}).$$
    In fact, we have $c_0+1 \in \Psi_p (F)$ for any free group $F$, since 
    \begin{align*}
        D_{c_0+2,p} (F) ( D_{c_0+1,p}(F) \cap \mathcal{R}(F)) = D_{c_0+2,p} (F) \mathcal{R}(F) = \prod_{sp^t \geq c_0+2} \gamma_s (F)^{p^t} \prod_{ip^j = c_0+1} \gamma_i (F)^{p^j} = D_{c_0+1,p}(F) .
    \end{align*}
\end{Example}

\subsection{Proofs of Theorems \ref{202510191826} and \ref{202510161105}}
\label{202512301553}

In this section, as corollaries of Theorem \ref{202510021437}, we prove Theorems \ref{202510191826} and \ref{202510161105}.

\begin{proof}[Proof of Theorem \ref{202510191826}]
    It suffices to investigate the set $\Gamma( \mathbf{N}_{c_0}^{\mathsf{o}} \to \Dim_{c_0,p}^{\mathsf{o}} )$, since the non-opposite part follows from Corollary \ref{202510141158}.
    For the proof, we apply Theorem \ref{202510021437} to $p >0$ and $\mathcal{R} = \gamma_{c_0+1}$.
    Then we obtain the left inclusion of the statement, i.e. $[c_0] \subset \Gamma( \mathbf{N}_{c_0}^{\mathsf{o}} \to \Dim_{c_0,p}^{\mathsf{o}} )$.
    To establish the right inclusion, we examine the theorem in more detail.
    We set
    $$S = [c_0] \cup \{ c \in \mathds{N} \mid D_{c+1,p}(\mathsf{F}_1)  = D_{c,p} ( \mathsf{F}_1 ) \}.$$
    By definition, we have $S = [c_0] \cup \Psi_{p} (\mathsf{F}_1)$, since $\mathcal{R} (\mathsf{F}_1) = \gamma_{c_0+1} (\mathsf{F}_1) = 1$.
    Thus, Theorem \ref{202510021437} implies that $\Gamma( \mathbf{N}_{c_0}^{\mathsf{o}} \to \Dim_{c_0,p}^{\mathsf{o}} ) \subset S^{\spadesuit}$.
    On the one hand, by using the fact recalled in Example \ref{202510121414}, we have
    $$
    D_{p^{r_0}+1,p}(\mathsf{F}_1) = \mathsf{F}_1^{p^{r_0+1}} \neq \mathsf{F}_1^{p^{r_0}} = D_{p^{r_0},p} ( \mathsf{F}_1 ) .
    $$
    Hence, $p^{r_0} \not\in S$, so $S^\spadesuit \subset [ p^{r_0}-1]$.
\end{proof}

\begin{proof}[Proof of Theorem \ref{202510161105}]
    It suffices to investigate the set $\Gamma (\Dim_{c_1,p}^{\mathsf{o}} \to \Dim_{c_0,p}^{\mathsf{o}})$, since the non-opposite part follows from Corollary \ref{202510141158}.
    Let $c_1 \in \mathds{N}$.
    In order to prove the theorem, we apply Theorem \ref{202510021437} to $p>0$ and $\mathcal{R} = D_{c_1+1,p}$.
    For any free group $F$ of rank $\geq 2$, by Proposition \ref{202510231737}, we have $$D_{c_0+2,p} (F) ( D_{c_0+1,p}(F) \cap D_{c_1+1,p} (F) ) \subset  D_{c_0+2,p}(F)  \subsetneq D_{c_0+1,p}(F),$$
    which implies $c_0 +1 \not\in \Psi_p (F)$.
    Thus, by Theorem \ref{202510021437}, we obtain $c_0+1 \not\in \Gamma ( \Rac^{\mathsf{o}} \to \Dim_{c_0,p}^{\mathsf{o}})$, so $\Gamma ( \Rac^{\mathsf{o}} \to \Dim_{c_0,p}^{\mathsf{o}}) = [c_0]$.
    Analogously, the case where $c_1 = \infty$ is deduced from the application of Theorem \ref{202510021437} to $\mathcal{R} \equiv 1$.
\end{proof}

\begin{remark}
    One may prove the analogous statement for $p = 0$, while it is, in fact, implied by Theorem \ref{202509021231} under the equality in Remark \ref{202510221959}.
    We note that the latter holds for a general ground ring $\mathds{k}$.
\end{remark}

\subsection{Proof of Theorem \ref{202510021437} for $p>0$}
\label{202510151145}

In this section, we prove Theorem \ref{202510021437} for positive characteristic.
We begin with an elementary statement that holds for any characteristic $p \geq 0$:
\begin{Lemma} \label{202510031338}
    Let $\mathfrak{g},\mathfrak{h}$ be non-negatively graded $p$-restricted Lie $\mathds{k}$-algebras.
    Let $f : \mathfrak{g} \to \mathfrak{h}$ be a homomorphism of graded $p$-restricted Lie algebras.
    For $d\in\mathds{N}$, if $f$ is an isomorphism in $\deg \leq d$, then the induced map $\mathcal{U}^{p}(\mathfrak{g}) \to \mathcal{U}^{p} (\mathfrak{h})$ is an isomorphism in $\deg \leq d$.
    Furthermore, if $f$ is not injective at $\deg = d+1$, then $f_\ast : \mathcal{U}^{p}(\mathfrak{g}) \to \mathcal{U}^{p}(\mathfrak{h})$ is not injective at $\deg = d+1$.
\end{Lemma}
\begin{proof}
    The first assertion follows from the definitions that, for $c \leq d$, $\mathcal{U}^{p}_{c}(\mathfrak{g})$ is the quotient of the degree $c$-part of the tensor algebra of $\mathfrak{g}$ by the homogeneous component generated by the relations consisting of elements of degree $\leq c \leq d$.
    
    We now prove the remaining assertion.
    By the assumptions, we can take $0 \neq y \in \mathfrak{g}_{d+1}$ such that $f(y) = 0$.
    We choose a well-ordered basis of $\mathfrak{g}$ containing $y$.
    As a consequence of the Poincar\'e-Birkhoff-Witt theorem for $p$-restricted Lie algebras (in both the cases $p=0$ or $p> 0$), $0 \neq y \in \mathcal{U}^{p}_{d+1}(\mathfrak{g})$, but $f_\ast (y) = 0$.
\end{proof}

We now consider the setting of Theorem \ref{202510021437}.
Let $\mathcal{R}$ be a radical functor for groups such that $\mathcal{R} \subset D_{c_0+1,p}$.
Let $G$ be a group.
Let
$$
G_1 = G/ D_{c_0+1,p} (G), \quad G_2 = G / \mathcal{R} (G).
$$
We describe $\mathfrak{L}(G_1 ; D_{\bullet,p})$ and $\mathfrak{L}(G_2 ; D_{\bullet,p})$ in terms of $G$.
To do that, we use some general properties of radical functors given in Section \ref{202510071642}.
By Lemma \ref{202510031417} applied to $\mathcal{R}_1 = D_{c+1,p}$, $\mathcal{R}_2 = D_{c,p}$ and $H = D_{c_0+1,p}(G)$, we obtain the following:
\begin{align*}
\mathfrak{L}(G_1 ; D_{\bullet,p}) = \bigoplus_{1 \leq c} D_{c,p}(G_1) / D_{c+1,p}(G_1) \cong \bigoplus_{1 \leq c \leq c_0} D_{c,p}(G) / D_{c+1,p}(G) .
\end{align*}
Similarly, we obtain
\begin{align} \label{202510151214}
    \mathfrak{L}(G_2 ; D_{\bullet,p}) \cong   \bigoplus_{1 \leq c}  D_{c,p}(G) /\left( D_{c+1,p}(G) (D_{c,p}(G) \cap \mathcal{R}(G)) \right) .
\end{align}
In particular, by the hypothesis $\mathcal{R} \subset D_{c_0+1,p}$, the components of  $\mathfrak{L}(G_2 ; D_{\bullet,p})$ in degrees $\leq c_0$ coincide with $\mathfrak{L}(G_1 ; D_{\bullet,p})$.

\begin{Lemma} \label{202510031436}
    For $c \leq c_0$, the induced map $\Aug ( G_2)^{c} / \Aug ( G_2)^{c+1} \to \Aug ( G_1)^{c} / \Aug ( G_1)^{c+1}$ is an isomorphism.
\end{Lemma}
\begin{proof}
    By applying Lemma \ref{202510031338} to the previous observation, $\mathcal{U}^{p}_{c}(\mathds{k} \otimes_{\mathds{Z}} \mathfrak{L}(G_2; D_{\bullet,p}) ) \to \mathcal{U}^{p}_{c} (\mathds{k} \otimes_{\mathds{Z}} \mathfrak{L}(G_1; D_{\bullet,p}) )$ is an isomorphism for $c \leq c_0$.
    By Theorem \ref{202510091518}, we obtain the desired statement.
\end{proof}

\begin{Lemma} \label{202510091734}
    $
    \{ c \in \mathds{N} \mid \Aug ( G_2)^{c} / \Aug ( G_2)^{c+1} \stackrel{\cong}{\longrightarrow} \Aug ( G_1)^{c} / \Aug ( G_1)^{c+1} \}^\spadesuit = \left( [c_0] \cup \Psi_p (G) \right)^\spadesuit 
    $.
\end{Lemma}
\begin{proof}
    By Lemma \ref{202510031436}, the set $[c_0]$ is contained in the left hand side.
    Let $m \in \mathds{N} \cup \{\infty\}$ be the supremum of the right hand side.
    Then, for $c \in \mathds{N}$ such that $c_0 < c \leq m$, we have $D_{c+1,p} (G) (D_{c,p} (G) \cap \mathcal{R}(G)) = D_{c,p}(G)$.
    Hence, the $c$-th component of $\mathfrak{L}(G_2; D_{\bullet,p})$ is trivial by (\ref{202510151214}), so that the surjective map $\mathfrak{L}(G_2; D_{\bullet,p}) \to \mathfrak{L}(G_1; D_{\bullet,p})$ gives an isomorphism in $\deg \leq m$.
    Moreover, the maximality of $m$ implies that it is not an isomorphism at $\deg = m+1$.
    Therefore, by Lemma \ref{202510031338}, the map $\mathcal{U}^{p} ( \mathds{k} \otimes_{\mathds{Z}} \mathfrak{L}(G_2; D_{\bullet,p}) ) \to \mathcal{U}^{p}(\mathds{k} \otimes_{\mathds{Z}}  \mathfrak{L}(G_1; D_{\bullet,p}))$ gives an isomorphism in $\deg \leq m$; and a non-isomorphism at $\deg = m+1$.
    By Theorem \ref{202510091518}, we obtain the statement.
\end{proof}

\begin{proof}[Proof of Theorem \ref{202510021437} for $p>0$]
    By Theorem \ref{202509291857}, it suffices to investigate the set $\mathrm{B} ( \Rac^{\mathsf{o}} \to \Dim_{c_0,p}^{\mathsf{o}} )^{\spadesuit}$.
    By the natural isomorphism in Theorem \ref{202408011134}, this is equal to the following set:
    \begin{align} \label{202510131354}
        \left\{ c \in \mathds{N} \mid \forall n \forall m\in \mathds{N}~\left( \Aug( \mathcal{C}_n^{\times m})^{c} /\Aug( \mathcal{C}_n^{\times m})^{c+1} \stackrel{\cong}{\longrightarrow} \Aug( \mathcal{D}_n^{\times m})^{c} /\Aug( \mathcal{D}_n^{\times m})^{c+1} \right) \right\}^\spadesuit 
    \end{align}
    where $\mathcal{C}_n = \mathsf{F}_n / \mathcal{R}(\mathsf{F}_n)$ and $\mathcal{D}_n = \mathsf{F}_n / D_{c_0+1,p}(\mathsf{F}_n)$.
    By applying Lemma \ref{202510091734} to $G= \mathsf{F}_n^{\times m}$ (especially, $G_1 = \mathcal{D}_n^{\times m}, G_2 = \mathcal{C}_n^{\times m}$) for each $n,m$, the set in (\ref{202510131354}) coincides with
    $$
    \left( [c_0] \cup \bigcap_{n,m\in\mathds{N}} \Psi_p ( \mathsf{F}_n^{\times m}) \right)^\spadesuit .
    $$
    Since a radical functor for groups preserves products of groups, as in Lemma \ref{202510131400}, we have $\Psi_{p} (\mathsf{F}_n^{\times m}) = \Psi_{p} (\mathsf{F}_n)$.
    Thus, this set coincides with
    the one in the statement.
\end{proof}

\subsection{Proof of Theorem \ref{202510021437} for $p=0$} \label{202512261437}

In this section, using a method parallel to that of the previous section, we prove Theorem \ref{202510021437} for characteristic zero.
We assume that $G$ is a product of free groups, and continue to use the notation $G,G_1,G_2$ as in Section \ref{202510151145}.
By Lemma \ref{202510091525}, 
$$
D_{c,0}(G_1)/D_{c+1,0}(G_1) \cong \sqrt{\gamma_c(G)D_{c_0+1,0} (G)} / \sqrt{\gamma_{c+1}(G)D_{c_0+1,0} (G)} , \quad c \in \mathds{N}^\ast.
$$
We have $D_{c,0} (G) = \sqrt{\gamma_c(G)} = \gamma_c(G)$, since $G$ is a product of free groups.
Hence, $D_{c,0}(G_1)/D_{c+1,0}(G_1)$ is isomorphic to $\sqrt{\gamma_c(G)} / \sqrt{\gamma_{c+1}(G)} = D_{c,0}(G)/D_{c+1,0}(G)$ for $c \leq c_0$; and trivial otherwise.
Hence, $\mathfrak{L}(G_1;D_{\bullet,0})$ coincides with the $\deg \leq c_0$ subspace of $\mathfrak{L}(G;D_{\bullet,0})$.
Likewise, we obtain
$$
D_{c,0}(G_2)/D_{c+1,0}(G_2) \cong \sqrt{\gamma_c(G) \mathcal{R} (G)} / \sqrt{\gamma_{c+1}(G) \mathcal{R} (G)} 
$$
which is isomorphic to $D_{c,0}(G)/D_{c+1,0}(G)$ for $c \leq c_0$ by the hypothesis $\mathcal{R} \subset D_{c_0+1,0}$.

\begin{Lemma} \label{202510131343}
    If $G$ is a product of free groups, then the induced map $\Aug ( G_2)^{c} / \Aug ( G_2)^{c+1} \to \Aug ( G_1)^{c} / \Aug ( G_1)^{c+1}$ is an isomorphism for $c \leq c_0$.
\end{Lemma}
\begin{proof}
    By applying Lemma \ref{202510031338} to the above observation, $\mathcal{U}^{0}_{c}(\mathds{k} \otimes_{\mathds{Z}} \mathfrak{L}(G_2; D_{\bullet,0}) ) \to \mathcal{U}^{0}_{c} (\mathds{k} \otimes_{\mathds{Z}} \mathfrak{L}(G_1; D_{\bullet,0}) )$ is an isomorphism for $c \leq c_0$.
    Hence, Theorem \ref{202510091518} leads to the statement.
\end{proof}
    
\begin{Lemma} \label{202510091849}
    If $G$ is a product of free groups, then we have
    $$
    \{ c \in \mathds{N} \mid \Aug ( G_2)^{c} / \Aug ( G_2)^{c+1} \stackrel{\cong}{\to} \Aug ( G_1)^{c} / \Aug ( G_1)^{c+1} \}^\spadesuit = \left( [c_0] \cup \Psi_0 (G) \right)^\spadesuit .
    $$
\end{Lemma}
\begin{proof}
    By Lemma \ref{202510131343}, the set $[c_0]$ is contained in the left hand side.
    The rest of the proof is parallel with that of Lemma \ref{202510091734}.
\end{proof}

\begin{proof}[Proof of Theorem \ref{202510021437} for $p=0$]
    The proof is analogous to that for $p>0$, as given in Section \ref{202510151145}.
    In particular, it can be proved by using Lemma \ref{202510091849}.
\end{proof}

\section{Non-existence of a core internalizer for analyticity}
\label{202509301437}

Let $\mathcal{C}$ be a Lawvere theory and $\mathcal{S}^{\omega}_{\mathcal{C}}$ denote the subclass of $\mathtt{L}_{\mathcal{C}}\mbox{-}\mathcal{M}\mathsf{od}$ consisting of isomorphism classes of analytic $\mathtt{L}_{\mathcal{C}}$-modules.
In this section, we investigate the existence of a core $\mathtt{L}_{\mathcal{C}}$-internalizer for the subclass $\mathcal{S}^{\omega}_{\mathcal{C}}$ (see Definition \ref{202512281236}).
We demonstrate that, for many Lawvere theories $\mathcal{C}$ of our interest, such a core $\mathtt{L}_{\mathcal{C}}$-internalizer does not exist, apart from trivial cases.

The following is the main theorem of this section:
\begin{theorem} \label{202508271151}
    Let $\mathcal{C}$ be a Lawvere theory satisfying condition (ZM*).
    The following are equivalent to each other:
    \begin{enumerate}
        \item There exists a core $\mathtt{L}_{\mathcal{C}}$-internalizer for $\mathcal{S}^{\omega}_{\mathcal{C}}$.
        \item For $n \in \mathds{N}$, the filtration $\Aug ( \mathcal{C}_{n}^{\times m})^{\bullet}$ consisting of powers of the augmentation ideal stabilizes uniformly with respect to $m \in \mathds{N}$.
        In other words, there exists $d_0 = d_0(n) \in \mathds{N}$ such that $\Aug ( \mathcal{C}_{n}^{\times m})^{d+1} = \Aug ( \mathcal{C}_{n}^{\times m})^{d+2}$ for $d \geq d_0$ and $m \in \mathds{N}$.
    \end{enumerate}
\end{theorem}

The proof is postponed to the end of this section.
Below we present some examples of $\mathcal{C}$ that admit no core $\mathtt{L}_{\mathcal{C}}$-internalizer for $\mathcal{S}^{\omega}_{\mathcal{C}}$.

\begin{remark} \label{202512251531}
    Before this, we observe that there exists an obvious core $\mathtt{L}_{\mathcal{C}}$-internalizer for $\mathcal{S}^{\omega}_{\mathcal{C}}$ when $\mathrm{D}(\mathcal{C}) = \{ 0\}$, i.e. all analytic functors are constant.
    In fact, it is given by $\mathtt{I}^{(0)}_{\mathcal{C}}$.
    Such examples arise from Examples \ref{202509241443} and \ref{202509301321}.
\end{remark}

Using the contraposition of the theorem, we derive a series of corollaries with examples.

\begin{Corollary}
    We assume that $\mathds{k} \neq 0$.
    If there exists $n \in \mathds{N}$ such that $\mathfrak{Q}_d (\mathcal{C}_n), ~d\in\mathds{N}$ are free abelian groups whose ranks do not converge to $0$ as $d \to \infty$, then the class $\mathcal{S}^{\omega}_{\mathcal{C}}$ admits no core $\mathtt{L}_{\mathcal{C}}$-internalizer.
\end{Corollary}
\begin{proof}
    Suppose that the class $\mathcal{S}^{\omega}_{\mathcal{C}}$ admits a core $\mathtt{L}_{\mathcal{C}}$-internalizer.
    By Theorem \ref{202508271151}, there exists $d_0 \in \mathds{N}$ such that, for $d \geq d_0$, we have $\Aug ( \mathcal{C}_n)^{d+1} = \Aug ( \mathcal{C}_n)^{d+2}$.
    By Lemma \ref{202508291109}, the quotient $\Aug (  \mathcal{C}_n)^{d+1} /\Aug ( \mathcal{C}_n)^{d+2}$ is a free $\mathds{k}$-module whose rank equals the rank of $\mathfrak{Q}_{d+1} ( \mathcal{C}_n)$.
    Hence, for $d \geq d_0$, the rank of $\mathfrak{Q}_{d+1} ( \mathcal{C}_n)$ is zero.
    This contradicts the assumption.
\end{proof}

\begin{Example}
    The category $\G^{\mathsf{o}}$, as well as $\mathcal{C} \in \{ \mathbf{N}_{c}^{\mathsf{o}} \mid c \in \mathds{N}^\ast \}$, satisfy the previous assumption, as shown in Section \ref{202510211823}.
    Hence, if $\mathds{k} \neq 0$, then the class $\mathcal{S}^{\omega}_{\mathcal{C}}$ admits no core $\mathtt{L}_{\mathcal{C}}$-internalizer.
\end{Example}

Analogous arguments hold for the Lawvere theory $\Dim_{c,p}^{\mathsf{o}}$:
\begin{Corollary}
    We assume that the ground ring $\mathds{k}$ is a field of positive characteristic $p$.
    For $\mathcal{C}= \Dim_{c,p}^{\mathsf{o}}$, the class $\mathcal{S}^{\omega}_{\mathcal{C}}$ admits no core $\mathtt{L}_{\mathcal{C}}$-internalizer.
\end{Corollary}
\begin{proof}
    By Theorem \ref{202508271151}, it suffices to prove that, for any $d_0 \in \mathds{N}$, there exists $m \in \mathds{N}$ such that $\Aug(\mathcal{C}_1^{\times m})^{d_0+1}/\Aug(\mathcal{C}_1^{\times m})^{d_0+2} \not\cong 0$.
    We will show that $m = d_0+1$ is enough.
    First, using Lemma \ref{202510031417}, we obtain
    $$\mathcal{C}_1 / D_{2,p}(\mathcal{C}_1) \cong \mathsf{F}_1/ D_{2,p}(\mathsf{F}_1) \cong \mathds{Z}/p .$$
    Hence, we can choose a nontrivial element $v \in \mathcal{C}_1 / D_{2,p}(\mathcal{C}_1)$.
    Let $v_j \in \mathcal{C}_1^{\times m} / D_{2,p}(\mathcal{C}_1^{\times m}) \subset \mathfrak{L}(\mathcal{C}_1^{\times m}; D_{\bullet,p})$ be the element induced by inserting $v$ into the $j$-th component.
    Since $v_j$'s are linearly independent in the $\mathds{F}_p$-vector space $\mathfrak{L}(\mathcal{C}_1^{\times m}; D_{\bullet,p})$, by the PBW theorem, we obtain 
    $$0 \neq v_1 \cdots v_{m=d_0+1} \in \mathcal{U}^{p}_{d_0+1} ( \mathds{k} \otimes_{\mathds{Z}} \mathfrak{L}(\mathcal{C}_1^{\times m}; D_{\bullet,p})) \cong \Aug(\mathcal{C}_1^{\times m})^{d_0+1}/\Aug(\mathcal{C}_1^{\times m})^{d_0+2}.$$
\end{proof}

\begin{Corollary}
    We assume that the ground ring $\mathds{k}$ is a field.
    If we have $\mathds{k} \otimes_{\mathds{Z}} \mathfrak{Q}_1 (\mathcal{C}_n) \not\cong 0$ for some $n$, then the class $\mathcal{S}^{\omega}_{\mathcal{C}}$ admits no core $\mathtt{L}_{\mathcal{C}}$-internalizer.
\end{Corollary}
\begin{proof}
    Suppose that the class $\mathcal{S}^{\omega}_{\mathcal{C}}$ admits a core $\mathtt{L}_{\mathcal{C}}$-internalizer.
    Let $n \in \mathds{N}$ be as in the statement, and set $M = \mathcal{C}_n$.
    By Theorem \ref{202508271151}, we can choose $d_0 \in \mathds{N}$ such that $\Aug ( M^{\times m})^{d+1} = \Aug ( M^{\times m})^{d+2}$ for $d \geq d_0$ and $m \in \mathds{N}$.
    Since $\mathds{k}$ is assumed to be a field, we can take a $\mathds{k}$-submodule $V_r \subset \Aug ( M)^{r}$ such that $\Aug ( M)^{r} = \Aug (M)^{r+1} \oplus V_r$.
    By the hypothesis, $V_r = 0$ for $r \geq d_0+1$; and $\Aug ( M)^{r} = \bigoplus_{r^\prime \geq r} V_{r^\prime}$.
    Since $\mathds{k}$ is a field, these decompositions give $$\Aug ( M^{\times m})^{d+1} \cong \bigoplus (V_{r_1} \otimes \cdots \otimes V_{r_m})$$ where $r_1 + \cdots + r_m \geq d+1$.
    By considering $\Aug ( M^{\times m})^{d_0+1} = \Aug ( M^{\times m})^{d_0+2}$ for $m=2$, we obtain $V_i \otimes V_j = 0$ where $i+j = d_0+1$.
    In particular, we obtain $V_1 = 0$, which is equivalent to $\Aug (M) = \Aug (M)^{2}$.
    This gives $\Aug(M) /\Aug(M)^{2} \cong 0$, which, by Lemma \ref{202512251348}, contradicts the assumption.
\end{proof}

\begin{Example} 
    Assume that $\mathds{k}$ is a field. 
    For $\mathcal{C} = \mathbf{M}_{R}^{\mathsf{o}}$ with $R$ a unital ring, if $\mathds{k} \otimes_{\mathds{Z}} R \not\cong 0$, then $\mathcal{S}^{\omega}_{\mathcal{C}}$ admits no core $\mathtt{L}_{\mathcal{C}}$-internalizer.
    By Remark \ref{202512251531}, the converse is also true.
\end{Example}

We conclude this section giving the proof of Theorem \ref{202508271151}:
\begin{proof}[Proof of Theorem \ref{202508271151}]
    By Proposition \ref{202405291746}, we know that (1) is equivalent to the condition that the polynomiality ideal filtration $\mathtt{I}^{(\bullet)}_{\mathcal{C}}$ of $\mathtt{L}_{\mathcal{C}}$ left-stabilizes in the sense of Definition \ref{202508270931}.
    We shall prove the condition (2) from this.
    Suppose that $\mathtt{I}^{(\bullet)}_{\mathcal{C}}$ left-stabilizes.
    Let $n \in \mathds{N}$.
    By the hypothesis, there exists $d_0$ such that, for $d \geq d_0$, we have $\mathtt{I}^{(d)}_{\mathcal{C}} (-,n) = \mathtt{I}^{(d+1)}_{\mathcal{C}} (-,n)$.
    In particular, this implies $\mathtt{I}^{(d)}_{\mathcal{C}} (m,n) = \mathtt{I}^{(d+1)}_{\mathcal{C}} (m,n)$.
    Applying Theorem \ref{202408011134}, we obtain $\Aug ( \mathcal{C}_{n}^{\times m})^{d+1} = \Aug ( \mathcal{C}_{n}^{\times m})^{d+2}$ for $d \geq d_0$.
    Note that the choice of $d_0$ does not depend on $m$.
    The converse is also clear from the previous discussion.
\end{proof}

\appendix

\section{Isomorphism-invariant properties and adjunctions}
\label{202510142110}

As mentioned in Section \ref{202510131439}, the adjunctions there arise from a general categorical mechanism, which we record here for use in the sequels; it is not needed for the results of this paper.

Let $\mathcal{B}$ be a category.
We begin by explaining the following identifications:
\begin{align} \label{202509291420}
    \{ \mathrm{Isom.-inv.~properties~for~}\mathcal{B} \} \cong \{ \mathrm{Subclasses~of~} \mathcal{B}/\cong \} \simeq \{ \mathrm{Full~subcategories~of~}\mathcal{B} \} .
\end{align}
An isomorphism-invariant property for the category $\mathcal{B}$ is a property of objects in $\mathcal{B}$, that is preserved by isomorphisms.
Such a property determines a unique subclass of $\mathcal{B}/\cong$ the class of isomorphism classes of objects in $\mathcal{B}$.
On the other hand, a subclass $\mathcal{S} \subset \mathcal{B}/\cong$ induces a full subcategory $\tilde{\mathcal{S}} \subset \mathcal{B}$ consisting of objects in $\mathcal{B}$ whose isomorphism classes lie in $\mathcal{S}$.
Conversely, any full subcategory of $\mathcal{B}$ arises in this way.
For subclasses $\mathcal{S}_1, \mathcal{S}_2$ of $\mathcal{B}$, we have $\mathcal{S}_1 = \mathcal{S}_2$ ($\mathcal{S}_1 \subset \mathcal{S}_2$, resp.) if and only if $\tilde{\mathcal{S}}_1 \simeq \tilde{\mathcal{S}}_2$ ($\tilde{\mathcal{S}}_1 \subset \tilde{\mathcal{S}}_2$, resp.).
These observations lead to the bijections in (\ref{202509291420}).

\begin{Example}
    The property of being of polynomial degree $\leq d$ is an isomorphism-invariant property of left $\mathtt{L}_{\mathcal{C}}$-modules.
    Under the correspondence (\ref{202509291420}), the property corresponds respectively to the full subcategory $\mathtt{L}_{\mathcal{C}}\mbox{-} \mathsf{Mod}^{\leq d}$ of $\mathtt{L}_{\mathcal{C}}\mbox{-}\mathsf{Mod}$ and to the subclass $\mathcal{S}^{d}_{\mathcal{C}}$ of $\mathtt{L}_{\mathcal{C}}\mbox{-}\mathcal{M}\mathsf{od}$, introduced in Definition \ref{202512241512}.
\end{Example}

Let $\mathcal{A}$ and $\mathcal{B}$ be categories, and let
$\mathbb{L} : \mathcal{A} \to \mathcal{B}$ be a functor with right adjoint $\mathbb{R}$.
The unit and counit of the adjunction define isomorphism-invariant properties of objects in $\mathcal{A}$ and $\mathcal{B}$ as follows, which in turn determine full subcategories of $\mathcal{A}$ and $\mathcal{B}$ via (\ref{202509291420}).

\begin{Defn}
    We denote by $\mathcal{B}_{\mathbb{L},\mathbb{R}} \subset \mathcal{B}$ the full subcategory of objects $B \in \mathcal{B}$ such that the counit $\mathbb{L}(\mathbb{R}(B)) \to B$ is an epimorphism.
    Likewise, let $\mathcal{A}_{\mathbb{L},\mathbb{R}} \subset\ \mathcal{A}$ denote the full subcategory of objects $A \in \mathcal{A}$ such that the unit $A \to \mathbb{R}(\mathbb{L}(A))$ is a monomorphism.
\end{Defn}

The adjunction $(\mathbb{L},\mathbb{R})$ restricts to a natural adjunction between
$\mathcal{A}_{\mathbb{L},\mathbb{R}}$ and $\mathcal{B}_{\mathbb{L},\mathbb{R}}$:
\begin{prop}[Refined adjunction] \label{202512241545}
    The restrictions of $\mathbb{L}$ and $\mathbb{R}$ induce the following adjunction:
    $$
    \begin{tikzcd}
            \mathbb{L} : \mathcal{A}_{\mathbb{L},\mathbb{R}} \arrow[r, shift right=1ex, ""{name=G}] & \mathcal{B}_{\mathbb{L},\mathbb{R}} : \mathbb{R} \arrow[l, shift right=1ex, ""{name=F}]
            \arrow[phantom, from=G, to=F, , "\scriptscriptstyle\boldsymbol{\top}"].
    \end{tikzcd}
    $$
    Furthermore, both of these are faithful.
\end{prop}
\begin{proof}
    An elementary argument on adjunctions shows that, for $A \in \mathcal{A}$ and $B = \mathbb{L}(A) \in \mathcal{B}$, the counit $\mathbb{L}(\mathbb{R}(B)) \to B$ is a (split) epimorphism.
    There is an analogous statement for the unit.
    Hence, the restrictions of $\mathbb{L}$ and $\mathbb{R}$ induce an adjunction in the statement.

    We prove that the restriction $\mathbb{R} : \mathcal{B}_{\mathbb{L},\mathbb{R}} \to \mathcal{A}_{\mathbb{L},\mathbb{R}}$ is faithful.
    Let $B,B^\prime \in \mathcal{B}_{\mathbb{L},\mathbb{R}}$.
    The composition with the counit $\mathbb{L}(\mathbb{R}(B)) \to B$ induces the following map:
    \begin{align*}
    \mathcal{B} ( B, B^\prime ) \to \mathcal{B}( \mathbb{L}(\mathbb{R}(B)), B^\prime ) \cong    \mathcal{A} ( \mathbb{R} (B), \mathbb{R} (B^\prime) ) .
    \end{align*}
    This map is injective since the counit at $B \in \mathcal{B}_{\mathbb{L},\mathbb{R}}$ is an epimorphism.
    Analogously, the functor $\mathbb{L}$ is faithful.
\end{proof}

Let $\mathtt{T}$ be a monad in $\mathsf{Mat}_{\mathds{k}}$.
Using the above constructions, we reinterpret the map $\mathcal{V} = \mathcal{V}_{\mathtt{T}}$ of Definition~\ref{202510121351} as follows.
The construction of $\mathcal{V}$ factors as a composite of two maps
$$
\mathcal{V} : ~\{ \mathrm{Two-sided~ideals~of~}\mathtt{T}\} \to \{\mathrm{Adjunctions~involving~}\mathtt{T}\mbox{-}\mathsf{Mod}\} \to \{\mathrm{Subclasses~of~\mathtt{T}\mbox{-}\mathcal{M}\mathsf{od}}\} .
$$
The first map assigns the adjunction in Proposition \ref{202509281759} to a two-sided ideal $\mathtt{J}$.
The second map is induced by the above discussion: it assigns to an adjunction $(\mathbb{L},\mathbb{R})$ the subclass of $\mathtt{T}\mbox{-}\mathcal{M}\mathsf{od}$ given by the isomorphism classes of the category $\mathcal{B}_{\mathbb{L},\mathbb{R}}$.
The image of the map $\mathcal{V}$ consists of closed subclasses (Theorem \ref{202609281358}).

\begin{Example}
    We apply this to $\mathbb{L}=\iota$ and $\mathbb{R}=\mathcal{V}(-;\mathtt{J})$ (Proposition \ref{202509281759}).
    Since $\iota$ is fully faithful, the unit is already an isomorphism, so $\mathcal{A}_{\mathbb{L},\mathbb{R}}=\mathtt{T}/\mathtt{J}\mbox{-}\mathsf{Mod}$; the refined adjunction of
    Proposition \ref{202512241545} therefore only restricts the $\mathtt{T}\mbox{-}\mathsf{Mod}$ side, to $\mathcal{B}_{\mathbb{L},\mathbb{R}}=\mathtt{T}\mbox{-}\mathsf{Mod}^{\mathtt{J}\text{-vg}}$, and (\ref{202512241357}) is exactly the
    resulting equivalence. 
    In this sense, $\mathcal{V}(\mathtt{J})$ corresponds to the isomorphism-invariant property
    ''being $\mathtt{J}$-vanishingly generated" via (\ref{202509291420}).
\end{Example}

\section*{Acknowledgements}

The author is supported by a KIAS Individual Grant MG093701 at Korea Institute for Advanced Study.
He would like to express his sincere gratitude to Teimuraz Pirashvili, Geoffrey Powell, Christine Vespa and Aurélien Djament for carefully reading the earlier versions of this series and providing him constructive comments.
He also thanks the anonymous referee, whose anticipation of the dichotomy for $\mathrm{D}(\mathcal{C})$ (Theorem \ref{2026009152033}) prompted the further generalization of the paper's other results.

\bibliography{reference}

@article{djament2016cohomologie,
  title={Cohomologie des foncteurs polynomiaux sur les groupes libres},
  author={Djament, Aur{\'e}lien and Pirashvili, Teimuraz and Vespa, Christine},
  journal={Documenta Mathematica},
  volume={21},
  pages={205--222},
  year={2016}
}

@article {kim2024analytic,
    AUTHOR = {Kim, Minkyu and Vespa, Christine},
     TITLE = {On analytic exponential functors on free groups},
   JOURNAL = {Journal of Homotopy and Related Structures},
    VOLUME = {},
      YEAR = {2025},
     PAGES = {},
      ISSN = {1512-2891},
   MRCLASS = {},
  MRNUMBER = {},
MRREVIEWER = {},
       DOI = {10.1007/s40062-025-00390-9},
       URL = {https://doi.org/10.1007/s40062-025-00390-9},
        note = {Published online: 2025-12-03}
}

@article {EML,
    AUTHOR = {Eilenberg, Samuel and Mac Lane, Saunders},
     TITLE = {On the groups {$H(\Pi,n)$}. {II}. {M}ethods of computation},
   JOURNAL = {Ann. of Math. (2)},
  FJOURNAL = {Annals of Mathematics. Second Series},
    VOLUME = {60},
      YEAR = {1954},
     PAGES = {49--139},
      ISSN = {0003-486X},
   MRCLASS = {56.0X},
  MRNUMBER = {65162},
MRREVIEWER = {P. J. Hilton},
       DOI = {10.2307/1969702},
       URL = {https://doi.org/10.2307/1969702},
}

@incollection{reutenauer2003free,
  title={Free {L}ie algebras},
  author={Reutenauer, Christophe},
  booktitle={Handbook of algebra},
  volume={3},
  pages={887--903},
  year={2003},
  publisher={Elsevier}
}

@article{betti1983variation,
  title={Variation through enrichment},
  author={Betti, Renato and Carboni, Aurelio and Street, Ross and Walters, Robert},
  journal={Journal of Pure and Applied Algebra},
  volume={29},
  pages={109--127},
  year={1983}
}

@book{franjou2003rational,
  title={Rational representations, the {S}teenrod algebra and functor homology},
  author={Franjou, Vincent and Friedlander, Eric M and Pirashvili, Teimuraz and Schwartz, Lionel},
  volume={16},
  year={2003},
  publisher={Soci{\'e}t{\'e} math{\'e}matique de France Paris}
}

@incollection {benabou2006introduction,
    AUTHOR = {B\'enabou, Jean},
     TITLE = {Introduction to bicategories},
 BOOKTITLE = {Reports of the {M}idwest {C}ategory {S}eminar},
    SERIES = {Lecture Notes in Math.},
    VOLUME = {No. 47},
     PAGES = {1--77},
 PUBLISHER = {Springer, Berlin-New York},
      YEAR = {1967},
   MRCLASS = {18.10},
  MRNUMBER = {220789},
MRREVIEWER = {J.\ R.\ Isbell},
}

@article{hartl2015polynomial,
  title={Polynomial functors from algebras over a set-operad and nonlinear {M}ackey functors},
  author={Hartl, Manfred and Pirashvili, Teimuraz and Vespa, Christine},
  journal={International Mathematics Research Notices},
  volume={2015},
  number={6},
  pages={1461--1554},
  year={2015},
  publisher={Oxford University Press}
}

@book {MR537126,
    AUTHOR = {Passi, Inder Bir S.},
     TITLE = {Group rings and their augmentation ideals},
    SERIES = {Lecture Notes in Mathematics},
    VOLUME = {715},
 PUBLISHER = {Springer, Berlin},
      YEAR = {1979},
     PAGES = {vi+137},
      ISBN = {3-540-09254-4},
   MRCLASS = {20C07},
  MRNUMBER = {537126},
MRREVIEWER = {E.\ Formanek},
}

@inproceedings{sandling1979augmentation,
  title={Augmentation quotients of group rings and symmetric powers},
  author={Sandling, Robert and Tahara, Ken-Ichi},
  booktitle={Mathematical Proceedings of the Cambridge Philosophical Society},
  volume={85},
  number={2},
  pages={247--252},
  year={1979},
  organization={Cambridge University Press}
}

@article{casacuberta1999singly,
  title={Singly generated radicals associated with varieties of groups},
  author={Casacuberta, Carles and Rodr{\'\i}guez, Jos{\'e} L and Scevenels, Dirk},
  journal={London Mathematical Society Lecture Note Series},
  pages={202--210},
  year={1999},
  publisher={Cambridge University Press}
}

@book {PHall1976,
    AUTHOR = {Hall, Jr., Marshall},
     TITLE = {The theory of groups},
      NOTE = {Reprinting of the 1968 edition},
 PUBLISHER = {Chelsea Publishing Co., New York},
      YEAR = {1976},
     PAGES = {xiii+434},
   MRCLASS = {20-02},
  MRNUMBER = {414669},
}

@article {PolyPira,
    AUTHOR = {Pirashvili, T. I.},
     TITLE = {Polynomial functors},
   JOURNAL = {Trudy Tbiliss. Mat. Inst. Razmadze Akad. Nauk Gruzin. SSR},
  FJOURNAL = {Akademiya Nauk Gruzinsko\u i\ SSR. Trudy Tbilisskogo
              Matematicheskogo Instituta im. A. M. Razmadze},
    VOLUME = {91},
      YEAR = {1988},
     PAGES = {55--66},
      ISSN = {0234-5838},
   MRCLASS = {18G50 (16E20)},
  MRNUMBER = {1029007},
MRREVIEWER = {Luchezar\ L.\ Avramov},
}

@article {Mitchell1972,
    AUTHOR = {Mitchell, Barry},
     TITLE = {Rings with several objects},
   JOURNAL = {Advances in Math.},
  FJOURNAL = {Advances in Mathematics},
    VOLUME = {8},
      YEAR = {1972},
     PAGES = {1--161},
      ISSN = {0001-8708},
   MRCLASS = {18H15},
  MRNUMBER = {294454},
MRREVIEWER = {J.-P.\ Lafon},
       DOI = {10.1016/0001-8708(72)90002-3},
       URL = {https://doi.org/10.1016/0001-8708(72)90002-3},
}

@article {Piras1993,
    AUTHOR = {Pirashvili, Teimuraz},
     TITLE = {Polynomial approximation of {${\rm Ext}$} and {${\rm Tor}$}
              groups in functor categories},
   JOURNAL = {Comm. Algebra},
  FJOURNAL = {Communications in Algebra},
    VOLUME = {21},
      YEAR = {1993},
    NUMBER = {5},
     PAGES = {1705--1719},
      ISSN = {0092-7872,1532-4125},
   MRCLASS = {18G15 (16E40)},
  MRNUMBER = {1213983},
MRREVIEWER = {Antonio\ M.\ Cegarra},
       DOI = {10.1080/00927879308824647},
       URL = {https://doi.org/10.1080/00927879308824647},
}

@article {Lawvere1963,
    AUTHOR = {Lawvere, F. William},
     TITLE = {Functorial semantics of algebraic theories},
   JOURNAL = {Proc. Nat. Acad. Sci. U.S.A.},
  FJOURNAL = {Proceedings of the National Academy of Sciences of the United
              States of America},
    VOLUME = {50},
      YEAR = {1963},
     PAGES = {869--872},
      ISSN = {0027-8424},
   MRCLASS = {18.10},
  MRNUMBER = {158921},
MRREVIEWER = {M.\ Artin},
       DOI = {10.1073/pnas.50.5.869},
       URL = {https://doi.org/10.1073/pnas.50.5.869},
}

@article{massuyeau2007finite,
  title={Finite-type invariants of 3-manifolds and the dimension subgroup problem},
  author={Massuyeau, Gw{\'e}na{\"e}l},
  journal={Journal of the London Mathematical Society},
  volume={75},
  number={3},
  pages={791--811},
  year={2007},
  publisher={Wiley Online Library}
}

@article {Quillen1968,
    AUTHOR = {Quillen, Daniel G.},
     TITLE = {On the associated graded ring of a group ring},
   JOURNAL = {J. Algebra},
  FJOURNAL = {Journal of Algebra},
    VOLUME = {10},
      YEAR = {1968},
     PAGES = {411--418},
      ISSN = {0021-8693},
   MRCLASS = {20.56 (17.00)},
  MRNUMBER = {231919},
MRREVIEWER = {J.\ Knopfmacher},
       DOI = {10.1016/0021-8693(68)90069-0},
       URL = {https://doi.org/10.1016/0021-8693(68)90069-0},
}

@article {HV2011,
    AUTHOR = {Hartl, Manfred and Vespa, Christine},
     TITLE = {Quadratic functors on pointed categories},
   JOURNAL = {Adv. Math.},
  FJOURNAL = {Advances in Mathematics},
    VOLUME = {226},
      YEAR = {2011},
    NUMBER = {5},
     PAGES = {3927--4010},
      ISSN = {0001-8708,1090-2082},
   MRCLASS = {18A25},
  MRNUMBER = {2770438},
MRREVIEWER = {Timothy\ Porter},
       DOI = {10.1016/j.aim.2010.11.008},
       URL = {https://doi.org/10.1016/j.aim.2010.11.008},
}

@article {BP1999,
    AUTHOR = {Baues, Hans-Joachim and Pirashvili, Teimuraz},
     TITLE = {Quadratic endofunctors of the category of groups},
   JOURNAL = {Adv. Math.},
  FJOURNAL = {Advances in Mathematics},
    VOLUME = {141},
      YEAR = {1999},
    NUMBER = {1},
     PAGES = {167--206},
      ISSN = {0001-8708,1090-2082},
   MRCLASS = {20J15 (18D10)},
  MRNUMBER = {1667150},
MRREVIEWER = {Timothy\ Porter},
       DOI = {10.1006/aima.1998.1784},
       URL = {https://doi.org/10.1006/aima.1998.1784},
}

@article {BJ1994,
    AUTHOR = {Baues, Hans Joachim},
     TITLE = {Quadratic functors and metastable homotopy},
   JOURNAL = {J. Pure Appl. Algebra},
  FJOURNAL = {Journal of Pure and Applied Algebra},
    VOLUME = {91},
      YEAR = {1994},
    NUMBER = {1-3},
     PAGES = {49--107},
      ISSN = {0022-4049,1873-1376},
   MRCLASS = {55U35 (55Q40)},
  MRNUMBER = {1255923},
MRREVIEWER = {Donald\ M.\ Davis},
       DOI = {10.1016/0022-4049(94)90135-X},
       URL = {https://doi.org/10.1016/0022-4049(94)90135-X},
}

@article{DTV2023,
	author = {Djament, Aurélien and Touzé, Antoine and Vespa, Christine},
	journal = {Annales Scientifiques de l'École Normale Supérieure},
	title = {Décompositions à la {S}teinberg sur une catégorie additive},
	doi = {10.24033/asens.2538},
	publisher = {Société mathématique de France},
	year = {2023},
	volume = {56},
	pages = {427-516}
}

@article{kimprop,
title = {{PROP}s associated to {L}awvere theories and their relation to polynomial functors},
journal = {Journal of Algebra},
volume = {709},
pages = {467-563},
year = {2027},
issn = {0021-8693},
doi = {https://doi.org/10.1016/j.jalgebra.2026.06.043},
url = {https://www.sciencedirect.com/science/article/pii/S0021869326004035},
author = {Minkyu Kim}
}
\bibliographystyle{alpha}
\end{document}